\documentclass{article}
\usepackage{graphicx}
\usepackage{amsmath}
\usepackage{amssymb}
\usepackage{xcolor}
\usepackage{float,algorithm}
\usepackage{algpseudocode}
\usepackage{tikz,hyperref}
\usepackage{stmaryrd}
\usetikzlibrary{calc}
\usepackage{comment}
\usepackage{enumitem}
\usepackage{pdflscape}
\usepackage{placeins}
\usepackage{authblk}

\usepackage{siunitx}
\usepackage{adjustbox}

\usepackage{bm}
\newtheorem{problem}{Problem}

\newcommand{\sprodf}[3][{\Omega}_f^0]{( {#2}, {#3} )_{#1}}
\newcommand{\sprods}[3][{\Omega}_s^0]{( {#2}, {#3} )_{#1}}

\newcommand{\bmh}[1]{\widehat{\bm{{#1}}}}

\newcommand{\bs}[1]{\boldsymbol{#1}}
\newcommand{\bx}{{\bs x}}

\renewcommand{\vec}[1]{{\bf #1}}
\renewcommand{\div}{{\rm div}}
\newcommand{\Div}{{\rm \bf div}}
\newcommand{\tr}{{\rm tr}}
\newcommand{\pder}[2]{\frac{\partial #1}{\partial #2}}
\def\eqd{\stackrel{\mathrm{def}}{=}}
\newcommand{\jump}[1]{\bigl[\hspace{-0.03in}\bigl[#1\bigr]\hspace{-0.03in}\bigr]}

\definecolor{darkblue}{RGB}{0,0,139}

\definecolor{gold}{RGB}{212,175,55}

\title{{A numerical benchmark for fluid--structure--contact interaction}
}

\author[1]{D.C.~Corti
}
\author[2]{J.~Fara
}
\author[1]{M.A.~Fernández
}
\author[3]{S.~Frei\thanks{Corresponding author: stefan.frei@uni-konstanz.de}}
\author[4]{T.~Knoke
}
\author[5,6]{S.~Schwarzacher
}
\author[2]{K.~T\r{u}ma
}
\author[4]{T.~Wick
}

\affil[1]{Sorbonne Universit\'e, CNRS, Laboratoire Jacques-Louis Lions, Inria, Paris, France}
\affil[2]{Mathematical Institute, Faculty of Mathematics and Physics, Charles University, Prague, Czech Republic}
\affil[3]{University of Konstanz, Department
 of Mathematics \& Statistics, Konstanz, Germany}
\affil[4]{Leibniz University Hannover, Institute of Applied Mathematics, Hannover, Germany}
\affil[5]{Department of Mathematical Analysis, Faculty of Mathematics and Physics,
Charles University, Prague, Czech Republic}
\affil[6]{Department of Mathematics, Uppsala University,
Uppsala, Sweden}

\begin{document}
\maketitle

\begin{abstract}
We propose a two-dimensional benchmark for fluid--structure--contact interaction  consisting of a deformable elastic disk falling under gravity within a viscous incompressible fluid and rebounding in the vicinity of the bottom wall. Solid deformability is essential, as rigid solids do not rebound in this framework. Besides this,
the setting is deliberately kept simple to facilitate reproduction. The configuration is particularly challenging due to the well-known no-contact paradox, which can lead to a contactless rebound and forces numerical methods to resolve a vanishingly thin fluid layer in the near-contact region, making the dynamics highly sensitive to the spatial and temporal discretizations. In addition to no-slip boundary and interface conditions, a reduced porous modeling of surface roughness, either on the disk boundary or on the bottom wall, is also considered; this circumvents the no-contact paradox and enables genuine contact. An energy balance law is derived theoretically for all three cases. Eight numerical methodologies, developed by five research groups and spanning different model formulations, numerical methods, and codes (including both fitted and unfitted discretizations), are applied to the benchmark at several levels of spatial and temporal refinement.
 Quantities of interest of varying complexity are collected and compared, showing close agreement during the falling phase and increased sensitivity in the near-contact and rebound regimes.
 The setting and the results provide a suitable reference for the systematic assessment of fluid--structure--contact interaction solvers.
 The time histories of all quantities of interest for every approach and refinement level are provided as supplementary material.
\end{abstract}


\section{Introduction}
When an elastic solid approaches another one in a viscous incompressible fluid, the fluid
trapped in the shrinking gap has to be squeezed out. The lubrication forces
this generates grow without bound as the gap closes, and they govern the
approach, the impact and the subsequent rebound. Solid contact in a fluid is
therefore a singular event in a precise sense: what decides whether, and how, a solid rebounds is set within a layer of fluid orders of magnitude thinner than the solid itself. Capturing this interplay has motivated a sustained research effort in physics, engineering and
mathematics (see, e.g., \cite{Brenner1961,CoxBrenner1967,Hesla,Hil07,MR2481302,TezduyarSathe2007,DosSantosEtAl2008,AstorinoGerbeauetal2009,MayerWallEtAl2013,Kamenskyetal2015,
Pironneau2016,FreiRichter2017,FormaggiaGattiZonca2021,
AgerSchottVuongPoppWall2019,BurmanFernandezFrei2020,vonWahletal2021,
GraSchSouTum22,BurmanFernandezFreiGerosa2022,ChampionFernandezetal2024,FarSchTum24,GerosaMarsden2024,champion:hal-05626925}),
and it is what makes contact and rebound so demanding for numerical simulation.
Yet a faithful reproduction of these phenomena is required in a wide range of
applications, including the opening and closure of native or
prosthetic heart valves~\cite{AstorinoGerbeauetal2009,Kamenskyetal2015,
TeraharaEtAl2020}, mechanical and microfluidic
valves~\cite{YuEtAl2017,WangEtAl2021,AissaBerraiesEtAl2025}, human
phonation~\cite{HoracekEtAl2005,ZhengEtAl2011,SvacekHoracek2021}, lubricated
contacts in bearings and related tribological systems~\cite{Knaufetal,
Bruyere2012} and tire-water-road
interaction~\cite{HermangeEtAl2019,NazariEtAl2020}.

Despite the growing interest in fluid-structure interaction (FSI) with contact over the past decade (see, e.g., 
~\cite{Richter2017,FreiHolmRichterWickYang2018,AgerSchottVuongPoppWall2019,ager-et-al-21,FormaggiaGattiZonca2021, BurmanFernandezFreiGerosa2022,
GraSchSouTum22,FarSchTum24,ChampionFernandezetal2024,GerosaMarsden2024,FreiKnokeSteinbachWenskeWick2026,champion:hal-05626925}), reference data remain scarce and the existing numerical approaches are rarely compared on a common setting. To fill this gap, we propose a benchmark for assessing the ability of a simulation package to reproduce contact and rebound behavior.
Inspired by the recent work reported in~\cite{vonWahletal2021, HagemeierTheveninRichter2021}, the setting consists of a deformable elastic disk falling freely in a viscous incompressible fluid, which is expected to rebound after contact with the bottom wall.
Deformability is critical here, since rigid solids do not rebound in fluids~\cite{GraSchSouTum22}. Given the singular nature of the problem, we keep the setting two-dimensional and of minimal complexity, in order to facilitate reproduction.

The proposed benchmark follows the tradition of previous prominent FSI test cases, most notably the Turek--Hron benchmark~\cite{HronTurek2006benchmark} (see also~\cite{BungartzSchaefer2006,BuSc10}), which at the time of its launch was a challenging computational task and has inspired substantial development in the community in the two decades since. Experimental reference data for a related laminar configuration were provided by Gomes et al.~\cite{GOMES201143}. A three-dimensional benchmark, including experimental results, was proposed more recently in the framework of an international call for participation by Hessenthaler et al.~\cite{Hessenthaleretal2017} (see also \cite{hessenthaler-et-al-17,jansson-et-al-17,landajuela-et-al-17}).

FSI with contact remains a delicate problem, and the accuracy and robustness of the available methods in such a singular setting are largely unexplored: the outcome may depend on the mathematical model and its numerical approximation alike. Approaches using unfitted meshes and Eulerian descriptions of the fluid are often favored for this kind of problem, due to the potential topology changes in the fluid domain,  but a thorough systematic numerical
comparison is, as far as we know, still missing. This also justifies the involvement of several research groups, since it is unusual for a single group to master all available approaches. Indeed, in our view, the term {\it benchmark} is only warranted when several groups participate in the comparison; see again~\cite{HronTurek2006benchmark} or, for example,~\cite{SchaeTu96,FLEMISCH2018239}.

What distinguishes this benchmark from earlier FSI test cases is the near-contact and contact regimes.
As one might expect, the results agree far better before near-contact than after it, where very small changes, on a range of scales, influence the subsequent dynamics; we give further insight into this below.
The fluid is described by the viscous incompressible Navier--Stokes equations, which, although debatable for some applications, is the
standard assumption. If, moreover, no-slip boundary conditions are enforced (likewise the most common choice), the so-called Cox--Brenner paradox~\cite{Brenner1961,CoxBrenner1967,cooley_o-neill-1969} arises: smooth rigid objects in an incompressible fluid can then never come into contact in finite time. Brenner showed this using lubrication theory~\cite{Brenner1961}, and Hesla and Hillairet later provided mathematical proofs of this no-contact  paradox~\cite{Hesla,Hil07} (see also~\cite{MR2481302} for extensions to the three-dimensional setting).

Numerical evidence indicates that the no-contact paradox also persists for deformable solids, leading to the theoretical possibility of a contactless rebound (see~\cite{Frei2016,Richter2013,GraSchSouTum22,FarSchTum24}). No experiment demonstrating a contactless rebound is available to date, and it is far from clear whether one is possible. 
Nevertheless, the no-slip model is expected to produce such a contactless rebound, no matter how thin the minimal fluid gap becomes.
How to treat the near-contact regime is thus left to the numerical method: resolving the thin fluid layer via mesh refinement, using ad hoc rebound laws, or modifying the assumptions in the vicinity of contact so as to avoid the no-contact paradox.

Here, we consider several modeling options and numerical strategies, and compare how they approximate the rebound. The computations fall into two classes, each represented by independent codes. The first concerns the contactless rebound, which several strategies recover under standard no-slip boundary conditions and which appears to be the behavior of the underlying continuous coupled problem.
The second relies on additional models of the fluid layer between contacting solids, developed in recent years to circumvent the no-contact paradox and motivated, {for} example, by surface roughness at the micro scale. Such models retain a thin fluid layer even when parts of the solid are in contact with the ground; this layer can be described by porous-media models (see~\cite{AgerSchottVuongPoppWall2019, BurmanFernandezFrei2020, BurmanFernandezFreiGerosa2022,GerosaMarsden2024}). In~\cite{ChampionFernandezetal2024} it was proved that the no-contact paradox is indeed circumvented, at least for rigid bodies, when the porous layer model of~\cite{BurmanFernandezFreiGerosa2022} replaces the no-slip condition on the ground. {Alongside the benchmark proposal with no-slip conditions, we therefore also investigate the effect of a porous layer on the effective behavior during and after contact, by imposing a Darcy law on either the bottom boundary or the interface.}

Five groups, each with its own code, have contributed results for this benchmark, obtained with eight different numerical methodologies in total, and compared them. A range of quantities of interests (QoIs), such as the motion of the {disk}, the contact time, the energy loss due to dissipation, and {others}, has been collected by each group at several levels of spatial and temporal refinement. The agreement between the results allows us to propose the setting, together with the collected QoIs, as a benchmark. Following the practice established by the Turek--Hron benchmark~\cite{HronTurek2006benchmark}, we have kept it rich enough to extend existing FSI benchmark knowledge, yet simple enough for many groups to participate and for a critical mass of comparable data to be collected. The QoIs are of varying complexity, so that later studies can select a subset of them: not every code can readily compute point values, line integrals and domain integrals. Finally, to make the comparison directly reusable, the time histories of all dynamic quantities are provided for each approach and refinement level as supplementary material~\cite{BenchmarkData}.

The rest of this paper is organized as follows. In Section~\ref{sec:setting}, we present the governing equations in three equivalent formulations, and state
an energy balance for the resulting coupled problems. The detailed
description of the benchmark, including the geometry, the material
parameters and the quantities of interest, is provided in Section~\ref{sec:benchmark:description}. Section~\ref{results} is devoted to the numerical results: we first list the
participating teams together with a short summary of their numerical approaches, and
then present and discuss the resulting QoIs in detail. Section~\ref{sec:conclusions} summarizes our findings. Finally, a sketch of the proof of the energy balance is reported in Part~\ref{app:proof} of
the appendix, while Part~\ref{sec:appendix}
provides details of the numerical approaches and codes used by the
different groups.

\section{Problem setting}
\label{sec:setting}

\begin{figure}[h!]
    \centering
    \begin{tikzpicture}[scale=2,>=stealth]

    \draw[thick] (0,0) rectangle (2,2);

    \node[above] at (1,1.7) {$\Gamma_{\mathrm{top}}$};
    \node[below] at (1,0.3) {$\Gamma_{\mathrm{bot}}$};
    \node[right, rotate=90] at (1.85,0.8) {$\Gamma_{\mathrm{wall}}$};

    \draw[<->] (0,2.2) -- (2,2.2)
        node[midway, above] {$2R$ };
    \draw[<->] (2.2,0) -- (2.2,2)
        node[midway, right] {$H$};

    \coordinate (C) at (1.0,1.3);
    \def\r{0.25}
    \draw[thick,fill=gray!60] (C) circle (\r);
    \node at (C) {+};

    \draw[->] (C) ++(25:0.35) node[right] {$c_{\rm s}$} -- ++(205:0.2);

    \node[left] at ($(C)+(-\r,0)$) {$\Gamma_{0}^{\rm i}$};

    \draw[<->] ($(C)+(-\r,-0.35)$) -- ++(2*\r,0)
        node[midway,below] {$d_{\rm s}$};

    \end{tikzpicture}
    \caption{Geometrical description of the benchmark setting.}
    \label{fig:domain}
\end{figure}
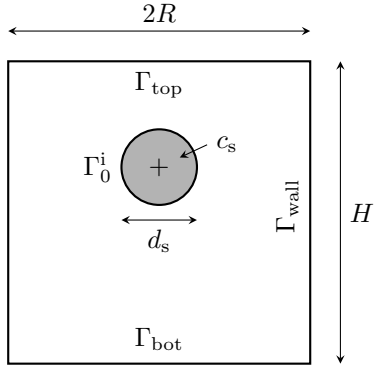

We investigate the simple setting visualized in Figure~\ref{fig:domain}. A deformable elastic disk of diameter $d_{\rm s}$ is placed within a viscous fluid in a rectangular box $\Omega\subset\mathbb{R}^2$. Under gravity $\vec{g}$, the disk falls until it approaches the ground, represented by the lower boundary $\Gamma_{\text{bo{t}}}$ of $\Omega$. We will investigate the near-contact and contact dynamics as well as the subsequent rebound.
A similar setting was proposed in~\cite{vonWahletal2021} in a rotationally symmetric three-dimensional configuration, with numerical results compared against experimental data~\cite{HagemeierTheveninRichter2021}. 

The full domain $\Omega$ is split into a time-dependent fluid sub-domain $\Omega^{\rm f}_t$ and a time-dependent solid sub-domain $\Omega^{\rm s}_t$ such that  
\begin{equation}\label{eq:omega-full}
\overline{\Omega}=\overline{\Omega^{\rm f}_t\cup\Omega^{\rm s}_t}.
\end{equation}

\subsection{Coupled problem in mixed Lagrangian-Eulerian formalism}
\label{sec.lageul}

The solid is assumed to be elastic, compressible and homogeneous, with reference configuration $\Omega^{\rm s}_0 \subset \mathbb{R}^2 $.
In a total Lagrangian formulation, the dynamics of the solid are described in terms of the motion map 
$\vec{X}: \Omega^{\rm s}_0 \times \mathbb{R}^+ \rightarrow \mathbb{R}^2$, so that $\vec{X}(\vec{x},t)$ provides the position at time $t \geq 0$ of the material point $\vec{x}\in \Omega^{\rm s}_0$. 
We can also introduce the instantaneous deformation $\vec{X}_t \eqd \vec{X}(\cdot,t)$ and the solid displacement $\vec{u}^{\rm s}: \Omega^{\rm s}_0 \times \mathbb{R}^+ \rightarrow \mathbb{R}^2$ as 
$\vec{u}^{\rm s}(\vec{x},t) \eqd \vec{X}(\vec{x},t) - \vec{x}$.
The dynamics of the solid are  governed by the linear momentum conservation equation 
\begin{equation}\label{eq:solid}
  \begin{split}
    \varrho^{\rm s} \pder{^2\vec{u}^{\rm s}}{t^2} - \Div \mathbb{P} = \varrho^{\rm s} \vec{g} 
  \end{split} \quad \text{in}\quad \Omega^{\rm s}_0,
\end{equation}
where $\vec{g}$ denotes the gravitational acceleration vector and $\varrho^{\rm s} >0$ is the solid density in the reference configuration, assumed to be constant, $\mathbb{P}$ is the first Piola--Kirchhoff stress tensor related to the second Piola--Kirchhoff stress tensor $\mathbb{S}$ through
$$
    \mathbb{P}=\mathbb{F}\,\mathbb{S},
$$
where $\mathbb{F} \eqd {\boldsymbol \nabla}\vec{X} =\mathbb{I}+{\boldsymbol \nabla}\vec{u}^{\rm s}$ is the deformation gradient. 
We also set $J \eqd \det \mathbb{F} > 0$ for the Jacobian. 
We assume that the material is hyperelastic, i.e., 
\begin{equation}\label{eq:hyper}
  \mathbb{S} = \frac{\partial W}{\partial \mathbb{E}}, 
\end{equation}
with $\mathbb{E} \eqd \frac12 (\mathbb{F}^{\rm T}\mathbb{F} - \mathbb{I})$  the Green-Lagrange strain tensor and $W:\mathbb{R}^{2\times 2} \to \mathbb{R}$ the strain energy density function. 
In this paper, we consider an isotropic geometrically non-linear Saint Venant--Kirchhoff law, so that 
$$
W\big(\mathbb{E} \big)=\frac{\lambda^{\rm s}}{2}(\tr\,\mathbb{E})^2+\mu^{\rm s}\left|\mathbb{E}\right|^2, 
\quad 
  \mathbb{S}  = \lambda^{\rm s} (\tr\, \mathbb{E})\, \mathbb{I} + 2 \mu^{\rm s} \mathbb{E},
$$
where    $\lambda^{\rm s},\mu^{\rm s}>0$ stand for Lamé constants  characterizing the material.

The current configuration of the solid, $\Omega^{\rm s}_t$, is given by the relation 
$
 \Omega^{\rm s}_t \eqd  \vec{X}_t(\Omega^{\rm s}_0).
$
Owing to \eqref{eq:omega-full}, the fluid domain is then defined as 
\begin{equation}\label{eq:geometry-coupling}
 \Omega^{\rm f}_t \eqd  \Omega  \backslash   \overline{\Omega^{\rm s}_t}
\end{equation}
and its associated space-time trajectory is given by
$$
 \mathcal S \eqd  \bigcup_{t \in \mathbb{R}^+} \Omega^{\rm f}_t \times \{t\}.
$$
The fluid filling the region $\Omega^{\rm f}_t$ is assumed to 
be Newtonian, incompressible and homogeneous. 
Its state is described by the unknown velocity and
pressure fields $\vec{v}^{\rm f}: \mathcal S\rightarrow \mathbb{R}^2$ and $p^{\rm f}: \mathcal S \rightarrow \mathbb{R}$. In an Eulerian formalism, the balance equations for the fluid are given by the following Navier--Stokes system:  
\begin{equation}\label{eq:fluid}
\left\{ 
\begin{aligned}
  \varrho^{\rm f} \left( \pder{\vec{v}^{\rm f}}{t}+ \vec{v}^{\rm f} \cdot {\boldsymbol \nabla}\vec{v}^{\rm f} \right) -
   \Div\,\mathbb{T}^{\rm f} = \varrho^{\rm f}\vec{g}& \quad \text{in}\quad \Omega^{\rm f}_t,\\
  \div \vec{v}^{\rm f} = 0& \quad \text{in}\quad \Omega^{\rm f}_t,
\end{aligned}
\right.
\end{equation}
where $\rho^{\rm f}>0$ denotes fluid density  and the fluid Cauchy stress tensor $\mathbb{T}^{\rm f}$ is given by
$$
  \mathbb{T}^{\rm f}=-p^{\rm f} \mathbb{I} + \mu^{\rm f} \left({\boldsymbol \nabla}\vec{v}^{\rm f} +
    ({\boldsymbol \nabla}\vec{v}^{\rm f})^{\top}\right),
$$
with $\mu^{\rm f}>0$ standing  for the dynamic viscosity of the fluid. 
We denote by $\Gamma_t^{\rm i} \eqd  \partial \Omega_t^{\rm f} \cap \partial \Omega_t^{\rm s} = \partial \Omega_t^{\rm s}$ 
the current configuration of the fluid-structure interface. In particular, we have $\Gamma_t^{\rm i} = \vec{X}_t(\Gamma^{\rm i}_0)$ with 
$\Gamma^{\rm i}_0 \eqd \partial \Omega^{\rm s}_0$. 
Besides the geometrical coupling \eqref{eq:geometry-coupling}, the mechanical interaction between the fluid and the solid 
is {imposed} on the fluid-solid interface via the following 
kinematic and dynamic coupling conditions: 
\begin{equation}\label{eq:fsi-coupling}
\left\{
\begin{aligned}
    \vec{v}^{\rm f}= \pder{\vec{u}^{\rm s}}{t}\circ \vec{X}_t^{-1}  \quad \text{on} \quad \Gamma^{\rm i}_t ,&\\
  \mathbb{P}\vec{n}_0^{\rm s}  = J 
    \mathbb{T}^{\rm f}  \circ \vec{X}_t  \mathbb{F}^{-\top} \vec{n}_0^{\rm s}  
  \quad \text{on} \quad \Gamma^{\rm i}_0,&
  \end{aligned}\right. 
\end{equation} 
where $\vec{n}_0^{\rm s}$ 
is the solid outward unit normal to $\Gamma_0^{\rm i}$.

The coupled fluid-structure interaction problem in mixed Lagrangian-Eulerian formalism \eqref{eq:solid}-\eqref{eq:fsi-coupling} needs to be completed with appropriate boundary and initial  conditions.  Here, we assume that the system is  initially at rest, so that $\vec{v}^{\rm f}(\cdot,0) = \vec{0} $, $\vec{u}^{\rm s}(\cdot,0) = \vec{0} $ and $\vec{v}^{\rm s}(\cdot,0) = \vec{0}$.
At the bottom $\Gamma_{\rm bot}$ and the walls $\Gamma_{\rm wall}$ we assume no-slip boundary condition; at the top $\Gamma_{\rm top}$ a free slip boundary condition is imposed, i.e.,
\begin{equation}\label{eq:bc-fluid}
\left\{ 
\begin{aligned}
    \vec{v}^{\rm f}=\vec{0}& \quad \text{on}\quad \Gamma_{\text{bo{t}}}\cup\Gamma_{\text{wa{ll}}},\\
    \vec{v}^{\rm f} \cdot \vec{n}  =0,\quad \bs {\tau}^\top \mathbb{T}^{\rm f}\vec{n} =0&\quad   \text{on} \quad  \Gamma_{\text{to{p}}},
\end{aligned}
\right.
\end{equation}
where $\vec{n}$ and $\boldsymbol{\tau}$ respectively denote the outer unit normal and tangent to $\partial \Omega_t^{\rm f} = \Gamma_t^{\rm i}\cup  \Gamma_{\text{bo{t}}}\cup\Gamma_{\text{wa{ll}}} \cup \Gamma_{\text{to{p}}}$. 





Alternatively, a given set of equations can be equivalently formulated in a different coordinate system, thereby leading to different numerical approaches. In the following  sections, we describe the Arbitrary Lagrangian--Eulerian approach, in which the fluid equations are formulated in an arbitrary reference frame, and a fully Eulerian approach, in which the solid equations are formulated in the  current configuration.



\subsection{Coupled problem with Arbitrary Lagrangian-Eulerian formalism in the fluid}
\label{subsec.ALE}

In the Arbitrary Lagrangian--Eulerian (ALE) approach (see, e.g., \cite{HughesLiuZimmermann1981,DoneaGiulianiHalleux1982,TurekHron2006,HronTurek2006benchmark,RichterWick2010,FreiRichterWick2016,landajuela-et-al-17,Richter2017,FarSchTum24}), 
the fluid  domain $\Omega_t^{\rm f}$ is assumed to be parametrized by a 
smooth injective map 
$  
\mathcal{A}: \Omega^{\rm f}_0 \times \mathbb{R}^+ \to \mathbb{R}^2$
as $\Omega_t^{\rm f} = \mathcal{A}(\Omega_0^{\rm f},t)$,
where $\Omega_0^{\rm f}$ denotes 
a given reference fluid domain.
The ALE map $\mathcal{A}$ is often 
given in terms of the fluid domain displacement $\vec{u}^{\rm f}: \Omega^{\rm f}_0 \times \mathbb{R}^+  \to \mathbb{R}^2$, 
as $\mathcal{A}(\vec{x},t) \eqd  \vec{x} + \vec{u}^{\rm f}(\vec{x},t)$ for  $\vec{x}\in \Omega_0^{\rm f}$, 
The corresponding fluid domain velocity and ALE deformation gradient and Jacobian are respectively defined as 
$$
\vec{w}^{\rm f} \eqd  \frac{ \partial \vec{u}^{\rm f} }{ \partial t},\quad 
    \mathbb{F}^{\rm f} \eqd  \mathbb{I} + \boldsymbol \nabla \vec{u}^{\rm f}, \quad J^{\rm f} \eqd  \det \mathbb{F}^{\rm f} > 0.
$$
Using the ALE map $\mathcal A$, the 
fluid equations \eqref{eq:fluid} 
can be pulled back to the fixed reference domain $\Omega^{\rm f}_0$, which yields the system 
\begin{equation}\label{eq:fluid-ale}
\left\{ 
\begin{aligned}
    \varrho^{\rm f} J^{\rm f} \left( \pder{\vec{v}^{\rm f}}{t} + {({\boldsymbol \nabla}\vec{v}^{\rm f})}\,(\mathbb{F}^{\rm f})^{-1}\bigl(\vec{v}^{\rm f} - {\vec{w}}^{\rm f}\bigr) \right) - \Div \!\left( J^{\rm f} \mathbb{T}^{\rm f} (\mathbb{F}^{\rm f})^{-{\top}} \right) = \varrho^{\rm f} J^{\rm f} \vec{g}& \quad \text{in}\quad \Omega^{\rm f}_0, \\
    \div \!\left( J^{\rm f} (\mathbb{F}^{\rm f})^{-1} \vec{v}^{\rm f} \right) = 0& \quad \text{in}\quad \Omega^{\rm f}_0.
\end{aligned}\right.
\end{equation}
Note that, for the sake of simplicity, the same notation has been used for the Eulerian fields and their pullbacks into $\Omega^{\rm f}_0$ by $\mathcal{A}$. In particular, the pullback of the fluid Cauchy stress tensor is given by
$$
{
    \mathbb{T}^{\rm f}
    =
    -p^{\rm f}\mathbb{I}
    +\mu^{\rm f}\left(
        {\boldsymbol \nabla}\vec{v}^{\rm f}(\mathbb{F}^{\rm f})^{-1}
        +(\mathbb{F}^{\rm f})^{-\top}({\boldsymbol \nabla}\vec{v}^{\rm f})^\top
    \right).
}
$$
 {The coupling of \eqref{eq:fluid-ale} with the solid equations \eqref{eq:solid}
in Lagrangian form, is operated 
through the following 
interface conditions   on $\Gamma^{\rm i}_0$:
\begin{equation}\label{eq:fsi-coupling-ale}
\left\{
\begin{aligned}
    \vec{v}^{\rm f}= \pder{\vec{u}^{\rm s}}{t}  \quad \text{on} \quad \Gamma^{\rm i}_0,&\\
    \mathbb{P}\vec{n}_0^{\rm s}  = J^{\rm f}
    \mathbb{T}^{\rm f}    (\mathbb{F}^{\rm f})^{-\top} \vec{n}_0^{\rm s}  
  \quad \text{on} \quad \Gamma^{\rm i}_0,&\\
  \vec{u}^{\rm f}= \mathcal L(\vec{u}^{\rm s}\vert_{\Gamma_0^{\rm i}}),&
  \end{aligned}\right. 
\end{equation} 
where $\mathcal L$ denotes an arbitrary smooth lifting operator from $\Gamma_0^{\rm i}$ into $\Omega^{\rm f}_0$.
}
A common and simple choice is to define  {$\mathcal L$ as a harmonic lifting operator:}
$$
\left\{ 
\begin{aligned}
    -\Delta \vec{u}^{\rm f} = \vec{0} &\quad \text{in}\quad \Omega^{\rm f}_0,\\ 
    \vec{u}^{\rm f} = \vec{u}^{\rm s}&\quad \text{on } \quad\Gamma^{\rm i}_0, \\
    \vec{u}^{\rm f} = \vec{0} &\quad \text{on}\quad \Gamma_{\rm bot}\cup\Gamma_{\rm wall}.
    \end{aligned}
    \right. 
$$
Other choices include {nonlinear harmonic,} pseudo-elasticity problems or biharmonic extensions, which may better preserve mesh quality for large deformations~\cite{Wi11,FarSchTum24,SteinTezduyarBenney2003}.

 {As regards the choice of the reference domain $\Omega_0^{\rm f}$, two variants of the approaches can be found in the literature.} In the {\it Full ALE} method~\cite{TurekHron2006,HronTurek2006benchmark,Wi11,Richter2017}, the reference domain $\Omega^{\rm f}_0$ is kept fixed throughout the entire simulation, and the ALE displacement $\vec{u}^{\rm f}$ accumulates all deformations from the initial configuration. In the {\it Updated ALE} method~\cite{SteinTezduyarBenney2003,FarSchTum24}, the computational domain is periodically replaced by {the deformed configuration} at a {remeshing} time $t_r$. The total deformation gradient is then decomposed multiplicatively as
\begin{equation}\label{eq:updated_ale_split}
    \mathbb{F} = \mathbb{F}_c\,\mathbb{F}_r,
\end{equation}
where $\mathbb{F}_r$ encodes the deformation accumulated up to $t_r$ and $\mathbb{F}_c = \mathbb{I} + \nabla_{{X_r}} \vec{u}_c$ describes the deformation relative to the last computational domain $\Omega^{\rm f}_{t_r}$. Both variants will be employed in the numerical methods presented in Appendix~\ref{sec:appendix}.


\paragraph{Weak formulation.}

The coupling conditions can be {incorporated variationally} by choosing global test and trial functions in $\Omega$,
$$
\boldsymbol{\phi}=\left\{\begin{aligned}
&\boldsymbol{\phi}^{\rm f}\quad \text{in}\quad \Omega_0^{\rm f}\\
&\boldsymbol{\phi}^{\rm s}\quad \text{in}\quad \Omega_0^{\rm s},\\
\end{aligned}\right.,\quad 
\vec{v}= \left\{\begin{aligned}
&\vec{v}^{\rm f}\quad \text{in}\quad \Omega_0^{\rm f}\\
&\vec{v}^{\rm s}\quad \text{in}\quad \Omega_0^{\rm s},\\
\end{aligned}\right. 
$$
{with ${\boldsymbol \phi}^{\rm f}={\boldsymbol \phi}^{\rm s}$ and $\vec{v}^{\rm f}=\vec{v}^{\rm s}$ on $\Gamma_0^{\rm i}$},
resulting in the following variational formulation: find $(\vec{v}, p^{\rm f}, \vec{u}^{\rm f}, \vec{u}^{\rm s})$ such that
\begin{align}
 \int_{\Omega^{\rm f}_0} \varrho^{\rm f} J^{\rm f} \left(\pder{\vec{v}}{t} + {({\boldsymbol \nabla}\vec{v})}(\mathbb{F}^{\rm f})^{-1}(\vec{v} - \vec{w}^{\rm f})\right)\cdot{\boldsymbol \phi}
  + \int_{\Omega^{\rm f}_0} J^{\rm f} \mathbb{T}^{\rm f} (&\mathbb{F}^{\rm f})^{-\rm T} : {\boldsymbol \nabla}{\boldsymbol \phi} \nonumber \\ 
  +\int_{\Omega^{\rm s}_0} \varrho^{\rm s} \pder{\vec{v}}{t}\cdot{\boldsymbol \phi}
  + \int_{\Omega^{\rm s}_0} \mathbb{P} : {\boldsymbol \nabla}{\boldsymbol \phi} &= \int_{\Omega^{\rm f}_0} \varrho^{\rm f} J^{\rm f} \vec{g}\cdot{\boldsymbol \phi} + \int_{\Omega^{\rm s}_0} \varrho^{\rm s} \vec{g}\cdot{\boldsymbol \phi},\label{varALE1}\\
  \int_{\Omega^{\rm f}_0} J^{\rm f} (\mathbb{F}^{\rm f})^{-1}\vec{v} \cdot {\boldsymbol \nabla} q = \int_{\Omega^{\rm f}_0} {\boldsymbol \nabla} \vec{u}^{\rm f} : {\boldsymbol \nabla} \boldsymbol{\psi}^{\rm f} =  \int_{\Omega^{\rm s}_0} \left(\pder{\vec{u}^{\rm s}}{t} - \vec{v}\right)\cdot\boldsymbol{\zeta}^{\rm s} &=0.\label{eq:ale::uv}
\end{align}
for all admissible test functions $({\boldsymbol \phi}, q, \boldsymbol{\psi}^{\rm f}, \boldsymbol{\zeta}^{\rm s})$.

\subsection{{Coupled problem in fully Eulerian formalism}}
\label{subsec:Eulerian}
In the fully Eulerian approach~\cite{Dunne2006,Richter2013,
RichterWick2010,FreiRichterWick2016,FreiRichter2017, GraSchSouTum22,FarSchTum24}, both fluid and solid are formulated in the current (Eulerian) configurations. After reformulating the solid equations in the current configuration, the
solid is described in $\Omega^{\rm s}_t$ by the velocity field $\vec{v}^{\rm s}$ and, in
the displacement-based variants considered here, also by an Eulerian
displacement field $\vec{u}^{\rm s}$, {which is defined by the relation $\vec{u}^{\rm s}(\vec{X},t) \eqd \vec{X}-\vec{x}(\vec{X},t)$, where $\vec{x}(\vec{X},t) \in \Omega_0^{\rm s}$ denotes the material point which corresponds to the current point $\vec{X}\in\Omega_t^{\rm s}$. The domain $\Omega_t^{\rm s}$, as well as the interface $\Gamma_t^{\rm i}$ are unknown a priori and depend on the displacement $\vec{u}^{\rm s}(\vec{X},t)$. They can be captured via
\begin{align*}
    \Omega_t^{\rm s} &\eqd  \left\{ \vec{X}\in\Omega: \, \vec{X}-\vec{u}^{\rm s}(\vec{X},t) \in \Omega_0^{\rm s} \right\}, \quad
     \Gamma_t^{\rm i} \eqd  \left\{ \vec{X}\in\Omega: \, \vec{X}-\vec{u}^{\rm s}(\vec{X},t) \in \Gamma_0^{\rm i} \right\},
\end{align*}
As in Section~\ref{sec.lageul}, the fluid domain $\Omega_t^{\rm f}$ is then given by $\Omega_t^{\rm f}=\Omega\setminus\overline{\Omega_t^{\rm s}}$.
}
The displacement field is transported with the
solid motion and satisfies the kinematic relation
\begin{equation}\label{eq:eulerian-kinematic}
\pder{\vec{u}^{\rm s}}{t} + \vec{v}^{\rm s} \cdot {\boldsymbol \nabla}\vec{u}^{\rm s} = \vec{v}^{\rm s} 
\quad \mbox{in}\quad \Omega_t^{\rm s}.
\end{equation}
The deformation gradient is then recovered from the Eulerian displacement as
$$
    \mathbb{F} = \left(\mathbb{I}-\nabla\vec{u}^{\rm s}\right)^{-1}.
$$
The solid momentum balance equation \eqref{eq:solid} formulated in the current configuration reads
\begin{equation}\label{eq:solid-eulerian}
\tilde\rho^{\rm s} \left( \pder{\vec{v}^{\rm s}}{t} +   \vec{v}^{\rm s} \cdot {\boldsymbol \nabla}\vec{v}^{\rm s}\right) - \Div \,\mathbb{T}^{\rm s} = \tilde\rho^{\rm s} \vec{g}
\quad \text{in}\quad \Omega_t^{\rm s},
\end{equation}
where $\mathbb{T}^{\rm s}$ is the Cauchy stress tensor of the solid and  $\tilde{\rho}^{\rm s}\eqd J^{-1}\rho^{\rm s}$. It is related to
the first Piola--Kirchhoff stress tensor $\mathbb{P}$ by the standard
Piola transformation
$$
    \mathbb{T}^{\rm s} \eqd  J^{-1} \mathbb{P}\,\mathbb{F}^{\rm T}
    =J^{-1} \mathbb{F}\mathbb{S}\mathbb{F}^{\rm T}
    =J^{-1} \left( \lambda^{\rm s} (\tr\, \mathbb{E})\, \mathbb{F}
    + 2 \mu^{\rm s} \mathbb{F}\,\mathbb{E} \right)\mathbb{F}^{\rm T}.
$$

{The coupling of Eulerian fluid and solid problems 
 \eqref{eq:fluid} and  \eqref{eq:solid-eulerian}
 is operated 
 through the following  conditions on the moving interface $\Gamma_t^{\rm i}$:
 \begin{equation}\label{eq:fsi-coupling-eulerian}
\left\{
\begin{aligned}
    \vec{v}^{\rm f}= \vec{v}^{\rm s}  \quad \text{on} \quad \Gamma^{\rm i}_t ,&\\
    \mathbb{T}^{\rm s} \vec{n} =  
    \mathbb{T}^{\rm f}   \vec{n}
  \quad \text{on} \quad \Gamma^{\rm i}_t ,&
  \end{aligned}\right. 
\end{equation} 
}

The coupling conditions can again be {incorporated variationally} by choosing global test and trial functions in $\Omega$, as in~\cite{Dunne2006,Richter2013,
RichterWick2010,FreiRichterWick2016,FreiRichter2017},
$$
\boldsymbol{\phi}=\left\{\begin{aligned}
&\boldsymbol{\phi}^{\rm f}\quad \text{in}\quad \Omega_t^{\rm f}\\
&\boldsymbol{\phi}^{\rm s}\quad \text{in}\quad \Omega_t^{\rm s}\\
\end{aligned}\right.,\quad 
\vec{v}= \left\{\begin{aligned}
&\vec{v}^{\rm f}\quad \text{in}\quad \Omega_t^{\rm f}\\
&\vec{v}^{\rm s}\quad \text{in}\quad \Omega_t^{\rm s}.\\
\end{aligned}\right. 
$$
{The traces of the global functions satisfy ${\boldsymbol \phi}^{\rm f}={\boldsymbol \phi}^{\rm s}$ and $\vec{v}^{\rm f}=\vec{v}^{\rm s}$ on $\Gamma_t^{\rm i}$.}
Setting $\tilde\varrho^{\rm f} \eqd \varrho^{\rm f}$, the resulting variational formulation reads: find $(\vec{v}, p^{\rm f}, {\vec{u}^{\rm s}})$ such that
\begin{align}
\begin{split}
\label{varfullyEulerian1}
  \sum_{i\in \{\rm f,s\}} \left(\int_{\Omega^i_t} \tilde\varrho^i \pder{\vec{v}^i}{t}\cdot {\boldsymbol \phi}^i + \int_{\Omega^i_t} \tilde\varrho^i  \vec{v}^i \cdot {\boldsymbol \nabla}\vec{v}^i  \cdot {\boldsymbol \phi}^i + \int_{\Omega^i_t} \mathbb{T}^i : {\boldsymbol \nabla}{\boldsymbol \phi}^i \right)  -\int_{\Omega^{\rm f}_t} q\,\div\,\vec{v}^{\rm f}  
  &=  \sum_{i\in \{f,s\}} \int_{\Omega^i_t} \tilde\varrho^i\vec{g}\cdot{\boldsymbol \phi}^i,
  \end{split}\\
  \int_{\Omega^{\rm s}_t} \left(\pder{\vec{u}^{\rm s}}{t} + \vec{v}^{\rm s} \cdot {{\boldsymbol \nabla}\vec{u}^{\rm s}  - \vec{v}^{\rm s}} \right)\cdot\vec{z}^{\rm s} &= 0,
\end{align}
for all admissible test functions $({\boldsymbol \phi}, q, \vec{z}^{\rm s})$.
{On the discrete level}, the continuity of velocities can also be enforced {weakly} by using Nitsche's method (see,  e.g.,  \cite{BurmanFernandezFrei2020, BurmanFernandezFreiGerosa2022, FreiKnokeSteinbachWenskeWick2026}).


\subsection{{Solid contact modeling}}
\label{subsec.contact}

As discussed in the introduction, the FSI model  \eqref{eq:solid}-\eqref{eq:bc-fluid} introduced in Section~\ref{sec.lageul} in combination with no-slip conditions does not allow for contact. To reproduce this behavior numerically, however, typically requires a very fine mesh resolution in the {near-contact} region. If the mesh is too coarse, the fluid forces might not be approximated accurately enough to prevent contact. As a result, the computation might either break down or yield non-physical results. 

\paragraph{Relaxed contact conditions.} A popular remedy to avoid this issue and to enable computations also on coarser meshes is to introduce no-penetration conditions that prevent solid-solid contact (see, e.g.~\cite{BurmanFernandezFrei2020, BurmanFernandezFreiGerosa2022, FormaggiaGattiZonca2021, FrKnStWeWi25}). To this end, we impose a minimal distance $\epsilon>0$ between the solid bodies. Numerically $\epsilon$ is typically chosen depending on the mesh size $h$ {(e.g., $\mathcal O(h)$)}. Denoting the initial distance from the interface $\Gamma_0$ to the ground by a function $g_0$, the resulting no-penetration conditions can be written in Lagrangian formalism as
\begin{equation}\label{eq:contact_ineq2}
\vec{u}^{\rm s}\cdot {\vec{n}_{\rm b}}  \leq g_\epsilon, \quad \lambda \leq 0, \quad \lambda (\vec{u}^{\rm s}\cdot {\vec{n}_{\rm b}} - g_\epsilon)  = 0 \quad \text{on}\quad {\Gamma_0^{\rm i}},
\end{equation}
where $g_\epsilon=g_0-\epsilon$ and $\lambda$ is {the Lagrange multiplier associated to the constraint (see,  e.g., ~\cite{ChoulyHildRenard2023}) and 
$\vec{n}_{\rm b}$ denotes the (constant) outward unit normal vector to the bottom wall.} Using the Alart--Curnier trick~\cite{AlartCurnier91,ChoulyHildRenard2023}, these conditions can be equivalently reformulated as {the non-linear} identity 
\begin{align}\label{AlartCurnier}
\lambda = - \big|  {\gamma_{\rm C}  (\vec{u}^{\rm s} \cdot \vec{n}_{\rm b} - g_\epsilon) - \lambda} \big|_+ \quad \text{on}\quad {\Gamma_0^{\rm i}}
\end{align}
for arbitrary $\gamma_{\rm C}>0$, where   $|\,x\,|_+\eqd \max\{x,0\}$ denotes the positive part of $x\in\mathbb{R}$. 

{As shown in \cite{BurmanFernandezFrei2020}, the addition of the contact constraint \eqref{eq:contact_ineq2} to the coupled problem \eqref{eq:solid}-\eqref{eq:bc-fluid}
 induces a perturbation of the dynamic coupling on the interface, through  the relation  
\begin{equation}\label{eq:mech-inconsistency}
\lambda = 
\vec{n}_{\rm b}^{\top} \left( \mathbb{P}    
- 
J
    \mathbb{T}^{\rm f}  \mathbb{F}^{-\top}   \right) \vec{n}_0^{\rm s}
\quad \mbox{on}\quad\Gamma_0^{\rm i}.
\end{equation}
This mechanical inconsistency will be  addressed  in Section~\ref{subsec.porous}.
Note that \eqref{eq:mech-inconsistency} and  
the original dynamic continuity in \eqref{eq:fsi-coupling} coincide in the absence of contact ($\lambda = 0$).
}

{
The non-linear identity 
\eqref{AlartCurnier} is the starting point 
for   Nitsche and Augmented Lagrangian formulations of contact. 
In Nitsche based approaches (see~\cite{BurmanFernandezFrei2020,ChoulyHildRenard2023}), the Lagrange multiplier is eliminated  by combining \eqref{AlartCurnier} with  \eqref{eq:mech-inconsistency}, 
whereas 
in augmented Lagrangian formulations (see, e.g., \cite{POULIOS201575,ChoulyHildRenard2023}), the
multiplier $\lambda$ serves as  the unknown associated 
to the  additional relation  \eqref{AlartCurnier}. 
}

While \eqref{eq:contact_ineq2} can be used as a \textit{numerical safeguard} to prevent solid-solid contact, when no contact is expected, it has also been proposed in cases, where solid-solid contact is expected. However, it has been observed in~\cite{BurmanFernandezFreiGerosa2022} that bounces get smaller and smaller under mesh refinement $h\to 0$, which is in accordance with the no-contact paradox. We will investigate this effect in detail in the numerical examples reported {in Section~\ref{results}} and will in fact find convergence to the contactless rebound, when $\epsilon(h)\to 0$ as $h\to 0$.

\subsection{{Porous media based roughness modeling with seepage}} \label{subsec.porous}
{In order to avoid the no-contact paradox and the mechanical inconsistency induced by \eqref{eq:mech-inconsistency},} different modifications of the contact model have been proposed for cases where solid-solid contact is expected, mostly based on the assumption of rough surfaces of the contacting solids. In~\cite{AgerSchottVuongPoppWall2019}, {a non-linear poroelastic 
modeling of  the surface roughness  has been investigated}, which is placed within a domain of thickness $\epsilon_{\rm p}>0$ between the contacting structures. Similarly, the works \cite{BurmanFernandezFrei2020} and~\cite{GerosaMarsden2024} propose to add an artificial fluid governed by a Navier-Stokes-Brinkman equation inside one of the solid bodies resp.\,between them. Finally,~\cite{BurmanFernandezFreiGerosa2022,champion:hal-05626925} propose a Darcy-type equation between the contacting solids, but in the limit case of infinitesimal thickness $\epsilon_{\rm p}\to 0$. 

In~\cite{BurmanFernandezFreiGerosa2022}, the FSI-contact model is enriched by one additional variable, namely the porous pressure $P_{\rm l}:\Gamma_{\text{bo{t}}}\times \mathbb{R}^+\rightarrow \mathbb{R}$, which is governed by the following system 
\begin{equation}
\label{eq:darcy}
\left\{ 
\begin{aligned}
- {\boldsymbol \nabla}_{\boldsymbol \tau} \cdot \left( \epsilon_{\rm p} K_{\boldsymbol \tau} {\boldsymbol \nabla}_{\boldsymbol \tau} P_{\rm l}  \right) =   \vec{v}^{\rm f} \cdot \vec{n} & \quad\mbox{on}\quad \Gamma_{\text{bo{t}}},\\
\epsilon_{p} K_{\boldsymbol \tau}  \boldsymbol \tau \cdot {\boldsymbol \nabla}_{\boldsymbol \tau} P_{\rm l} =  0  &\quad\mbox{on}\quad \partial\Gamma_{\text{bo{t}}}.
\end{aligned}\right. 
\end{equation}
Here,  $K_{ \boldsymbol \tau}$ and $K_{\boldsymbol n}$ respectively describe the conductivity of the porous layer in tangential and normal directions. The porous layer is coupled to the fluid by the following equations based on the Beavers--Joseph--Saffman conditions
\begin{equation}
\label{eq:darcy-coupling}
\mathbb{T}^{\rm f}   \vec{n} =- \left( P_{\rm l}  + \frac{\epsilon_{p} K_n^{-1} }{4} \vec{v}^{\rm f} \cdot \vec{n}\right) \vec{n}
-\frac{\alpha}{\sqrt{K_{\boldsymbol \tau}\epsilon_{\rm p}}} \vec{v}^{\rm f}\cdot \bs{\tau} \bs{\tau}
\quad\mbox{on}\quad \Gamma_{\text{bot}}.
\end{equation}
As discussed  in the introduction, it has been shown  in~\cite{ChampionFernandezetal2024} that  the resulting FSI-contact model enables contact in finite time.

{
 We will also consider the alternative model reported in 
 \cite{champion:hal-05626925}, in which the porous layer is attached to the solid boundary $\Gamma_t^{\rm i}$, instead of to the bottom wall $\Gamma_{\rm bot}$. 
 The dynamics of the porous pressure are given by 
 \begin{equation}
\label{eq:darcy-moving}
- {\boldsymbol \nabla}_{\boldsymbol \tau} \cdot \left( \epsilon_{\rm p} K_{\boldsymbol \tau} {\boldsymbol \nabla}_{\boldsymbol \tau} P_{\rm l}  \right) =   \vec{v}^{\rm f} \cdot \vec{n} - 
\vec{v}^{\rm s} \cdot \vec{n}
 \quad\mbox{on}\quad \Gamma_t^{\rm i}, 
\end{equation}
complemented with the following  dynamic coupling conditions for the fluid and the solid: 
\begin{eqnarray} 
\begin{aligned}\label{eq:por:int}
&  \mathbb{T}^{\rm f}  \vec{n}  = 
- \left( P_{\rm l} \vec{n}  + \frac{\alpha}{\sqrt{\epsilon_{\rm p} K_{\boldsymbol \tau}}} (\vec{v}^{\rm f} -\vec{v}^{\rm s} )\cdot \bs \tau \bs \tau   \right)  \quad \text {on}\quad \Gamma_t^{\rm i},\\
& \int_{\Gamma_0^{\rm i}} (\mathbb{P}\vec{n}_0^{\rm s}  - \lambda \vec{n}_{\rm b}  )  
\cdot {\boldsymbol \phi}^{\rm s}  =   \int_{\Gamma_t^{\rm i}     }  
\left(P_{\rm l} \vec{n} -  \frac{\alpha}{\sqrt{\epsilon_{\rm p} K_{\boldsymbol \tau}}}  (\vec{v}^{\rm s} - \vec{v}^{\rm f} )\cdot {\bs \tau} \bs \tau \right) \cdot  {\boldsymbol \phi}^{\rm s}  \circ \vec{X}_t^{-1}
\end{aligned}
\end{eqnarray}
for any admissible displacement ${\boldsymbol \phi}^{\rm s}:\Omega^{\rm s}_0\rightarrow \mathbb{R}^2$.}


%
%

\subsection{{Summary of the models and energy balance}}
\label{sec.energy}

{
In the above sections, we have introduced three different fluid-structure-contact interaction modeling options: 
\begin{itemize}
\item Model I (full no-slip): \eqref{eq:solid}-\eqref{eq:bc-fluid} or the alternative formulations \eqref{eq:fluid-ale}-\eqref{eq:fsi-coupling-ale} or \eqref{eq:eulerian-kinematic}-\eqref{eq:fsi-coupling-eulerian}, possibly with
 the contact condition \eqref{eq:contact_ineq2};
 \item Model II (no-slip on the interface coupling 
 and Darcy on the bottom wall): 
 \eqref{eq:solid}-\eqref{eq:fsi-coupling}, 
 \eqref{eq:contact_ineq2}, 
\eqref{eq:darcy} and \eqref{eq:darcy-coupling};
\item Model III (no-slip on the   bottom wall and  
interface coupling via Darcy):
  \eqref{eq:solid}-\eqref{eq:fluid}, 
  \eqref{eq:contact_ineq2}, 
  \eqref{eq:darcy-moving} and \eqref{eq:por:int}.
\end{itemize}
These three coupled problems obey a total energy balance. We
state it here for two fundamental reasons: (i) it identifies the mechanisms, through which the
system can dissipate energy; and (ii) it provides a scalar output quantity whose time history is a convenient
and formulation-independent basis of comparison for the numerical approaches 
considered in this benchmark (see Section~\ref{results} below). }

{
In what follows, the solutions of these coupled problems are assumed to be
smooth enough for the manipulations to be licit. 
The total energy of the  system is defined as 
\begin{equation}\label{eq:total-energy}
    E(t) \eqd 
    {\frac{\varrho^{\rm f}}{2}\int_{\Omega^{\rm f}_t} \bigl|\vec{v}^{\rm f}(t)\bigr|^2} 
    + {\frac{\varrho^{\rm s}}{2}\int_{\Omega^{\rm s}_0} \Bigl|\pder{\vec{u}^{\rm s}}{t}(t)\Bigr|^2} 
    + {\int_{\Omega^{\rm s}_0} W\bigl(\mathbb{E}(\vec{u}^{\rm s}(t))\bigr)}
    - {\varrho^{\rm s}\int_{\Omega^{\rm s}_0} \vec{g}\cdot\bigl(\vec{I}+\vec{u}^{\rm s}(t)\bigr)} 
    + {\varrho^{\rm f}\int_{\Omega^{\rm s}_t} \vec{g}\cdot\vec{I}} ,
 \end{equation}
 where $\vec{I}$ denotes the identity map in $\mathbb{R}^2$, i.e., $\vec{I}(\vec{x}) = \vec{x}.$
 The first and second terms correspond to the fluid
 and solid kinetic energies, respectively. 
 The third stands for the solid elastic energy, 
 while the last two terms represent the solid and 
 displaced-fluid potential energies.
   Note  that, although
\eqref{eq:total-energy} is written in the mixed Lagrangian--Eulerian variables
of Section~\ref{sec.lageul}, the energy is a physical quantity independent of
the formalism; in the fully Eulerian setting of Section~\ref{subsec:Eulerian}
the last three terms are simply evaluated as
$$\int_{\Omega^{\rm s}_t} J^{-1} W(\mathbb{E}),\quad -\varrho^{\rm s}\int_{\Omega^{\rm s}_t} J^{-1}\vec{g}\cdot\vec{I}, \quad 
\varrho^{\rm f}\int_{\Omega^{\rm s}_t}\vec{g}\cdot\vec{I}.$$ }

{
In the porous models of Section~\ref{subsec.porous} the boundary carrying the
layer is permeable, so that the fluid convective term of \eqref{eq:fluid} yields
a contribution of indefinite sign (viz., energy can be fed into the system wherever
the flow re-enters the fluid domain). We therefore assume (see, e.g.,  \cite{ChampionFernandezetal2024,MoghadamEtAl2011}), that the traction conditions
\eqref{eq:darcy-coupling} and \eqref{eq:por:int} are supplemented with the
back-flow stabilization term
$\frac{\varrho^{\rm f}}{2}\bigl|\vec{v}^{\rm f}\cdot\vec{n}\bigr|_{-}\vec{v}^{\rm f}$,
respectively
$\frac{\varrho^{\rm f}}{2}\bigl|(\vec{v}^{\rm f}-\vec{v}^{\rm s})\cdot\vec{n}\bigr|_{-}\vec{v}^{\rm f}$,
where $|\,x\,|_{-}\eqd \max\{-x,0\}$ denotes the negative part, so that the
stabilization acts only where the fluid enters the domain. Under this
assumption, the three coupled problems satisfy the energy balance
\begin{equation}\label{eq:energy-balance}
  \frac{{\rm d}}{{\rm d}t} E (t) = -2\mu^{\rm f} \int_{\Omega^{\rm f}_t}
  \bigl|\bs{\varepsilon}(\vec{v}^{\rm f})\bigr|^2 - \mathcal{D}(t),
\end{equation}
where the interfacial  dissipation $\mathcal{D}\geq 0$ is given by
\begin{equation}\label{eq:interface-dissipation}
  \mathcal{D} \eqd
  \begin{cases}
    0, & \text{(Model I)}\\
    \begin{aligned}[t]
      &\epsilon_{\rm p} K_{\boldsymbol \tau} \int_{\Gamma_{\text{bo{t}}}}\! \bigl|{\boldsymbol \nabla}_{\boldsymbol \tau} P_{\rm l}\bigr|^2
      + \frac{\epsilon_{\rm p} K_n^{-1}}{4} \int_{\Gamma_{\rm bot}} \bigl|\vec{v}^{\rm f}\!\cdot\vec{n}\bigr|^2+ \frac{\alpha}{\sqrt{K_{\boldsymbol \tau}\epsilon_p}} \int_{\Gamma_{\text{bo{t}}}} \bigl|\vec{v}^{\rm f}\cdot\bs{\tau}\bigr|^2\\
      &\quad 
      + \frac{\varrho^{\rm f}}{2} \int_{\Gamma_{\rm bot}}
 \bigl|\vec{v}^{\rm f}\!\cdot\vec{n}\bigr|_+ 
        \bigl|\vec{v}^{\rm f}\bigr|^2,
    \end{aligned}
      & \text{(Model II)}\\[6mm]
    \begin{aligned}[t]
      &\epsilon_{\rm p} K_{\boldsymbol \tau} \int_{\Gamma^{\rm i}_t }\bigl|{\boldsymbol \nabla}_{\boldsymbol \tau} P_{\rm l}\bigr|^2
      + \frac{\alpha}{\sqrt{K_{\boldsymbol \tau}\epsilon_p}} \int_{\Gamma^{\rm i}_t } \bigl|(\vec{v}^{\rm f}-\vec{v}^{\rm s})\cdot\bs{\tau}\bigr|^2+ \frac{\varrho^{\rm f}}{2} \int_{\Gamma^{\rm i}_t }
         \bigl|(\vec{v}^{\rm f}- \vec{v}^{\rm s})\cdot\vec{n}\bigr|_{+}  
        \bigl|\vec{v}^{\rm f}\bigr|^2.
    \end{aligned}
      & \text{(Model III)}
  \end{cases}
\end{equation}
 In particular, the
system is unconditionally dissipative, $E(t)\leq E(0)$ for all $t\geq 0$, irrespective of the density
ratio and of the magnitude of the interface displacement.
A sketch of the proof of \eqref{eq:energy-balance}
is given in Appendix~\ref{app:proof}.
}

\section{Benchmark description}
\label{sec:benchmark:description}
In this section, we describe the benchmark configuration and the quantities of interest.

\subsection{Configuration}

\paragraph{Domain description.}
 The computational domain is the rectangle $\Omega = (-R,R)\times(0,H)$ with $H = 0.08\,$m and $R = 0.04\,$m, so that $\Gamma_{\mathrm{bot}} = (-R,R)\times\{0\}$. The diameter of the solid is $d_{\rm s} = 2r_{\rm s} = 0.022\,$m, and $c_{\rm s} = (0, h_0 + r_{\rm s})$ with $h_0 = 0.039\,$m denotes the initial position of the disk center, i.e., $h_0$ is the initial gap between the disk and $\Gamma_{\mathrm{bot}}$.

\paragraph{Fluid properties.}
We consider a mixture consisting of glycerine and water at equal volume fractions. The constant density of the fluid is $\varrho^{\rm f}=1141$\,kg/m$^3$, the dynamic shear viscosity is  $\mu^{\rm f}=0.008$\,Pa\,s, and the gravitational acceleration  $\vec{g}=(0,-9.807)$\,m/s$^2$. 

\paragraph{Solid properties.}
 We assume that $\varrho^{\mathrm s} = 1361\,\mathrm{kg/m^3}$, $\mu^{\mathrm s} = 666.711\,$kPa and $\lambda^{\mathrm s} = 10.0\,$MPa, corresponding to a Young modulus $E = 1.958\,$MPa and a Poisson ratio $\nu = 0.4687$ under the plane strain convention.

\paragraph{Porous layer properties.} The parameters of the porous medium  in Section~\ref{subsec.porous} are chosen as $\epsilon_{\rm p}= 10^{-6}$m, $K_{\boldsymbol n} =K_{\boldsymbol \tau} = 10^{-5} {\rm m^3 s}/{\rm kg}$ and $\alpha = 1$ kg$^\frac{1}{2}$s$^{-\frac{1}{2}}$.

\subsection{Quantities of interest}
\label{sec.qoi}

{We propose to evaluate the following quantities of varying complexity {(all the units are expressed in the SI system)}:}
\begin{center}
\begin{tabular}{|lp{13.8cm}|}
\hline
$t^{*}$& We define $t^{*}$ as the time after release when ${\rm dist}({\Gamma_t^{\rm i}},\Gamma_{\rm bot})=r_s$. 
\\
$v^{*}$& The velocity of the disk in the $y$-direction at $t=t^{*}$ .\\
$f^{*}$& The vertical component of the force acting on the disk at time $t^{*}$. 
\\
$t_{\rm cont}$& The time of the first {near-contact event} (i.e., the time when the solid is closest to the bottom before the first rebound).\\  
$d_{\rm cont}$& The minimum of ${\rm dist}(\Gamma_t^{\rm i},\Gamma_{\rm bot})$ at $t=t_{\rm cont}$.\\  
{$f_{\rm max}$}& {The temporal maximum of the vertical component of the force acting on the disk before {the first near-contact event}, i.e.\ within the interval $[0,t_{\rm cont}]$.} \\ 
{$\max_{t} p_{\rm bc}$}& The maximum pressure at the center on the bottom wall over time.\\ 
$t_{\rm jump}$& The time of the maximum of ${\rm dist}({\Gamma_t^{\rm i}},\Gamma_{\rm bot})$ after {$t_{\rm cont}$}. \\ 
$d_{\rm jump}$& The maximum of ${\rm dist}({\Gamma_t^{\rm i}},\Gamma_{\rm bot})$ after {$t_{\rm cont}$}, i.e., the size of the {rebound}. \\ 
$E_{\rm el,s}$& The elastic energy in the solid over time $t$. \\ 
$E_{\rm k,s}$& The kinetic energy in the solid over time $t$.\\ 
$E_{\rm k,f}$& The kinetic energy in the fluid over time $t$. \\ 
\hline
\end{tabular}
\end{center}
Note that the first three quantities are computed significantly before the {first near-contact event} and can thus serve as a first simple test for an FSI code to {evaluate} if the falling period is captured accurately.

The elastic and kinetic energies are worth recording for every time $t$ {to investigate} the energy exchange. {Owing to \eqref{eq:total-energy}, these quantities} are defined as follows:
$$
E_{\rm el,s} \eqd \int_{\Omega_0^{\rm s}} \left(\frac{\lambda^{\rm s}}{2}(\tr\,\mathbb{E})^2+\mu^{\rm s}\left|\mathbb{E}\right|^2\right),\quad  
E_{\rm k,s} \eqd \int_{\Omega_0^{\rm s}} \frac12{\varrho^{\rm s}} \left|\vec{v}^{\rm s}\right|^2,\quad
E_{\rm k,f} \eqd \int_{\Omega_t^{\rm f}} \frac12{\varrho^{\rm f}} |\vec{v}^{\rm f}|^2.
$$
Moreover, we can define the gravitational potential energy of the solid relative to the fluid by means of 
$$
    E_{\rm grav} \eqd   -{\varrho^{\rm s}\int_{\Omega^{\rm s}_0} \vec{g}\cdot\bigl(\vec{I}+\vec{u}^{\rm s}(t)\bigr)} 
    + {\varrho^{\rm f}\int_{\Omega^{\rm s}_t} \vec{g}\cdot\vec{I}}
    =
    \int_{\Omega_t^{\rm s}} ({\varrho^{\rm f}}-\tilde{\varrho}^{\rm s}) \vec{g} \cdot \vec{I}.
$$
The vertical component of the force acting on the disk, $f$, is defined by
$$
f \eqd - \int_{\Gamma_t^{\rm i}} {\mathbb{T}^{\rm f}}\vec{n}\cdot\vec{e}_2,
$$
 where $\vec{e}_2\eqd (0,1)$ and $\vec{n}$ denotes {again} the { fluid outward unit normal to $\Gamma_t^{\rm i}$}.

\section{Results}\label{results}

This section collects the numerical results obtained by the participating groups.  Besides the scalar values and figures reported below, the full time histories of the dynamic quantities (gap, drag force, pressure at center on the bottom wall and  energies) are provided for each approach and each refinement level as supplementary material~\cite{BenchmarkData}.

\subsection{Participating teams}

Five different teams have participated in the computation of the benchmark results, using different  numerical approaches and software packages. We give here a brief summary of the numerical approaches used. Further details are  provided in Appendix~\ref{sec:appendix}.

\begin{itemize}
\item[[W]] Wick uses a {plain vanilla} {variational-monolithic ALE approach} for Model~I (see Section~\ref{subsec.ALE}) with neither contact nor remeshing. The implementation is
  based on the FSI template~\cite{Wi13_fsi_with_deal}, which {builds} on the open-source software library \texttt{deal.II}~\cite{dealII96}. The codes are available open-source on github\footnote{\url{https://github.com/tommeswick/fsi}, \url{https://github.com/FlorianStahlhuth/FSI-FallingBall}}.
\item[[FST]] Fara, Schwarzacher \& T\r{u}ma use a more advanced monolithic {updated ALE approach} for Model~I (see Section~\ref{subsec.ALE}), with fitted fluid and solid meshes and a periodically applied remeshing strategy tailored to the {contactless-rebound regime}. The remeshing step preserves the fluid--solid interface, {{refines the mesh in the gap between the potentially contacting boundaries}}, improves mesh quality and is based on the ADmesh remeshing strategy described in~\cite{FarSchTum24}. The implementation {builds} on \texttt{FEniCS}~\cite{FEniCS}.
\item[[CF]] Corti \& Fern\'andez use a mixed Lagrangian–Eulerian coordinate approach for Model~I (see Section~\ref{sec.lageul}), with a CutFEM approximation in the fluid and fitted solid discretization, both coupled via Nitsche's method. A relaxed contact approach is applied imposing a minimal gap $\epsilon=ch$ to the lower ground using an augmented Lagrangian formulation. This yields the results [CF1].  Alternatively, instead of using no-slip conditions, [CF2] considers a Darcy model of surface roughness on $\Gamma_{\text{bo{t}}}$  (Model~II), as described in Section~\ref{subsec.porous}. Finally, [CF3] uses the Darcy layer on $\Gamma_{t}^{\rm i}$ (Model~III).
The implementation is based on the open-source \texttt{FELiScE}\footnote{\url{https://gitlab.inria.fr/felisce/felisce}} finite element library developed  at Inria.

\item[[Fr]] Frei uses a {fully Eulerian approach} for Model~I (see Section~\ref{subsec:Eulerian}). The discretization is based on a fitted locally modified finite element approach resolving the interface and a coupling approach based on Nitsche's method. The implementation {builds} on \texttt{Gascoigne3d}~\cite{Gascoigne3d}. As in [CF1] a relaxed contact approach is applied in the first approach [Fr1]. [Fr2] introduces a Darcy layer on $\Gamma_{\text{bo{t}}}$, as [CF2] (Model~II).

\item[[KWF]] Knoke, Wick and Frei also use a {fully Eulerian approach} for Model~I  (see Section~\ref{subsec:Eulerian}). The discretization is, however, based on an unfitted CutFEM approach (i.e., not resolving the interface) with Nitsche type coupling. As for the previous teams, a relaxed contact approach is used. The implementation is based on~\cite{FrKnStWeWi25}, which {builds} on the open-source software library \texttt{deal.II}~\cite{dealII96}.

\end{itemize}

\subsection{Quantitative Results}

{In Table~\ref{tab:res_qi}, we report the results {in the first nine QoIs} obtained with the eight approaches. For each approach, three different results were provided, which correspond to different meshes and possibly different time steps. As most groups use non-uniform meshes and time steps, we provide the minimum and maximum edge sizes $h_{\min}$ and $h_{\max}$ of the meshes and the minimum and maximum time steps $\delta t_{\min}$ and $\delta t_{\max}$, respectively, in the first {columns} of Table~\ref{tab:res_qi}. We observe large differences in the resolution, ranging from edge sizes of order $10^{-3}$ in [W] to minimum edge sizes of order $10^{-7}$ used by {[FST]}, and similarly for the time steps. {We will see below that the use of a highly refined mesh in the {near-contact} region, as well as {very small time steps} in the time interval around {the near-contact event}, is necessary to obtain accurate results.} The purpose of the results [W] is to illustrate that the benchmark can, in principle, be solved with a plain vanilla ALE-FSI method using relatively coarse meshes and coarse time step sizes
even though it leads to a substantial loss of accuracy, in particular with regard to the post-contact dynamics.}


\clearpage

\begin{landscape}
\begin{table}[h!]
\begin{adjustbox}{width=\linewidth}
\sisetup{round-mode=places,round-precision=4}
\begin{tabular}{rrr|rr|rrrrrrrrr}
&$h_{\min}$ &$h_{\max}$ &$\delta t_{\min}$ & $\delta t_{\max}$ &$t^*$ &$v^*$ &$f^*$ 
& $t_{\text{cont}}$ &
$d_{\text{cont}}$ &{$f_{\text{max}}$} & $\max_{t} p_{\rm bc}$ & $t_{\text{jump}}$ & $d_{\text{jump}}$ \\
\hline
[W] & \num[round-precision=2]{1.45137e-03} & \num[round-precision=2]{5.63648e-03} & \num[round-precision=2]{1.0e-04} 
& \num[round-precision=2]{1.00000e-02} 
& \num{3.00000e-01} & \num{-1.50900e-01}  &\num{4.55575e+00}  
& \num{3.74600e-01} &*0  
&\num{4.89810} & \num{4.04233e+04  } 
& \num{4.32700e-01} & \num{1.65750e-03} \\

[W] & \num[round-precision=2]{7.25687e-04} & \num[round-precision=2]{2.81824e-03} & \num[round-precision=2]{1.0e-04} & \num[round-precision=2]{1.00000e-02} & \num{3.00000e-01} & \num{-1.60936e-01} & \num{4.52285e+00} 
&\num{3.62400e-01} & \num{4.08344e-07} 
& \num{5.11397} & \num{9.11619e+04} 
&\num{4.78600e-01} & \num{4.41583e-03} \\

[W] & \num[round-precision=2]{3.62843e-04} & \num[round-precision=2]{1.40912e-03} & \num[round-precision=2]{1.0e-04} & \num[round-precision=2]{1.00000e-02} & \num{2.90000e-01} & \num{-1.61817e-01} & \num{4.48887e+00} 
& \num{3.60300e-01} & \num{1.73952e-06} 
&\num{5.06838} & \num{1.57929e+05} 
& \num{4.80000e-01} & \num{5.42029e-03} \\





\hline
\hline 

[CF1] &  \num[round-precision=2]{1.54e-5} &  \num[round-precision=2]{8e-4} &   \num[round-precision=2]{2e-4} &  \num[round-precision=2]{2e-4}  &  
\num{2.906e-1} &  \num{-1.7e-1}  & \num{4.96617}  
&\num{3.598e-1} &  *\num{4.8399e-7} 
&\num{4.36459e1} &  \num{6.74347e4} 
&\num{4.074e-1} & \num{6.29284e-4}\\

[CF1] &  \num[round-precision=2]{7.93e-6} &  \num[round-precision=2]{4e-4} &   \num[round-precision=2]{1e-4} &  \num[round-precision=2]{1e-4}  & 
\num{2.898e-1}  & \num{-1.7e-1}   & \num{4.96778}   
&   \num{3.591e-1} & *\num{1.62125e-7}   
&\num{6.33615e1} &\num{1.09394e5}
&\num{3.92e-1} & \num{3.52178e-4} \\
  
[CF1] &  \num[round-precision=2]{3.963e-6} &  \num[round-precision=2]{2e-4} &   \num[round-precision=2]{5e-5} &  \num[round-precision=2]{5e-5}  & \num{2.8925e-1} & \num{-1.7e-1} &  \num{4.96805} 
&  \num {3.5875e-1} & \num{6.41584e-6}  
&\num{6.4072e+1} &   \num{9.66673e4} 
&\num{4.265e-1} & \num{1.18102e-3} \\

\hline
[Fr1] &\num[round-precision=2]{6.25e-05} & \num[round-precision=2]{6.25e-04} &\num[round-precision=2]{1.1094e-5}&\num[round-precision=2]{1.42e-3} &\num{2.89058166118421e-1} &\num{-1.71457312911184e-1} &\num{4.94775338404605} 
&\num{3.57674e-1}
&*\num{7.70217e-06} 
&\num{3.23475e1} &\num{3.29293e+4} 
&\num{4.503955e-1} &\num{2.45388e-3}
\\ 

[Fr1] &\num[round-precision=2]{4.6875e-05} & \num[round-precision=2]{4.6875e-04} &\num[round-precision=2]{7.8125e-6} &\num[round-precision=2]{1e-3} &\num{2.88905792861322e-1} &\num{-1.71103216500878e-1} &\num{4.95097266822703} 
&\num{3.5782e-1}
&*\num{5.61241e-06} 
&\num{4.58057e1} &\num{5.01904e+4} 
&\num{4.31133e-1}  &\num{1.29332e-3}\\

[Fr1] &\num[round-precision=2]{3.125e-05} & \num[round-precision=2]{3.125e-04} &\num[round-precision=2]{5.5469e-6}&\num[round-precision=2]{7.1e-4} &\num{2.88786793075021e-1} &\num{-1.70912584501237e-1} &\num{4.95395937345425} 
&\num{3.57768e-1}
&*\num{4.05171e-06} 
&\num{5.55129e1} &\num{6.43345e+4} 
&\num{4.17757e-1} &\num{9.3074e-4}\\
\hline

[KWF] & \num[round-precision=2]{3.75000e-05} & \num[round-precision=2]{1.87500e-03} & \num[round-precision=2]{3.12500e-05} & \num[round-precision=2]{1.00000e-03} & \num{2.89000e-01} & \num{-1.7169e-01} & \num{4.9558e+00}
& \num{3.547812e-01} & \num{1.6723024e-05} 
&\num{6.7414473e1} & \num{1.8293e+05} 
& \num{4.27e-01} & \num{1.3472692e-03} \\

[KWF] & \num[round-precision=2]{3.75000e-05} & \num[round-precision=2]{1.87500e-03} & \num[round-precision=2]{1.5625e-05} & \num[round-precision=2]{5.00000e-04} & \num{2.89000e-01} & \num{-1.7115e-01} & \num{4.9534e+00} 
& \num{3.5575e-01} & \num{1.6701549e-05} 
&\num{6.7703983e1} & \num{1.8374e+05} 
& \num{4.285e-01} & \num{1.3399363e-03} \\

[KWF] & \num[round-precision=2]{3.75000e-05} & \num[round-precision=2]{1.87500e-03} & \num[round-precision=2]{7.8125e-06} & \num[round-precision=2]{2.50000e-04} & \num{2.89500e-01} & \num{-1.7092e-01} & \num{4.9549e+00} 
& \num{3.580078e-01} & \num{1.1499965e-05} 
&\num{6.7855230e+01} & \num{1.8406e+05}
& \num{4.3375e-01} & \num{1.5198394e-03} \\
\hline


[FST] & \num[round-precision=2]{3.05900313e-07} & \num[round-precision=2]{1.26378917e-03} & \num[round-precision=2]{1.56200e-05} & \num[round-precision=2]{1.00000e-03} & \num{2.88195e-01} & \num{-1.70360e-01} & \num{4.96789e+00}
& \num{3.57328e-01} & \num{5.93000e-06} 
&\num{6.392584677e1} & \num{1.00057e+05} 
& \num{4.29109e-01} & \num{1.34009e-03} \\

[FST] & \num[round-precision=2]{2.67923617e-07} & \num[round-precision=2]{9.70524954e-04} & \num[round-precision=2]{1.56200e-05} & \num[round-precision=2]{1.00000e-03} & \num{2.88195e-01} & \num{-1.70340e-01} & \num{4.96829e+00}
& \num{3.57344e-01} & \num{5.92000e-06} 
&\num{6.493203607e1} & \num{1.00007e+05} 
& \num{4.29047e-01} & \num{1.33163e-03} \\

[FST] & \num[round-precision=2]{2.18676954e-07} & \num[round-precision=2]{8.51448139e-04} & \num[round-precision=2]{1.56200e-05} & \num[round-precision=2]{1.00000e-03} & \num{2.88194e-01} & \num{-1.70340e-01} & \num{4.96833e+00}
& \num{3.57359e-01} & \num{5.92000e-06}
&\num{6.436219052e1} & \num{9.99630e+04} 
& \num{4.28828e-01} & \num{1.32122e-03} \\
\hline

\hline 
\hline

[CF2] & \num[round-precision=2]{1.54e-5} & \num[round-precision=2]{8e-4} & \num[round-precision=2]{2e-4} &  \num[round-precision=2]{2e-4} & \num{2.9060e-1} & \num{-1.7e-1} & \num{4.96614} 
&\num{3.59e-1} &*\num{9.56058e-07}
&\num{1.57562e1} & \num{1.03379e4} 
& \num{4.698e-1} &  \num{3.88207e-3} \\

[CF2] & \num[round-precision=2]{7.93e-6} & \num[round-precision=2]{4e-4} & \num[round-precision=2]{1e-4} &  \num[round-precision=2]{1e-4} & \num{2.898e-1} & \num{-1.7e-1} 
& \num{4.99865} 
&\num{3.585e-1} &*\num{4.86374e-07}
&\num{1.74445e1} &\num{1.23205e4} 
& \num{4.737e-1} & \num{4.48394e-3} \\

[CF2] & \num[round-precision=2]{3.963e-6} & \num[round-precision=2]{2e-4} & \num[round-precision=2]{5e-5} &  \num[round-precision=2]{5e-5} & \num{2.8955e-1} & \num{-1.7e-1} & \num{4.96892} 
&\num{3.5845e-1} &*\num{2.26498e-07}
&\num{1.8153e1} &\num{1.2773e4} 
& \num{4.7575e-1} & \num{4.79279e-3} \\

\hline

[CF3] & \num[round-precision=2]{1.54e-5} & \num[round-precision=2]{8e-4} & \num[round-precision=2]{2e-4} & \num[round-precision=2]{2e-4} & \num{2.90e-1} & \num{-1.69e-1} & \num{4.97176} 
&\num{3.592e-1} &*\num{2.19345e-07}
&\num{7.29143} & \num{ 9.47679e2} 
& \num{4.762e-1} & \num{4.61781e-3} \\

[CF3] & \num[round-precision=2]{7.93e-6} & \num[round-precision=2]{4e-4} & \num[round-precision=2]{1e-4} &  \num[round-precision=2]{1e-4} & \num{2.893e-1} & \num{-1.7e-1} & \num{4.97407} 
&\num{3.582e-1} &*\num{4.86374e-07}
&\num{7.36221} & \num{9.76773e2} 
& \num{4.7880e-1} & \num{5.1582e-3} \\

[CF3] & \num[round-precision=2]{3.963e-6} & \num[round-precision=2]{2e-4} & \num[round-precision=2]{5e-5} &  \num[round-precision=2]{5e-5} & \num{2.8915e-1} & \num{-1.68e-1} & \num{4.97463} 
&\num{3.5815e-1} &*\num{2.50339e-07}
&\num{7.39539} & \num{9.86732e2}
& \num{4.8040e-1} & \num{5.40013e-3}\\

\hline

[Fr2] &\num[round-precision=2]{6.25e-05} & \num[round-precision=2]{6.25e-04} &\num[round-precision=2]{1.1094e-5}&\num[round-precision=2]{1.42e-3} &\num{2.89058166118421e-1} &\num{-1.71457312911184e-1} &\num{4.94775462582237} 
&\num{3.5754e-1}
&*\num{7.79704e-06} 
&\num{2.10829e1} &\num{1.72724e+4}  
&\num{4.7164e-1} &\num{4.39023e-3}
\\

[Fr2] &\num[round-precision=2]{4.6875e-05} & \num[round-precision=2]{4.6875e-04} &\num[round-precision=2]{7.8125e-6} &\num[round-precision=2]{1e-3} &\num{2.88905792861322e-1} &\num{-1.71103216500878e-1} &\num{4.95097266822703} 
&\num{3.5757e-1}
&*\num{5.71701e-06} 
&\num{2.18529e1} &\num{1.93103e+4}  
&\num{4.73359e-1} 
&\num{4.40132e-3}\\

[Fr2] &\num[round-precision=2]{3.125e-05} & \num[round-precision=2]{3.125e-04} &\num[round-precision=2]{5.5469e-6}&\num[round-precision=2]{7.1e-4} &\num{2.88786793075021e-1} &\num{-1.70913326463314e-1} &\num{4.95395937345425} 
&\num{3.57529e-1} 
&*\num{4.14201e-06} 
&\num{2.1539e1} &\num{1.91484e+4} 
&\num{4.72072e-1} &\num{4.45864e-3}\\
\end{tabular}
\end{adjustbox}

\caption{Quantitative simulation results obtained with the eight approaches and the first nine quantities of interest described in Section~\ref{sec.qoi}. A star in front of the {minimum distance} $d_{\rm cont}$ indicates that a no-penetration condition has become active, such that the distance is essentially determined by the {(weakly)} imposed gap {width $\epsilon$}. For the coarsest [W] result, the star indicates that the disk reached the bottom boundary. If no star is indicated, a contactless rebound has taken place. \label{tab:res_qi}
}
\end{table}
\end{landscape}

{The table is divided into three parts separated by double lines. We start with the simple ALE approach [W] {(with no-slip boundary conditions)}, followed by four numerically more advanced methods that also use no-slip conditions on the bottom, and finally three approaches in which a surface Darcy model was used to resolve the no-contact paradox.}

\paragraph{Pre-contact dynamics.}
{As expected, the first three quantities $t^*$, $v^*$ and $f^*$, which are recorded significantly before {the near-contact event}, are relatively close for all approaches on the finest meshes. Due to the coarse time step of $\delta t =10^{-2}$ used by [W] at this instant, the results differ slightly from the remaining approaches. For {the others}, the time $t^*$, {at which the distance} $r_s$ is reached, {lies within the interval} [$2.8819\cdot 10^{-1}$s, $2.8955\cdot 10^{-1}$s] {considering the finest computations of each group}. The deviations of order {$10^{-3}$\,s} {likely reflect differences in both spatial and temporal discretization}.  
With the exception of the results [W] and [CF3], the velocity $v^*$ at time $t^*$ is determined with similar accuracy, ranging between $-1.7092\cdot 10^{-1}$ and $-1.70\cdot 10^{-1}$ {$\frac{\rm m}{\rm s}$}. The small deviation in [CF3] is due to the fact that the no-slip condition on the interface is modified by the porous layer on the interface \eqref{eq:darcy-moving}-\eqref{eq:por:int}, which has a visible effect even before {the near-contact event}. Similarly, the drag $f^*$ takes values in the range $(4.9540, 4.9746) {\frac{\rm kg}{{\rm s}^2}}$ among all approaches except [W]. Among all approaches, [FST] uses the finest meshes and thus produces results that differ the least between the three mesh levels.}


\paragraph{{Near-contact and contact dynamics.}}
{Similarly, the time {of the first near-contact or contact event} $t_{\rm cont}$ agrees very well between the approaches [CF1], [CF2], [CF3], [Fr1], [Fr2], [KWF] and [FST], ranging between $3.5736\cdot 10^{-1}$ and $3.5875\cdot 10^{-1}$s on the finest mesh levels.} 

{Concerning the {minimum distance at this event}, $d_{\text{cont}}$, we note that for several approaches a no-penetration condition became active; {for [W] on the coarsest mesh the disk actually reached the bottom boundary.} These cases are marked by a star in front of the value of $d_{\text{cont}}$. In the former cases, the distance is mainly determined by the (weakly) imposed gap $\epsilon=ch$. For the Darcy-contact approaches, it is expected that the no-penetration condition becomes active due to the contact expected for the corresponding continuous model. The condition also becomes active for [Fr1], where the mesh resolution in the {near-contact} region is not accurate enough to prevent {spurious} contact (as the pressure is not resolved sufficiently). For [CF1], this is the case on the two coarsest meshes. {For [FST] and [KWF] on all mesh levels, and for [CF1] on the finest mesh level}, a contactless rebound is observed. While for [FST] and [CF1] the minimum distance agrees well and is around $6\,\mu$m, the distance predicted by {[KWF]} is significantly larger.}

{The maximum drag force $f_{\rm max}$ before {the first near-contact or contact event} is of particular importance, as it {characterizes the lubrication forces that prevent or permit contact}. Again, we observe good agreement between the no-slip approaches, in particular [FST] and [CF1] on the finest mesh level, where $f_{\rm max}\approx 6.4\cdot 10^1$. [Fr1] is also approaching a similar value, but from the three values provided one can conclude that the simulation is not yet fully converged. {Consistent with} the theoretical results in~\cite{ChampionFernandezetal2024}, the {maximum drag force} is much smaller for the Darcy-contact approaches. Here, the results of [CF2] and [Fr2] (using the Darcy model on the bottom $\Gamma_{\text{bo{t}}}$) show reasonable agreement with values $f_{\rm max}\approx 2\cdot 10^1$. We will investigate the behavior of the drag force plotted over the distance to $\Gamma_{\text{bo{t}}}$ in more detail in the following subsection. [CF3] shows a significantly smaller drag force. This is due to the fact that the surface Darcy model is considered on the moving interface $\Gamma_{t}^{\rm i}$, such that fluid can escape through $\Gamma_{t}^{\rm i}$ itself. This results in a significantly lower pressure on $\Gamma_{t}^{\rm i}$, which is the dominant contribution to the drag force.}

{The maximum pressure point value at the central bottom $\max_t p_{bc}$ is more difficult to approximate, since (a) it would get singular in case of full solid-solid contact and (b) the computation of point values is known to be less stable compared to integrals such as the drag force $f_{\rm max}$, which is essentially the pressure force integrated over $\Gamma_t^{\rm i}$. We observe that all no-slip approaches show a similar order of magnitude around $10^5$, while quantitatively [Fr1] and [KWF] differ by almost a factor of 2 from [CF1] and [FST] on the finest mesh level. On the other hand, the results [CF1] and [FST] show excellent agreement, with values of $9.6667\cdot 10^{4}$ resp.\,$9.9963\cdot 10^{4} \frac{{\rm kg}}{{\rm m s}^2}$ on the finest meshes. For the Darcy contact approaches the maximum pressure is significantly smaller, by around one order of magnitude for [CF2] and [Fr2] and by roughly another magnitude for [CF3]. The deviation of the latter can again be explained by the fact that the Darcy model is placed on the moving  interface.}

\paragraph{{Post-event dynamics.}}
{The quantities after {the near-contact or contact event} are the most difficult to predict, as they are very sensitive to the {preceding dynamics}. For this reason, the approach [W] {deviates substantially from} the other no-slip approaches. The rebound height $d_{\text{jump}}$ of [CF1], [Fr1], [KWF] and [FST] ranges between $9.3074\cdot 10^{-4}$ and $1.5198\cdot 10^{-3}$m. As expected, we obtain a much larger rebound for the Darcy-contact approaches (around four times larger), ranging from $4.4586\cdot 10^{-3}$ to $5.4001\cdot 10^{-3}$m on the finest mesh levels. We will give further insights into these {near-contact, contact,} and {post-event} dynamics in the next subsection. Accordingly, the time {to reach the maximum rebound height} is longer for the Darcy-contact approaches than for the no-slip model.}


\subsection{Graphical visualizations}


In this subsection, we compare the results of the different groups visually by plotting different quantities over time. Compared to Table~\ref{tab:res_qi}, we only report the {\it finest} result of each approach.

\paragraph{Minimum distance over time.}

\begin{figure}[h!]
\centering
\includegraphics[width=0.70\textwidth]{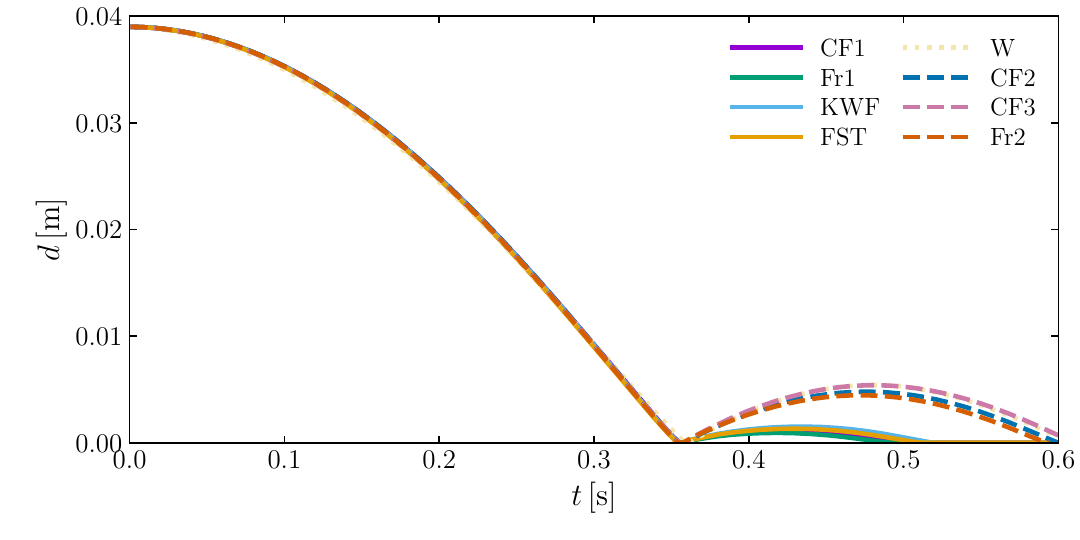}
\caption{Minimum gap between the disk and the lower wall over the full time
interval for the eight numerical approaches.}
\label{fig:mindist}
\end{figure}

Figure~\ref{fig:mindist} shows that all methods predict nearly the same falling
phase. The differences appear only in the near-contact regime and in the first
rebound. [CF1], [Fr1], [KWF] and [FST] form a close group and predict a relatively low
first rebound. In contrast, [W], [CF2], [CF3] and [Fr2] predict a visibly larger rebound.
\begin{figure}[h!]
\centering
\includegraphics[width=0.33\textwidth]{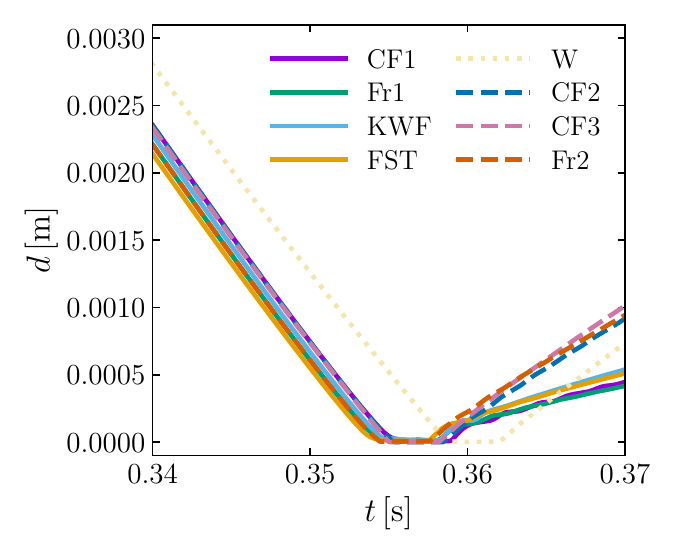}\hfill
\includegraphics[width=0.33\textwidth]{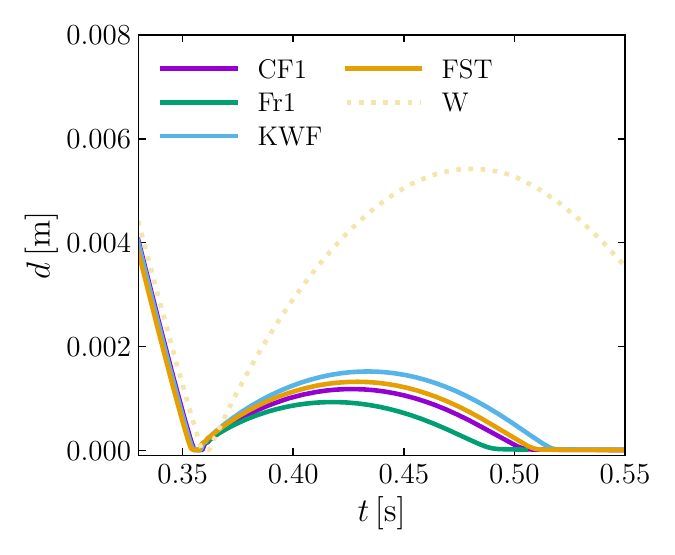}\hfill
\includegraphics[width=0.33\textwidth]{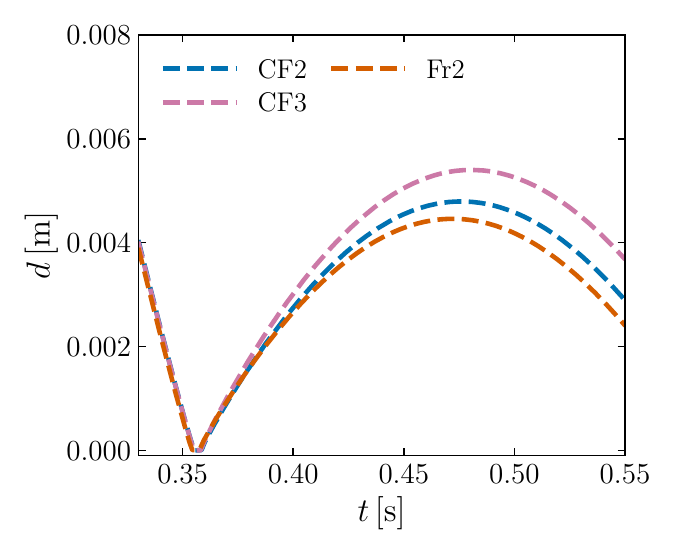}

\caption{Minimum gap between the disk and the lower wall. Left: zoom into the
near-contact regime. Middle and right: first rebound after the near-contact event.}
\label{fig:mindist_zoom}
\end{figure}
For [CF2], [CF3] and [Fr2], the larger rebound is associated with the Darcy-contact model,
which weakens the lubrication effect near the lower wall, as fluid can now escape through the porous layer. As a consequence, the disk is slowed down less, which leads to a larger compression of the disk and hence a larger amount of elastic energy can be transferred into kinetic energy for the rebound. The result [W] also
shows a larger rebound, but for a different numerical reason. The vanilla full
ALE method accumulates the mesh deformation from the initial configuration.
Consequently, the mesh elements below the disk become severely compressed and the
thin fluid layer is not resolved sufficiently accurately.

Nevertheless, we consider [W] a valuable result: it demonstrates that a
vanilla full ALE formulation can pass through an almost-contact event and
produce a rebound. However, the comparison also shows that vanilla ALE alone
is not sufficient to obtain a quantitatively reliable result. The thin fluid
layer is not resolved accurately enough, leading to an inaccurate rebound.
Obtaining a reliable solution therefore requires additional numerical effort to resolve the thin fluid layer accurately.


Figure~\ref{fig:mindist_zoom} provides a closer view of the near-contact event
and the first rebound. In particular, [W] shows a visibly delayed time of minimal distance $t_{\rm cont}$
compared with the other approaches.
The zoom also shows that the rebound is not instantaneous. After entering the
near-contact regime, the disk remains close to the lower wall over a short but
finite time interval, during which the fluid pressure deforms the disk and
reverses its motion.


\paragraph{Energies.}

\begin{figure}[h!]
\centering
\hspace*{4mm}FST\par\vspace{0.0pt}
\includegraphics[width=0.60\textwidth]{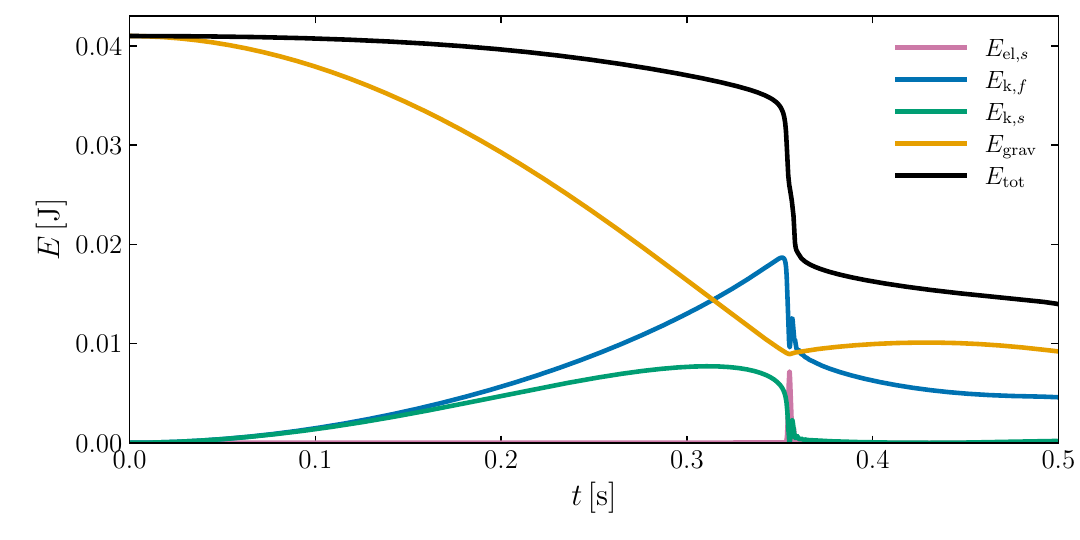}
\caption{Energy contributions for the [FST] simulation over the full time
interval. Gravitational energy is converted into kinetic energy during the
fall. During the finite near-contact interval, part of the kinetic energy is
temporarily stored as elastic energy of the solid and is subsequently released
during rebound. The total mechanical energy decreases throughout the
computation.}
\label{fig:energy_FST}
\end{figure}

The energy plots in Figures~\ref{fig:energy_FST} and
\ref{fig:energy_zoom_all} illustrate the transfer of energy during
fall, the near-contact or contact event, and rebound. Initially, most mechanical energy is stored as
gravitational potential energy $E_{\mathrm{grav}}$. As the disk falls,
$E_{\mathrm{grav}}$ decreases and is converted into kinetic energy of the fluid
and the solid, $E_{\mathrm{k},f}$ and $E_{\mathrm{k},s}$.

During the short near-contact interval, the reversal of motion is not instantaneous: part
of the kinetic energy is temporarily stored as elastic energy of the solid,
which appears as a sharp peak of $E_{\mathrm{el},s}$. Afterwards, this elastic
energy is released and converted back mainly into solid kinetic energy, leading
to the rebound. The [FST] result in Figure~\ref{fig:energy_FST} illustrates this
mechanism {over the whole time interval}.

{Owing to \eqref{eq:energy-balance}, the total mechanical energy
$
E_{\mathrm{tot}}
=
E_{\mathrm{grav}}
+
E_{\mathrm{k,f}}
+
E_{\mathrm{k,s}}
+
E_{\mathrm{el,s}}
$
is non-increasing at the continuous level. This fundamental  property  should be    preserved at the discrete level by the numerical schemes.} This is observed for [FST], [CF1], [CF2], [CF3] and, up to very minor
numerical imperfections, also for [KWF]; see
Figure~\ref{fig:energy_zoom_all}. [Fr1] and [Fr2] show a small local
non-monotonicity of $E_{\mathrm{tot}}$ near the event, where the energy balance seems to be temporarily perturbed by the contact force. Note, however, that after the event, the total energy takes very similar values to those in the other approaches. The strongest violation is
observed for [W], {where the fluid kinetic energy $E_{\mathrm{k,f}}$ increases in an unphysical way before {the near-contact event}, due to the extreme mesh deformations of the ALE method, which enter through the determinant $J^{\rm f}$ of the transformation gradient in the evaluation of $E_{\mathrm{k,f}}$ in the reference frame
$
E_{\rm k, f} = \int_{\Omega_0^{\rm f}} \frac12 J^{\rm f}\rho^{\rm f} |\vec{v}^{\rm f}|^2.
$}
Comparing the approaches using no-slip conditions for the bottom [CF1], [KWF], [Fr1], [FST] and the Darcy-contact approaches, we observe that the elastic energy is slightly higher for the latter group. This additional elastic energy is then transferred into the kinetic energies after the event. In particular, the solid kinetic energy after the event is larger, as seen from the green lines at the final time shown $t=0.38$\,s, which leads to larger rebounds.

\newcommand{\energyplot}[2]{%
\begin{minipage}[t]{0.49\textwidth}
\centering
\hspace*{4mm}#1\par\vspace{-1pt}
\includegraphics[width=\linewidth]{#2}
\end{minipage}}

\begin{figure}[p]
\centering
\energyplot{CF1}{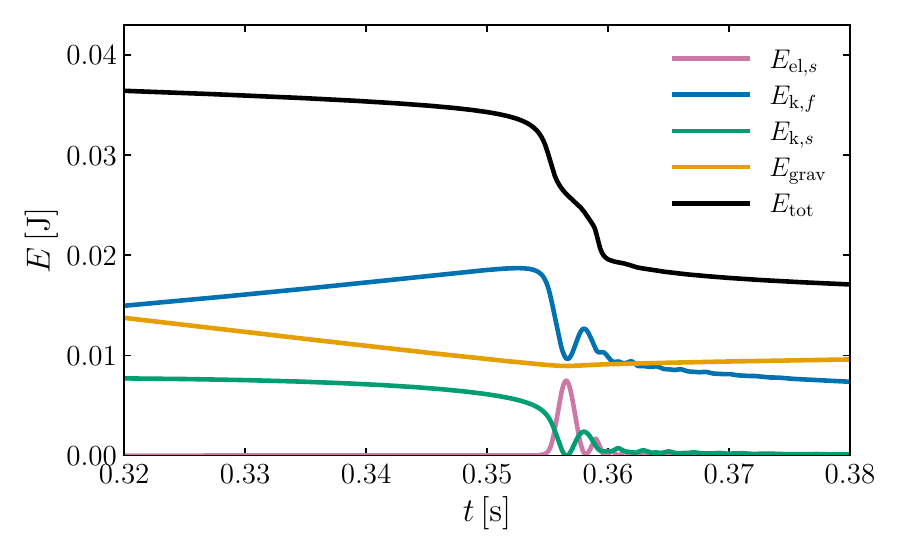}\hfill
\energyplot{Fr1}{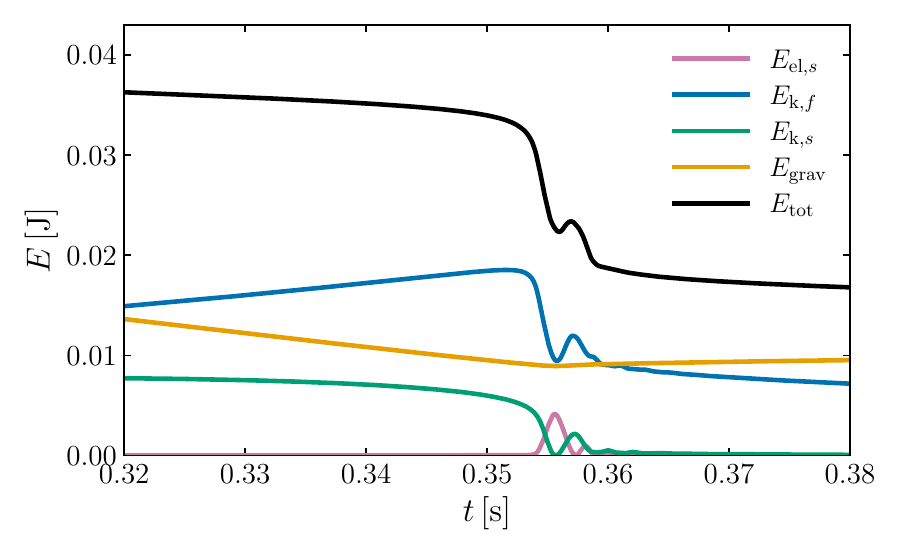}\\[0.4em]
\energyplot{KWF}{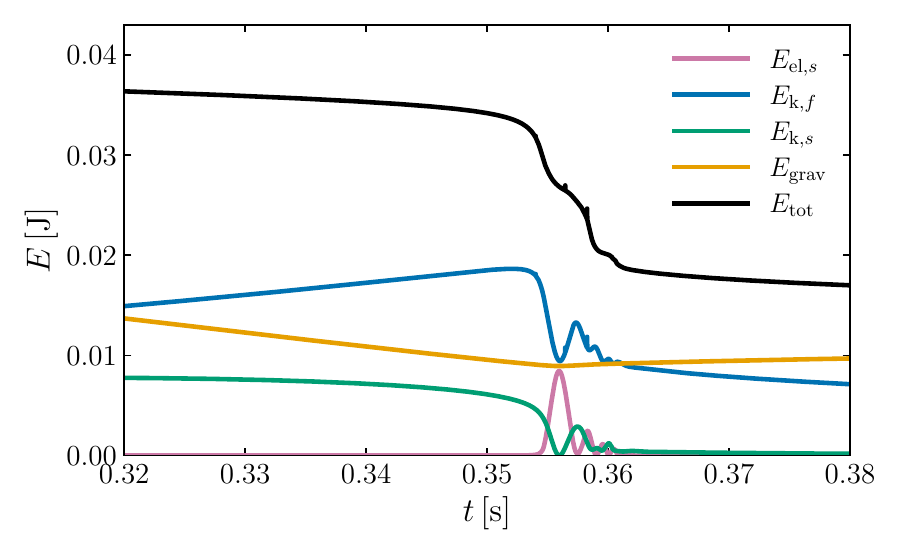}\hfill
\energyplot{FST}{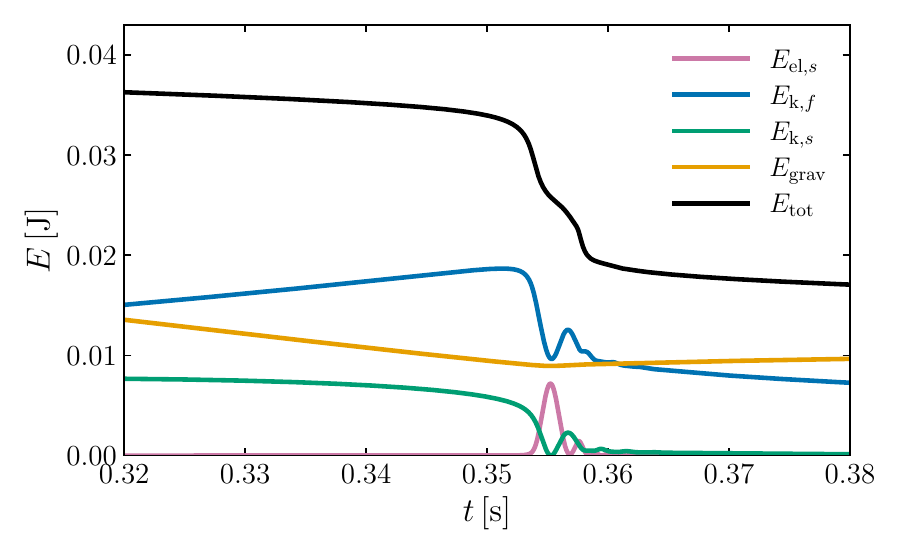}\\[0.4em]
\energyplot{W}{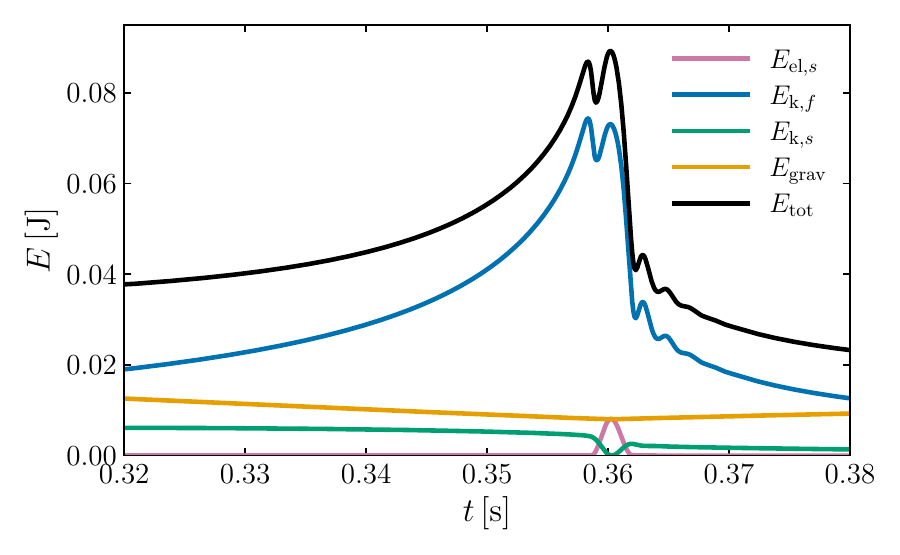}\hfill
\energyplot{CF2}{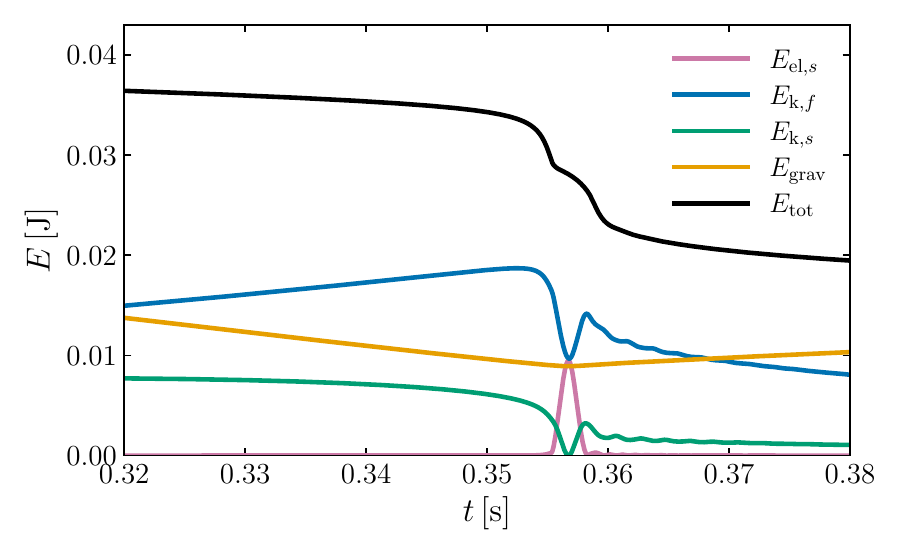}\\[0.4em]
\energyplot{CF3}{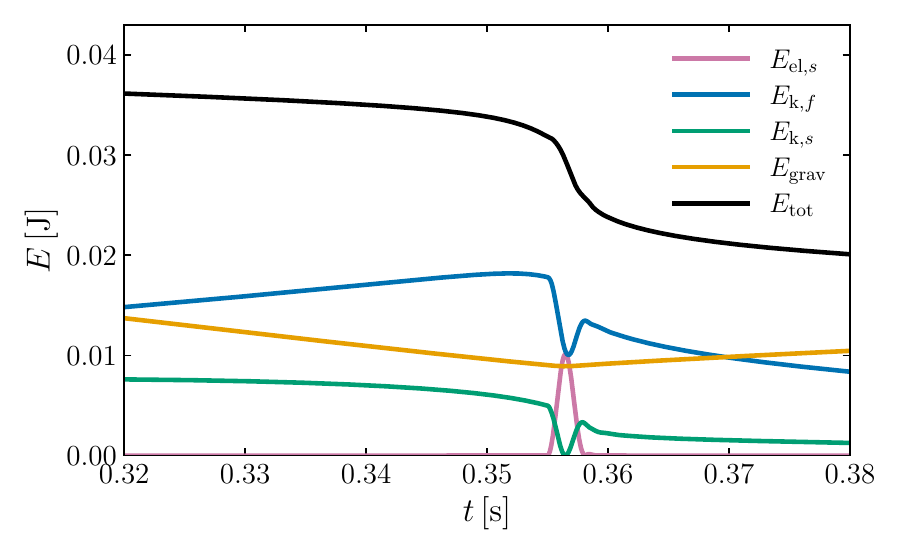}\hfill
\energyplot{Fr2}{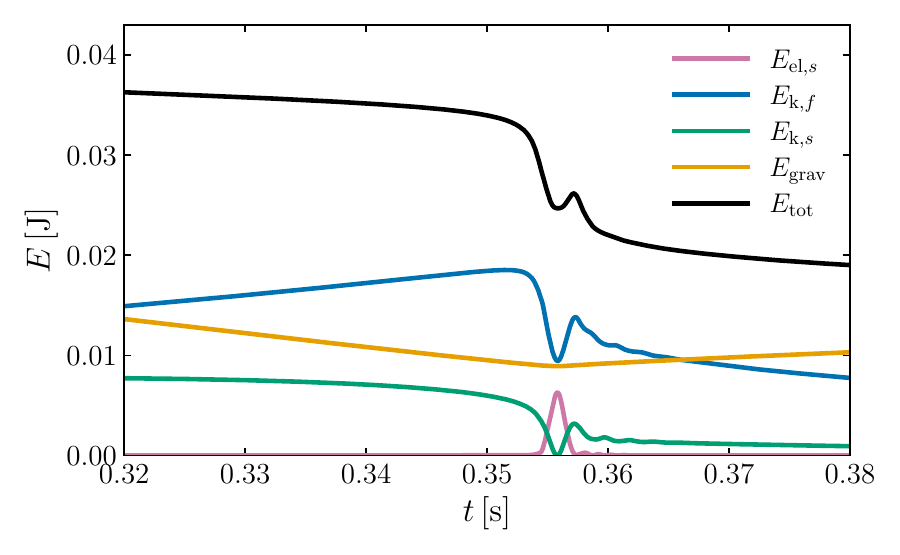}\\[0.4em]
\caption{Energy contributions in the near-contact regime for the different
approaches. All methods show the same basic mechanism of energy transfer, but
the behavior of the total energy differs. [CF1], [CF2], [CF3] and [FST] show a clean monotone decay of
$E_{\mathrm{tot}}$, [Fr1], [Fr2], and [KWF] show small
local violations of monotonicity near the event, and [W] shows the strongest
nonphysical increase of total energy. Note that the [W] panel uses a different
vertical scale.}
\label{fig:energy_zoom_all}
\end{figure}

\begin{figure}[h!]
\centering
\hspace*{-53mm}CF1 vs.\ CF2\par\vspace{-1pt}
\includegraphics[width=0.86\textwidth]{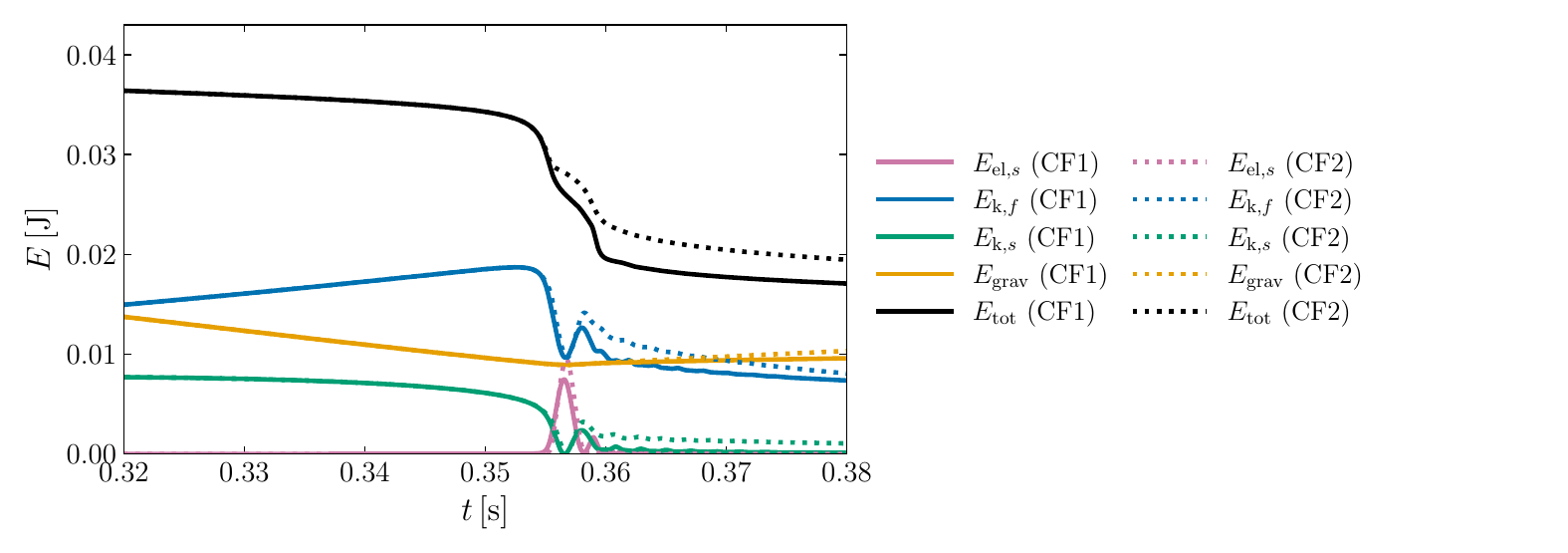}
\caption{Direct comparison of the energy contributions for the no-slip case
[CF1] (solid lines) and the Darcy-like contact case [CF2] (dotted lines). Both
approaches show the same qualitative energy transfer mechanism. Near and after
the near-contact or contact event, [CF2] exhibits slightly larger elastic and
post-event kinetic energies, consistent with the stronger rebound.}
\label{fig:energy_CF1_CF2_compare}
\end{figure}

A direct comparison of the CutFEM no-slip case [CF1] and the Darcy-like
contact case [CF2] is shown in Figure~\ref{fig:energy_CF1_CF2_compare}. The
overall energy evolution is very similar in both simulations. Before the near-contact event,
the decrease of gravitational energy is converted into fluid and solid kinetic
energy in essentially the same way. Near the event, the Darcy-like model [CF2]
shows a slightly larger elastic response and slightly larger post-event
kinetic energies, in particular of the solid, which is consistent with the
larger rebound observed in the distance plot in Figure~\ref{fig:mindist}.

\paragraph{Drag force.}

\begin{figure}[!h]
\centering
\includegraphics[width=0.49\textwidth]{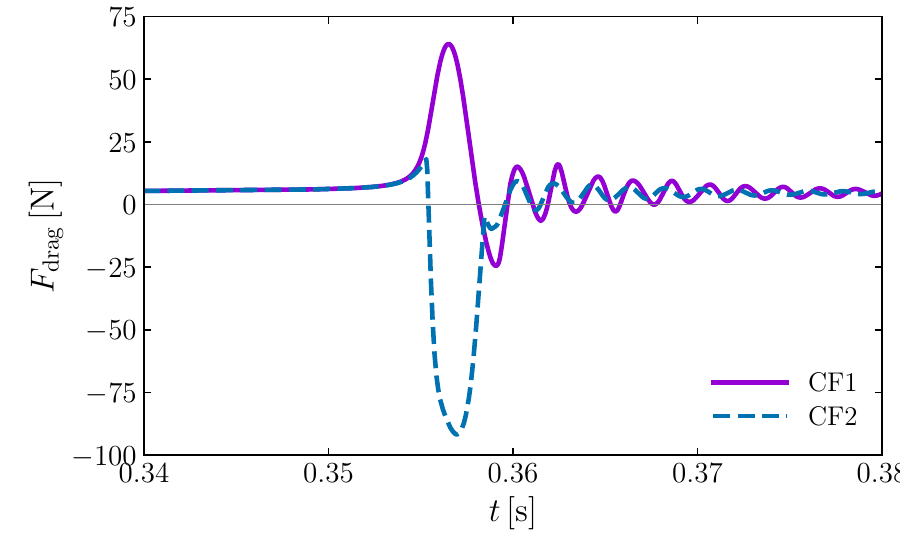}\hfill
\includegraphics[width=0.49\textwidth]{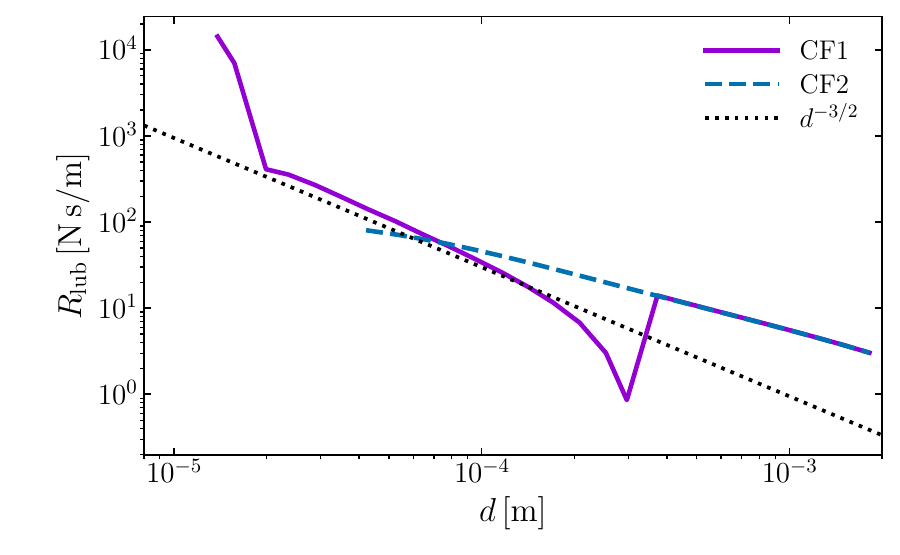}
\caption{Comparison of the no-slip result [CF1] and the Darcy-like result
[CF2]. Left: drag force during the near-contact event. Right: effective
resistance $R_{\mathrm{lub}}=F_{\mathrm{excess}}/V_n$ during approach,
plotted against the minimum gap $d$. The dotted black line shows the reference
scaling $R_{\mathrm{lub}}\sim d^{-3/2}$.}
\label{fig:drag}
\end{figure}

Figure~\ref{fig:drag} compares the drag force and the effective lubrication
resistance for the no-slip result [CF1] and the Darcy-like result [CF2].
As the disk approaches the lower wall, the drag grows rapidly due to the thin
fluid layer. During rebound, the force changes sign and subsequently relaxes
towards its far-field value.

To account for the changing velocity during the fall, we consider
\[
R_{\mathrm{lub}}=\frac{F_{\mathrm{excess}}}{V_n},
\]
where $F_{\mathrm{excess}}$ is the drag force after subtracting a far-field
contribution and $V_n$ is the normal approach speed. For the
two-dimensional no-slip setting, the expected scaling for a rigid disk not
changing shape is
$$
R_{\mathrm{lub}}\sim d^{-3/2}
\quad\text{as } d\to0.
$$
However, as is shown in~\cite{GraSchSouTum22}, for a contactless rebound it is
{crucial that the disk is changing shape}, such that the effective
lubrication resistance changes its scaling to $d^{-\alpha}$ for $\alpha$
strictly larger than $\frac{3}{2}$. Indeed, this shape-induced asymmetry in the
fluid resistance is what allows the elastic energy to be transferred back into
kinetic energy.

The direct comparison in Figure~\ref{fig:drag} shows that [CF1] and [CF2]
behave similarly during most of the fall but deviate in the near-contact regime, where
the modified lower-wall condition in [CF2] changes the lubrication resistance
and permits contact. {In the right plot, we observe that for [CF1] the drag force follows the ${\cal O}(d^{-3/2})$ line up to a very small distance and then increases significantly again, due to the modified shape. Such an effect is not observed for the Darcy-type approach [CF2], where the drag force remains bounded.}

\begin{figure}[!h]
\centering
\includegraphics[width=0.49\textwidth]{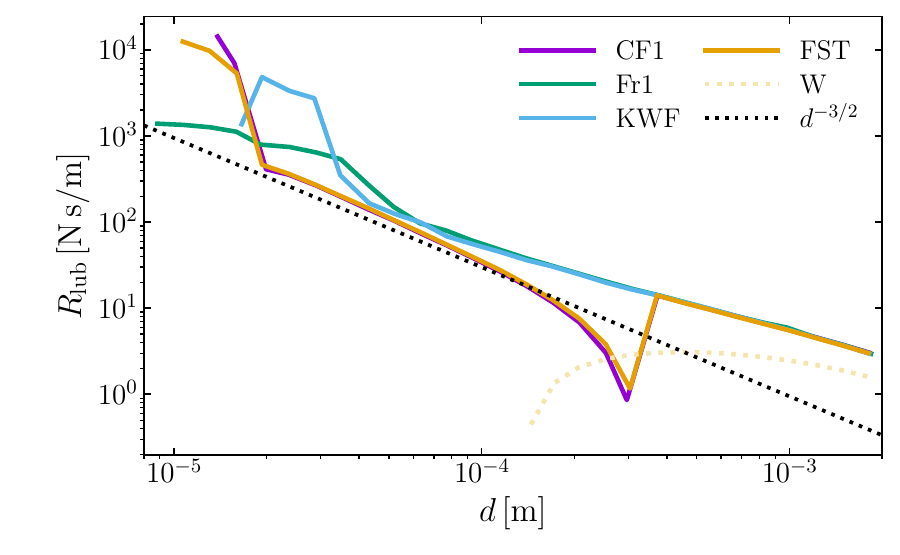}\hfill
\includegraphics[width=0.49\textwidth]{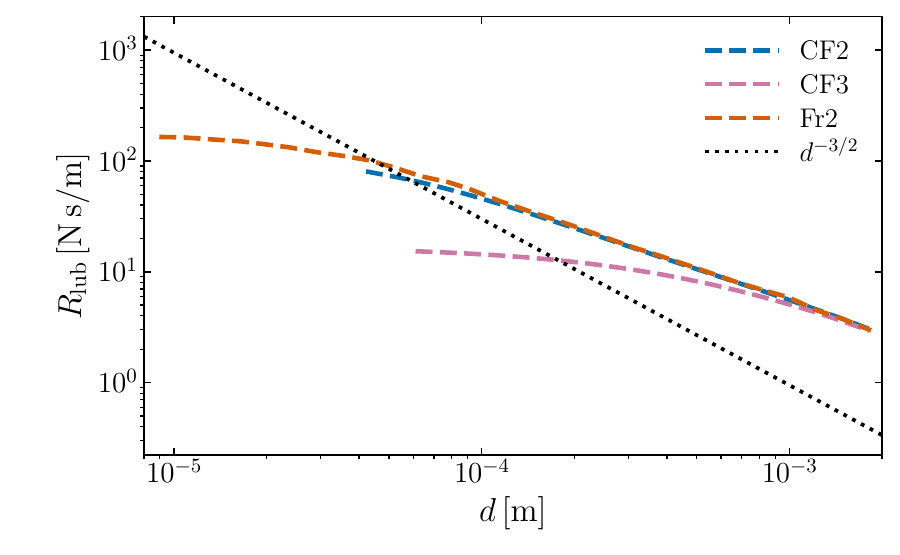}
\caption{Effective resistance during approach, plotted against the minimum
gap $d$. Left: no-slip approaches [CF1], [Fr1], [KWF], [FST] and [W].
Right: Darcy-contact approaches [CF2], [CF3] and [Fr2]. The dotted black lines
show the reference behavior proportional to $d^{-3/2}$.}
\label{fig:drag_resistance}
\end{figure}

The no-slip results [CF1], [Fr1], [KWF] and [FST] in the left panel of
Figure~\ref{fig:drag_resistance} show the same overall increase of the
lubrication resistance as the gap decreases. Over an intermediate near-contact
range, their behavior is broadly compatible with the reference $d^{-3/2}$
trend. Furthermore, [CF1], [KWF] and [FST] show a clear increase of the
effective exponent near the point of closest approach, reflecting the flattening of the disk.
Nevertheless, at the smallest gaps, the estimate becomes sensitive because the
normal approach velocity tends to zero and the quotient
$F_{\mathrm{excess}}/V_n$ is therefore strongly affected by small numerical
errors.

The approach [W] deviates substantially from this trend. The vanilla full ALE
formulation does not resolve the lubrication regime sufficiently accurately,
since the mesh elements in the thin gap become severely compressed. As a
result, the expected increase of the lubrication resistance is not recovered.
This is consistent with the larger rebound observed in the minimum-distance
comparison.

The Darcy-contact results [CF2], [CF3] and [Fr2] are shown in the right panel
of Figure~\ref{fig:drag_resistance}. 
As has been shown
in~\cite{ChampionFernandezetal2024}, the drag force remains bounded as
$d\to0$, which enables contact.

\FloatBarrier
\section{Conclusion}
\label{sec:conclusions}

{We have introduced a benchmark problem for {FSI with near-contact, possible contact, and rebound} with the objective of keeping it simple enough to be reproducible by different numerical codes. The quantities of interest defined {before the near-contact or contact event} can even be reproduced by pure FSI codes, which enables verification of {the FSI component independently of the contact handling}. 
}

In our numerical study, we have investigated three models for FSI with contact: {one using} a no-slip condition on the bottom wall, and {two using} a Darcy description of the surface roughness, imposed either on the bottom wall or on the moving interface. {As expected, the models} agree very well {before the near-contact or contact event}, while {after this event} significantly different dynamics are observed; in particular, {the rebound is much higher for the Darcy model}. Moreover, we have also numerically {investigated the effective lubrication resistance as a function of the minimum gap}. We observe {the contactless-rebound mechanism described in~\cite{GraSchSouTum22} for the no-slip model, while the Darcy-contact model enables  contact through a bounded resistance, as predicted in~\cite{ChampionFernandezetal2024} in the rigid-body  case}.

{Furthermore, we have compared very different numerical methodologies for both FSI and contact simulation. {We observe that even though a rebound can be obtained using a simple vanilla ALE method with no-slip conditions, {refinement strategies and specialized mesh treatment substantially improve the quantitative accuracy of the near-contact and rebound dynamics.} Among the methods with {dedicated treatment of the thin-gap or contact regime},} the quantitative differences {within each model class are substantially smaller than those observed for the vanilla ALE result}.
}

The agreement of the selected quantities of interests (QoIs) is satisfactory with regard to the expected convergence behavior of fluid and solid velocity as well as the fluid pressure.
Further, we observe qualitatively in all approaches how kinetic energy is transferred into elastic energy and then, with substantial losses due to dissipation, is returned to kinetic energy. Moreover, {the mechanisms described for contactless rebound in~\cite{GraSchSouTum22} and for contact with Darcy seepage in~\cite{ChampionFernandezetal2024} are supported by} the analysis of the fluid forces on the elastic solid, in particular the dynamics of its effective lubrication resistance, which changes drastically during the {near-contact or contact phase} and also differs significantly between the two approaches.

{In summary, following our experience in producing the numerical results, we suggest the following hierarchy of checkpoints or quality signatures as a guideline for future repetitions of the benchmark:}
\begin{enumerate}[itemsep=4pt, topsep=5pt, parsep=0pt]
    \item Recover the pre-contact QoI parameters $t^*$, $v^*$ and, in particular, $f^*$, as the correct resolution of the fluid forces on the elastic object is obviously essential here;
    \item Compare the QoIs at the {near-contact or contact event} and {there}after {within each model class};
    \item Analyze the energy transfer from kinetic to elastic energy and back in relation to the total energy, which should generally be {non-increasing};
    \item Analyze the total drag force {and the effective lubrication resistance as functions of the minimum gap}.
\end{enumerate}

{To determine whether {and under which conditions} the contactless rebound that we observe using no-slip conditions is physical, or whether the Darcy-contact model represents the correct physics more accurately, we plan to compare the models with experimental results in the future.
} {This is our suggestion for the next step towards {physically appropriate modelling of fluid effects on elastic rebound}. The analysis of the fluid impact on rebounds as such remains {an important open problem} in fluid--structure interaction in computational mathematics and beyond.}





\appendix 

\section{Sketch of the proof of the energy balance\label{app:proof}}

Relation \eqref{eq:energy-balance} is obtained by testing the fluid momentum
equation \eqref{eq:fluid} with $\vec{v}^{\rm f}$, the solid momentum equation
\eqref{eq:solid} with $\partial_t \vec{u}^{\rm s}$ and the Darcy equation
\eqref{eq:darcy}, respectively \eqref{eq:darcy-moving}, with $P_{\rm l}$, and by
summing the resulting relations. The inertial terms give the time derivative of
the two kinetic energies, the fluid one after applying the Reynolds transport
theorem in $\Omega^{\rm f}_t$ and using $\div\,\vec{v}^{\rm f}=0$, while the elastic energy results  from \eqref{eq:hyper}, which yields
$\int_{\Omega^{\rm s}_0}\mathbb{P}:\nabla\partial_t\vec{u}^{\rm s} = \frac{{\rm d}}{{\rm
d}t}\int_{\Omega^{\rm s}_0}W^{\rm s}(\mathbb{E})$. The gravitational power of the solid is
directly the derivative of the corresponding potential energy, the reference
configuration being fixed. That of the fluid, on the contrary, has to be tested
with the velocity in the current configuration; since $\vec{g}$ is constant and
$\div\,\vec{v}^{\rm f} = 0$, we have
$\vec{g}\cdot\vec{v}^{\rm f} = \div\bigl((\vec{g}\cdot\vec{I})\,\vec{v}^{\rm f}\bigr)$, so
that, by using that $\vec{v}^{\rm f}\cdot\vec{n}=0$ on $\partial\Omega$, the kinematic condition in
\eqref{eq:fsi-coupling} and 
the Reynolds transport formula in $\Omega^{\rm s}_t$, we get  
$$
  \varrho^{\rm f}\int_{\Omega^{\rm f}_t} \vec{g}\cdot\vec{v}^{\rm f}
  = \varrho^{\rm f}\int_{\partial\Omega^{\rm f}_t} \bigl(\vec{g}\cdot\vec{I}\bigr)\,
    \vec{v}^{\rm f}\cdot\vec{n}
  = \varrho^{\rm f}\int_{\Gamma^{\rm i}_t } \bigl(\vec{g}\cdot\vec{I}\bigr)\,
    \vec{v}^{\rm s}\cdot\vec{n}
  = - \frac{{\rm d}}{{\rm d}t}\left( \varrho^{\rm f}\int_{\Omega^{\rm s}_t}
    \vec{g}\cdot\vec{I}\right).
$$
 The
gravitational power of the fluid is hence carried entirely by the fluid-solid
interface $\Gamma_t^{\rm i}$, which is why the associated potential energy appears in
\eqref{eq:total-energy} as an integral over the current solid configuration.

It remains to
collect the interface and boundary terms. Owing to the coupling conditions \eqref{eq:fsi-coupling} and \eqref{eq:mech-inconsistency}, respectively
\eqref{eq:por:int}; the interface contributions  produced by the   momentum equations 
yield the term $\int_{\Gamma^{\rm i}_0}\lambda\,\partial_t\vec{u}^{\rm s}\cdot\vec{n}_{\rm
b}$ to the
balance. Owing to the complementarity conditions \eqref{eq:contact_ineq2}, this term vanishes:  $\lambda=0$ on the
inactive set, whereas $\vec{u}^{\rm s}\cdot\vec{n}_{\rm b}=g_\epsilon$, and hence
$\partial_t(\vec{u}^{\rm s}\cdot\vec{n}_{\rm b})=0$, on the active set. The work of the
contact forces is therefore zero, so that \eqref{eq:energy-balance} holds
unchanged in the presence of the contact model.
On the other hand, the terms involving $P_{\rm l}$ (either on $\Gamma_{0}^{\rm i}$ for Model II or on $\Gamma_t^{\rm i}$ for Model III) combine with the Darcy equation
into the tangential dissipation of \eqref{eq:interface-dissipation}; and the
Beavers-Joseph-Saffman friction and the resistance opposed to the normal
seepage are non-negative. Finally, as regards the convective term, it yields a
  boundary   contribution
($\frac{\varrho^{\rm f}}{2}\int_{\Gamma_{\rm bot}}(\vec{v}^{\rm f}\cdot\vec{n})\,|\vec{v}^{\rm f}|^2$ for Model II or $\frac{\varrho^{\rm f}}{2}\int_{\Gamma_t^{\rm i}}(\vec{v}^{\rm f}\cdot\vec{n})\,|\vec{v}^{\rm f}|^2$ for Model III) of indefinite
sign. This with the back-flow stabilization results in 
the last terms of \eqref{eq:interface-dissipation}, which are non-negative. Hence
$\mathcal{D}\geq 0$.

\section{Description of the numerical approaches}
\label{sec:appendix}


\subsubsection*{[W] Wick}
A simple plain vanilla ``full ALE'' method is used, see Section~\ref{subsec.ALE}, where
the ALE mapping is realized with a nonlinear harmonic mesh motion model, i.e., $-\nabla\cdot (\frac{\alpha_u}{J}{\boldsymbol \nabla}u_f)$ with the choice $\alpha_u=10^{-1}$. Here, $u_f$ is the fluid displacement to move the fluid mesh in the ALE context. {A small} parameter $\alpha_u$ avoids back-coupling from the artificial fluid displacements into the physical solid displacements and {is needed} due to a variational-monolithic setup of our FSI system and globally defined velocity $v$ and displacement field $u$. 
Then, we have


\begin{problem}[Variational Monolithic ALE-FSI Problem]\label{prob:wick-full-ale}
     For almost all $t \in I$, find $\bmh{v} \in \bmh{\mathcal{V}}^0$, $\bmh{u} \in\bmh{\mathcal{V}}^1$ and $\widehat{p} \in L^2(\widehat{\Omega})$ satisfying $\bmh{v}(0) = \bmh{v}^0$ and $\bmh{u}(0) = \bmh{u}^0$ such that: 
    {\begin{equation}
        \begin{aligned}
            \label{ALE-FSI-prob}
            &\begin{alignedat}{3}
            \int_{\Omega^{\rm f}_0} \varrho^{\rm f} J^{\rm f} &\left(\pder{\vec{v}}{t} + {({\boldsymbol \nabla}\vec{v})}(\mathbb{F}^{\rm f})^{-1}(\vec{v} - \vec{w}^{\rm f})\right)\cdot{\boldsymbol \phi}
   + \int_{\Omega^{\rm f}_0} J^{\rm f} \mathbb{T}^{\rm f} (\mathbb{F}^{\rm f})^{-\rm T} : {\boldsymbol \nabla}{\boldsymbol \phi}  +\int_{\Omega^{\rm s}_0} \varrho^{\rm s} \pder{\vec{v}}{t}\cdot{\boldsymbol \phi} \\ 
  &
   + \int_{\Omega^{\rm s}_0} \mathbb{P} : {\boldsymbol \nabla}{\boldsymbol \phi}
+\int_{\Omega^{\rm f}_0} J^{\rm f} (\mathbb{F}^{\rm f})^{-1}\vec{v} \cdot {\boldsymbol \nabla} q
   +\sprods{\partial_t \vec{u}^{\rm s} - \vec{v}}{\psi} \\
          &\qquad \quad+ \sprodf{\alpha J^{-1} {\nabla} \vec{u}^{\rm f}}{{\nabla} \vec{\psi}} - \left< \alpha {J}^{-1} {\nabla} \vec{u}^{\rm f} \vec{n}, \psi \right>_{\Gamma_{\text{int}}}
                =\int_{\Omega^{\rm f}_0} \varrho^{\rm f} J^{\rm f} \vec{g}\cdot{\boldsymbol \phi} + \int_{\Omega^{\rm s}_0} \varrho^{\rm s} \vec{g}\cdot{\boldsymbol \phi}
            \end{alignedat}
        \end{aligned}
    \end{equation}}
    for all test functions $\phi,q,\psi$.
    Here $\vec{u}$ is a global displacement consisting of $\vec{u}^{\rm f}$ resp.\,$\vec{u}^{\rm s}$ in the subdomains and similarly $\psi$ is a global test function. 
\end{problem}
For details, we refer to~\cite{Wi11,Wi13_fsi_with_deal,Stahlhuth_Hergl_Wick_2026}. The discretization is fully implicit and
uses the Rothe method with finite differences in time and finite elements in space.

\paragraph{Temporal discretization} More specifically, a One-Step-$\theta$ scheme is employed for $\theta=1.0$ (first order implicit backward Euler). Therein, for a higher temporal solution during contact different time step sizes are employed. 
From $0\leq t\leq 0.32$ and $0.5\leq t\leq 0.6$, the time step size $\delta t=10^{-2}$ is used; between $0.32\leq t\leq 0.5$, the time step is $10^{-4}$. We notice that the numerical solvers work robustly with these choices and additional adaptive schemes are not employed. In our temporal refinement studies, we halve the time step sizes to $\delta t=5\cdot 10^{-3}$, resp. $\delta t=5\cdot 10^{-5}$, and on the finest time step level, we have $\delta t=2.5\cdot 10^{-3}$, resp. $\delta t=2.5\cdot 10^{-5}$.

\paragraph{Spatial discretization} 
The spatial discretization is by a conforming finite element Galerkin
method on quadrilateral meshes with the inf-sup-stable $(v,u,p)\in Q_c^2/Q_c^2/P_{dc}^1$ element. The initial mesh corresponds to Level 0 and is refined uniformly once
resulting in 336 elements and 
6644 (2818+2818+1008) degrees of freedom (DoFs) and used for our first computation. Next, for spatial computational convergence studies, the initial mesh is also twice and three times uniformely refined, resulting into 
1344 elements and 26052 (11010+11010+4032) DoFs and 
5376 elements and 103172 (43522+43522+16128) DoFs, respectively. As a side notice, no spatial adaptivity or local mesh refinement was used.

\paragraph{Numerical solution}
For the resulting nonlinear algebraic problems in
each time step a Newton iteration with backtracking line search is used. The linear subproblems are solved in fully coupled form by a direct solver (UMFPACK). 

\paragraph{Open-source programming codes}
More details and the
implementation are described
  in the FSI template~\cite{Wi13_fsi_with_deal}, which is based on the open-source software library deal.II~\cite{dealII96}. The codes are available open-source on github, see~\url{https://github.com/tommeswick/fsi}
  and \url{https://github.com/FlorianStahlhuth/FSI-FallingBall} ~\cite{Stahlhuth_Hergl_Wick_2026}.

\begin{figure}[!h]
    \centering
    \includegraphics[width=0.49\linewidth]{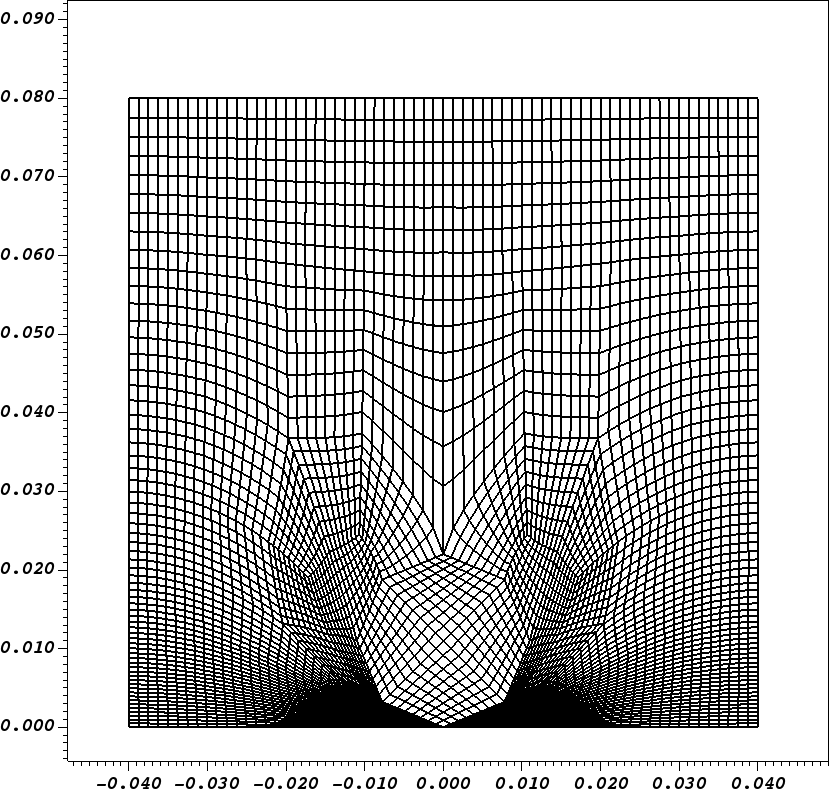}
    \includegraphics[width=0.47\linewidth]{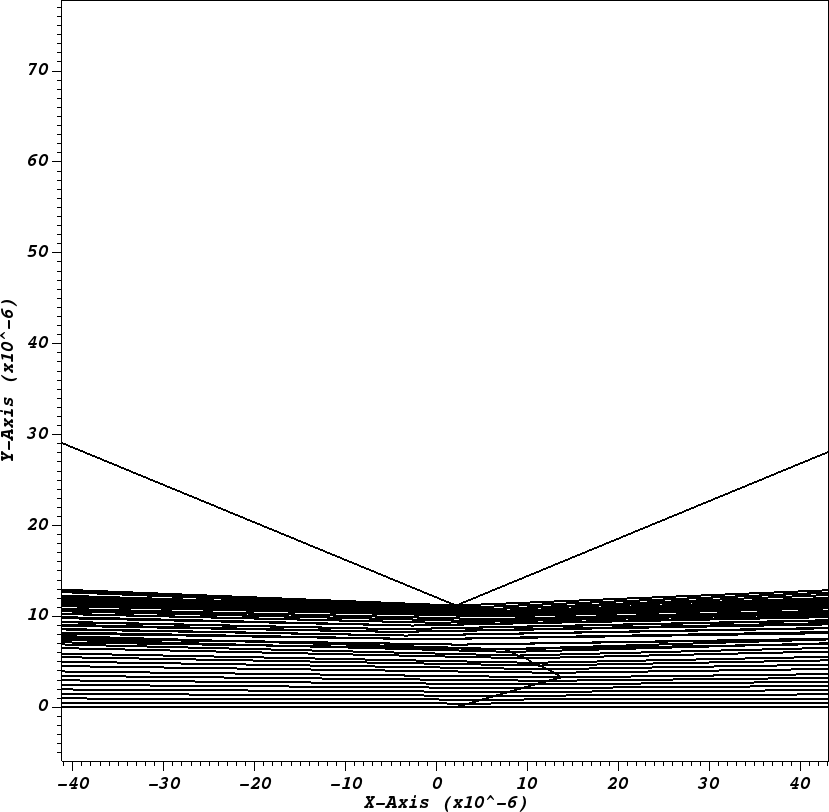}
    \caption{{[W] Left: Global mesh at the near-contact event deformed with the ALE mesh motion. Right: Zoom-in at the centre bottom.}}
    \label{fig:placeholder}
\end{figure}

\subsubsection*{[FST] Fara, Schwarzacher, T\r{u}ma}

Fara, Schwarzacher and T\r{u}ma use a monolithic fitted Updated ALE
method. The fluid and solid meshes conform at the interface, and velocity,
displacement and pressure are solved simultaneously in the ALE weak
formulation~\eqref{varALE1}. No contact model is used; the near-contact regime
is resolved by refining the thin fluid layer between the disk and the lower
wall.

\paragraph{Updated ALE and remeshing.}
The implementation follows~\cite{FarSchTum24}. Unlike Wick's full ALE
formulation~\eqref{ALE-FSI-prob}, which retains the initial reference
configuration throughout the simulation, our method replaces the 
{reference configuration}
by {the current} deformed configuration at remeshing times. In the reported
computations, remeshing is performed every 20 time steps and whenever the
time-step size changes. It is therefore a full ALE method between two
remeshing events, without remeshing at every time step. 
At a remeshing time $t_r$, the deformation is split as
$$
\mathbb{F}=\mathbb{F}_c\mathbb{F}_r,
$$
see~\eqref{eq:updated_ale_split}, where
$\mathbb{F}_c=\mathbb{I}+\nabla\mathbf{u}_c$ is the deformation gradient
relative to the current computational domain and $\mathbb{F}_r$ stores the
previous deformation history.

The fluid mesh displacement is obtained from a pseudo-elastic extension of the
solid displacement~\cite{HronTurek2006benchmark,FarSchTum24}:
$$
0=\mathcal{A}(\mathbf{u}^{\rm f},\boldsymbol{\phi})
\eqd{}\int_{\Omega^{\rm f}}
\left[\lambda_a\,\div\mathbf{u}^{\rm f}\,\div\boldsymbol{\phi} +
\mu_a\left(\nabla\mathbf{u}^{\rm f}+
(\nabla\mathbf{u}^{\rm f})^{\rm T}\right):
\nabla\boldsymbol{\phi}\right]\,\mathrm{d}x ,
$$
with interface values inherited from the solid displacement. In the present
computations, $\lambda_a=0$ and $\mu_a=2.91$, so that only the artificial
shear term is active. The artificial mesh equation is restricted to the fluid
domain and does not affect the physical solid response.

\begin{algorithm}[h!]
\caption{Remeshing in the Updated ALE method.}
\label{alg:FST-remeshing}
\begin{algorithmic}[1]
\Require $\mathbf{v}^n$, $\mathbf{u}_c^n$, $\mathcal{T}$,
$\mathbb{F}_r^n$
\If{$n\bmod 20=0$ \textbf{or} $\delta t^n\neq\delta t^{n-1}$}
    \State $\mathcal{T}_m \gets
    \texttt{move\_mesh}(\mathcal{T},\mathbf{u}_c^n)$
    \State $\mathbb{F}_r^n \gets
    (\mathbb{I}+\nabla\mathbf{u}_c^n)\mathbb{F}_r^n$
    \State $\mathcal{T}_{\mathrm{new}}\gets
    \texttt{adapt\_mesh}(\mathcal{T}_m)$
    \State $(\mathbf{v}^n,\mathbb{F}_r^n)\gets
    \texttt{transfer}
    ((\mathbf{v}^n,\mathbb{F}_r^n),
    \mathcal{T}_m,\mathcal{T}_{\mathrm{new}})$
    \State $\mathcal{T}\gets\mathcal{T}_{\mathrm{new}}$
    \State $\mathbf{u}_c^n\gets\mathbf{0}$
\EndIf
\State \Return $\mathbf{v}^n$, $\mathbf{u}_c^n$, $\mathcal{T}$,
$\mathbb{F}_r^n$
\end{algorithmic}
\end{algorithm}

Algorithm~\ref{alg:FST-remeshing} summarizes the remeshing procedure. The mesh
is first moved by the current ALE displacement into the deformed configuration.
This changes only the vertex coordinates, while the finite-element
coefficients remain associated with the moved mesh. The deformation history is
then updated according to
$
    \mathbb{F}_r^n \gets
    (\mathbb{I}+\nabla\mathbf{u}_c^n)\mathbb{F}_r^n.
$

The triangular mesh is remeshed using ADmesh~\cite{ADmesh}, developed for the
Updated ALE strategy~\cite{FarSchTum24}. ADmesh preserves the fluid--solid
interface and subdomain structure and improves the mesh by local edge
flipping, reduction, splitting, and vertex movement. In the narrow gap, the
requested local element size is guided by an Eikonal-type distance computation
so that several element layers are retained between the disk and the lower
wall.

After remeshing, $\mathbf{v}^n$ and $\mathbb{F}_r^n$ are transferred from the
moved mesh to the new mesh. The velocity is transferred either by standard
interpolation or by a divergence-free projection, which reduces the discrete
divergence error introduced by the transfer. The deformation-history tensor is
transferred by standard interpolation. Finally, the new mesh becomes the
computational reference and $\mathbf{u}_c^n$ is reset to zero.

\begin{figure}
    \centering
    \includegraphics[width=0.49\textwidth]{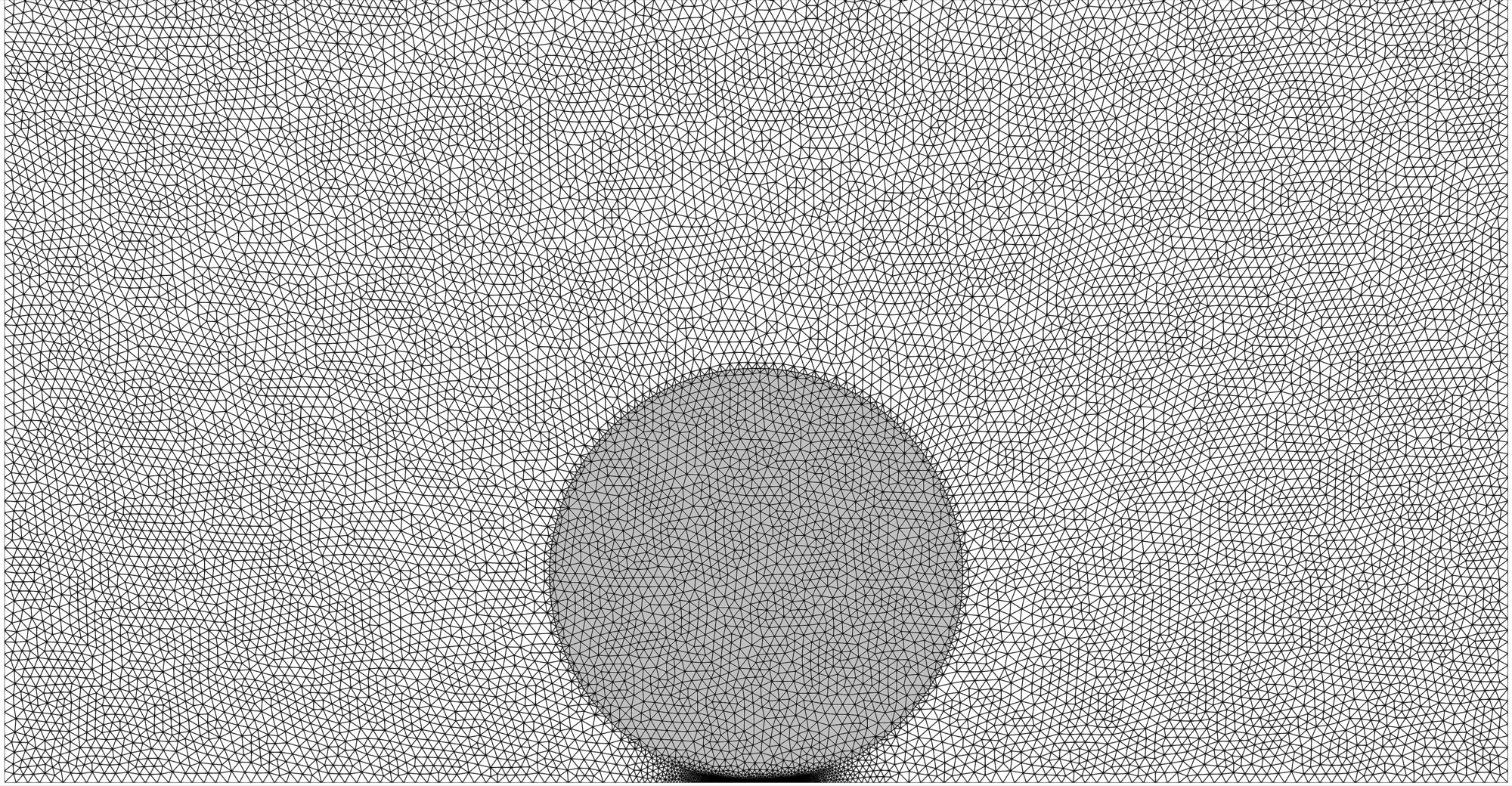}\hfill
    \includegraphics[width=0.49\textwidth]{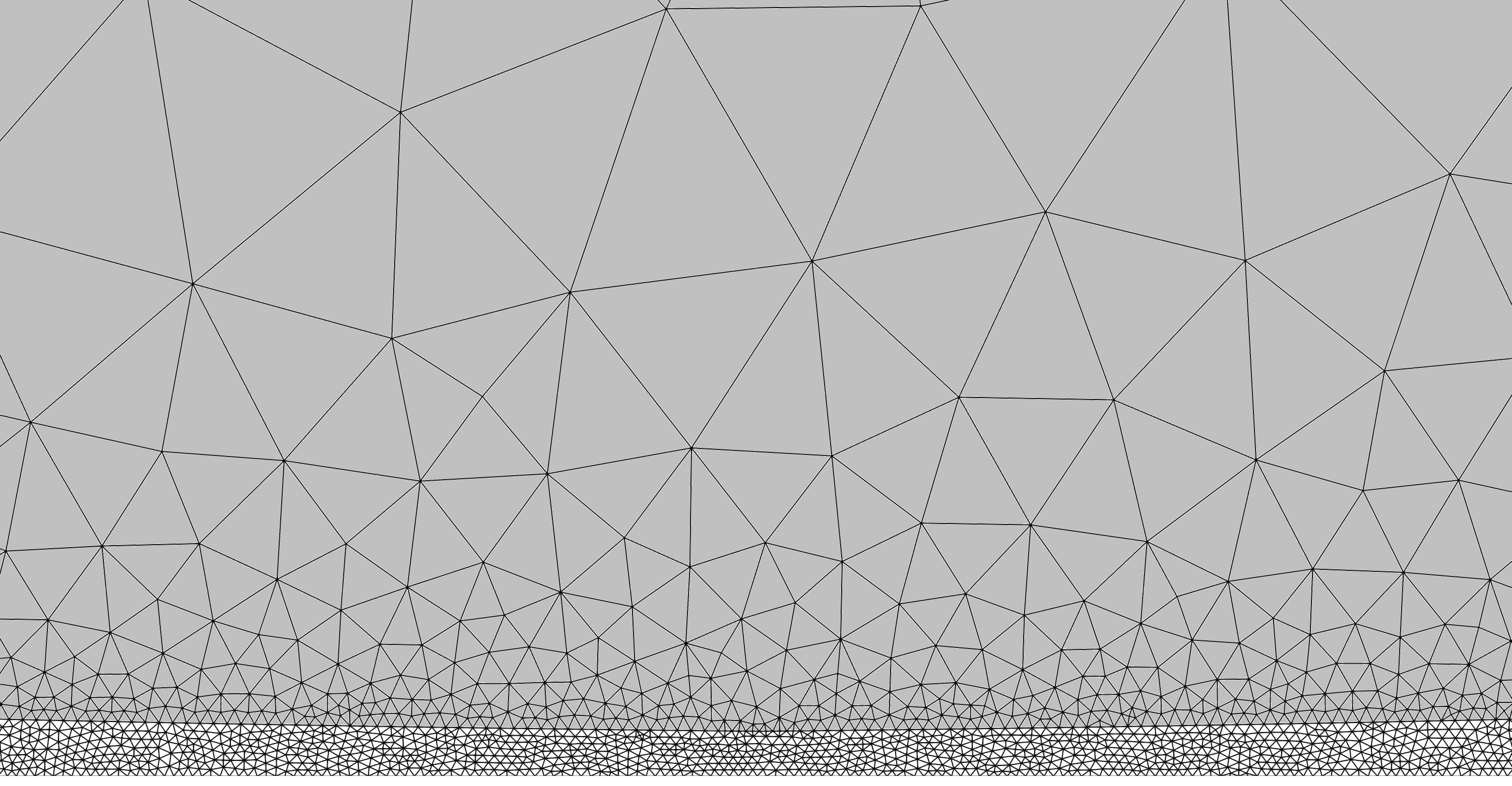}
    \caption{Left: global mesh in the current ALE configuration at
    $t=0.35736\,\mathrm{s}$, when the disk is closest to the lower boundary.
    Right: detail of a $0.2\,\mathrm{mm}\times0.1\,\mathrm{mm}$ region. The
    disk is shown in gray and the remaining fluid layer in white. The minimum
    gap is $5.92\,\mu\mathrm{m}$ and the fluid layer is resolved by about eight
    element layers.}
    \label{fig:FST_mesh_ball_colored}
\end{figure}

\paragraph{Spatial discretization.}
The method uses fitted triangular meshes and
$
    (\mathbf{v},\mathbf{u},p)\in P_2^2/P_2^2/P_1.
$
The velocity and displacement spaces are continuous across the fitted
fluid--solid interface. The pressure is also defined on the complete mesh, with
a Laplacian regularization in the solid subdomain. The circular interface is
represented by 300 boundary points. Three spatial resolutions are used. Due to
remeshing and local refinement, the numbers of degrees of freedom range
approximately from $2.0\times10^5$ to $4.4\times10^5$, from
$3.1\times10^5$ to $5.8\times10^5$, and from $5.5\times10^5$ to
$7.2\times10^5$, respectively. {An illustration of the mesh at the near-contact event is given in Figure~\ref{fig:FST_mesh_ball_colored}}. The implementation is based on
FEniCS~\cite{FEniCS}.

\paragraph{Temporal discretization.}
Time integration uses implicit backward Euler with
$\delta t_0=10^{-3}\,\mathrm{s}$. The time step is halved whenever the minimum
gap between the disk and the lower wall crosses one of the thresholds
\[
    \{10^{-2},\,4\cdot 10^{-3}, 10^{-3},\,4\times10^{-4},\,10^{-4},\,
      4\times10^{-5}\}\,\mathrm{m}.
\]
Thus,
\[
    \delta t_{\min}
    =
    \frac{\delta t_0}{64}
    =
    1.5625\times10^{-5}\,\mathrm{s}.
\]
Each change of the time-step size also triggers a reference update and
remeshing, independently of the regular 20-step update interval.

\paragraph{Numerical solution and force evaluation.}
The nonlinear system is solved monolithically by Newton's method, using a
relative and absolute tolerance of $10^{-8}$ and at most 20 Newton iterations
per time step. The linearized systems are solved by a direct LU factorization
with MUMPS. The drag force is evaluated from the fluid weak form using the
Babu\v{s}ka--Miller trick~\cite{babuska-miller-1984a,
babuska-miller-1984b,giles-larson-levenstam-suli-1997}, rather than by direct
integration of the discrete boundary traction.

\subsubsection*{[CF] Corti, Fern\'andez 
}

\begin{figure}[h!]
\centering 
\includegraphics[width=0.45\textwidth]{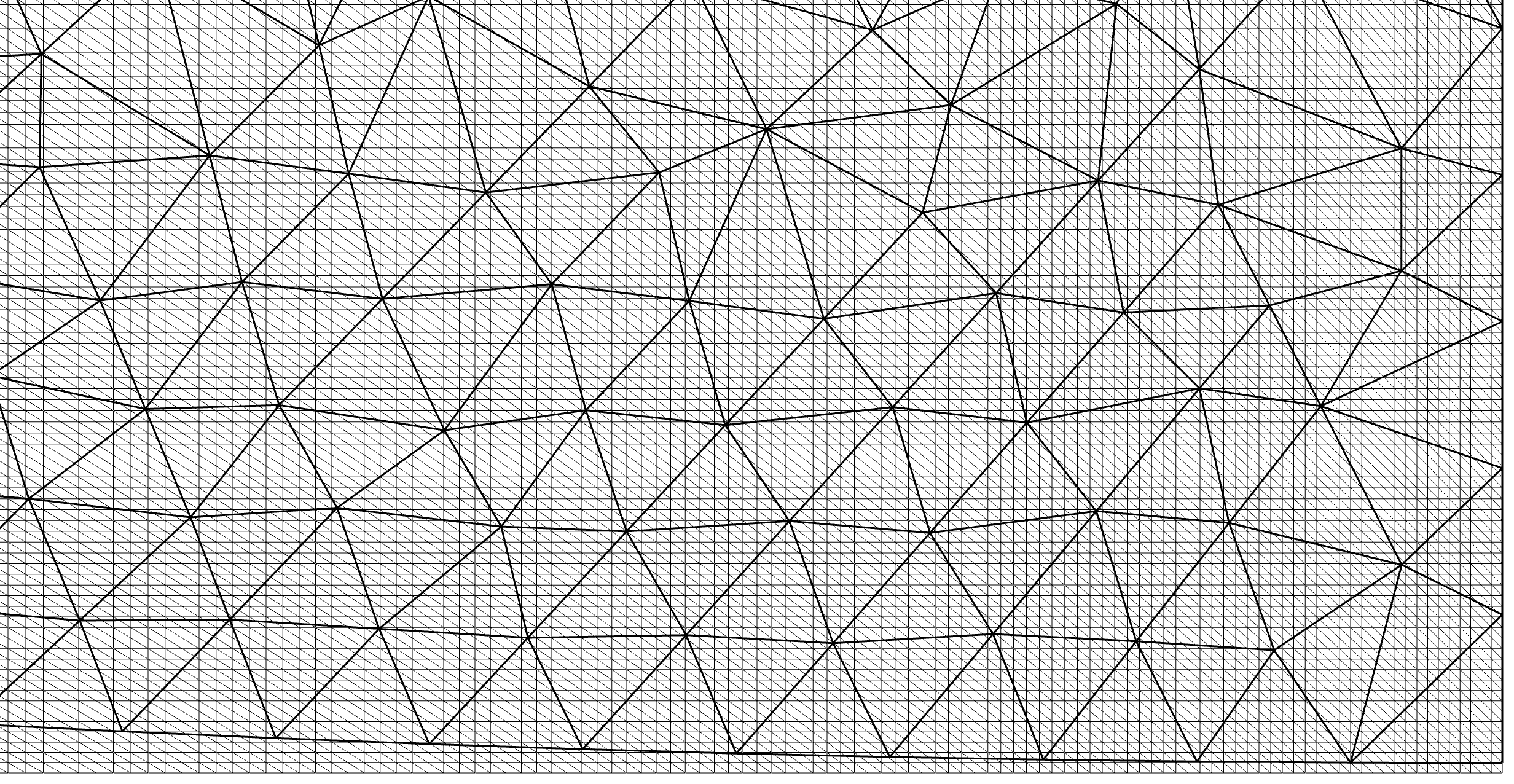}
\qquad 
\includegraphics[width=0.48\textwidth]{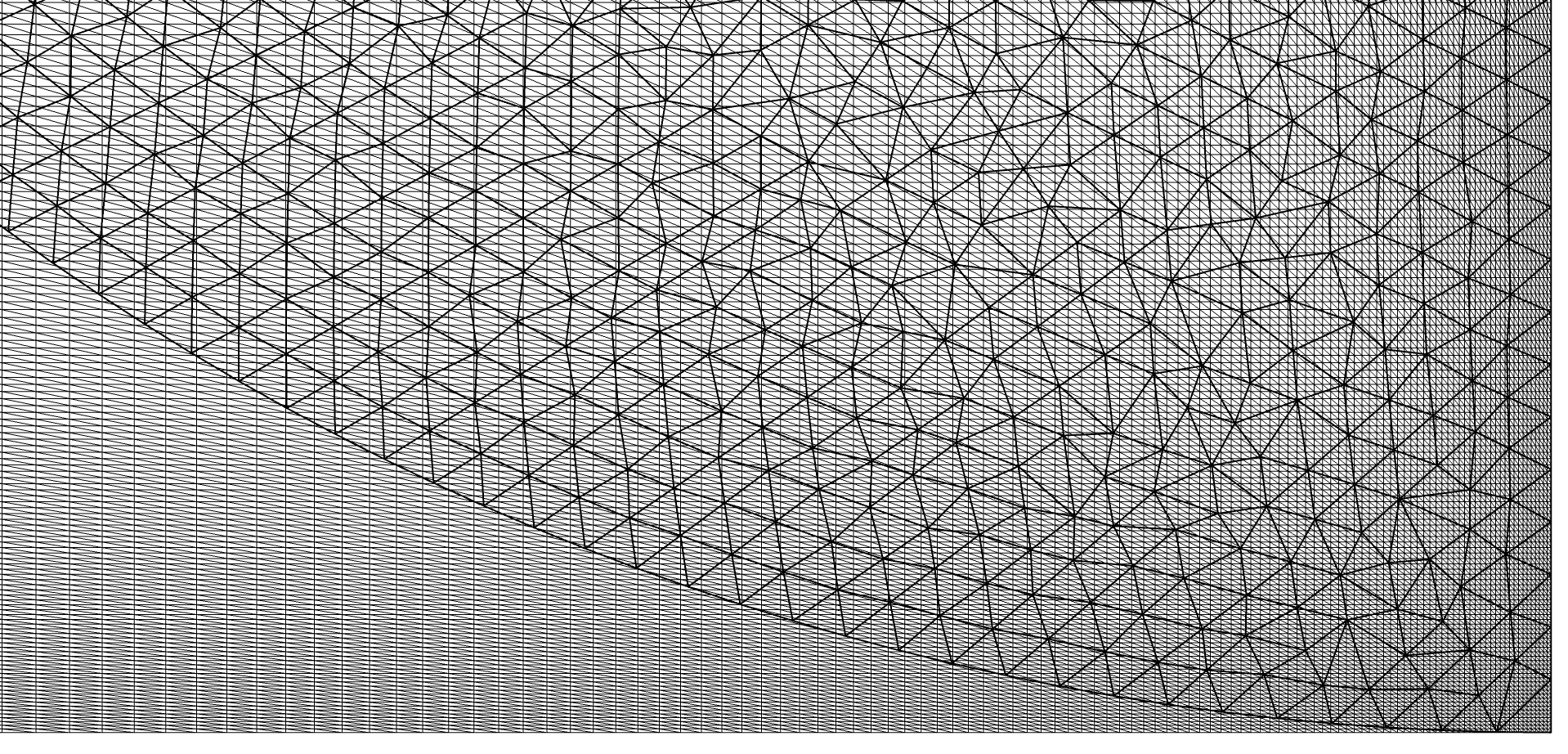}
\caption{Left: 2x1 mm close-up view of the finest meshes. Right: 22x11 mm close-up view of the coarsest meshes. {Note that the fluid mesh is barely visible, as it is much finer than the solid mesh.}}
\end{figure}

 We consider the  mixed-coordinate Lagrangian-Eulerian modeling 
 approach of Section~\ref{sec.lageul}. A fitted finite element  method on the discrete reference  domain $\Omega_{0,H}^{\rm s}$ is used for the solid, whereas, for the fluid, an unfitted  finite element approach in the whole fluid computational domain $\Omega$ is employed.  {Piece-wise  affine finite elements are used for all the discrete variables $\vec{v}^{\mathrm f,n}_h, p^{\mathrm f,n}_h, \vec{u}^{\mathrm s,n}_H$, $\vec{v}^{\mathrm s,n}_H$ and $\lambda_H^n$.}

Let 
${\vec u}^{\mathrm s,n}_H \in {\vec V}_H^{\rm s}$, for $n \geq 0$, denote a finite element
approximation of the solid displacement at time $t_n \eqd n\delta t$, 
where $\delta t >0$ stands for the time-step length. The
approximated solid motion map is given as $\vec X_H^n \eqd {\vec I} + {\vec u}^{\mathrm s,n}_H$.
We can then define the discrete domains
$
\Omega^{\rm s}[\vec X_H^n] \eqd \vec X_H^n\big(\Omega^{\rm s}_{0,H}\big)$, $
\Omega^{\rm f}[\vec X_H^n] \eqd \Omega_h \setminus \Omega^{\rm s}[\vec X_H^n]$, 
 $\Gamma^{\rm i}[\vec X_H^n] \eqd \vec X_H^n\big(\Gamma_{0,H}^{\rm i}\big)
= \partial \Omega^{\rm s}[\vec X_H^n].
$
For the numerical approximation of Model~I, we consider a CutFEM methodology for the fluid in combination with a  Nitsche treatment of the interface kinematic and
dynamic couplings (see, e.g., \cite{burman-fernandez-14,MassingLarsonLoggRognes2015,alauzet:hal-01149225,ZoncaVergaraFormaggia2018,SchottAgerWall2019}). The discrete fluid velocity and pressure are 
defined on the fixed whole computational domain $\Omega$.
As regards the time discretization,
a semi-implicit backward-Euler scheme is used in the fluid, whereas for the
solid we consider the mid-point scheme. To this purpose, we shall 
make use of the standard notation
$
\partial_{\delta t} x^n \eqd \frac{1}{\delta t}\bigl(x^n - x^{n-1}\bigr)$ and $
x^{n-\frac12} \eqd \frac{1}{2}\bigl(x^n + x^{n-1}\bigr).
$
The resulting numerical method is presented in 
Algorithm~\ref{alg:fsi-unfitted-mf}, corresponding to the [CF1] results of Section~\ref{results}. 
\begin{algorithm}[h!]
{
For $n\geq 1$:
\begin{enumerate}
\item Find $\bigl({\vec v}^{\mathrm f,n}_h, p^{\mathrm f,n}_h, {\vec u}^{\mathrm s,n}_H, {\vec v}^{\mathrm s,n}_H,\lambda_H^n\bigr)
  \in {\vec V}_h^{\rm f} \times Q_h^{\rm f} \times {\vec V}_H^{\rm s} \times {\vec V}_H^{\rm s} \times \Lambda_H$
  with ${\vec v}^{\mathrm s,n-\frac12}_H = \partial_{\delta t} {\vec u}^{\mathrm s,n}_H$ and such that
\begin{equation}\label{eq:discrete-semi-fsi}
\begin{aligned}
& \int_{\Omega^{\rm f}[\vec X_H^{n-1}]} \varrho^{\rm f}\, \partial_{\delta t} {\vec v}_h^{\mathrm f,n} \cdot {\boldsymbol \phi}_h^{\rm f}
  + \int_{\Omega^{\rm f}[\vec X_H^{n-1}]} \varrho^{\rm f}\, \bigl(\vec v_h^{f,n-1}\cdot\nabla\bigr) {\vec v}^{\mathrm f,n}_h \cdot {\boldsymbol \phi}_h^{\rm f}+ \int_{\Omega^{\rm f}[\vec X_H^{n-1}]} \mathbb{T}^{\rm f}\bigl({\vec v}^{\mathrm f,n}_h, p^{\mathrm f,n}_h\bigr) : \mathbb{D}({\boldsymbol \phi}_h^{\rm f})\\
& \quad 
  + \int_{\Omega^{\rm f}[\vec X_H^{n-1}]} q_h\,\div {\vec v}^{\mathrm f,n}_h  + \int_{\Omega^{\rm s}_{0,H}} \varrho^{\rm s}\, \partial_{\delta t}{\vec v}^{\mathrm s,n}_H \cdot {\boldsymbol \phi}_H^{\rm s}
  + \int_{\Omega^{\rm s}_{0,H}} \mathbb{S}\bigl({\vec u}^{s,n-\frac12}_H\bigr) : D_{\vec u^{\rm s}}\mathbb{E}\bigl({\vec u}^{s,n-\frac12}_H\bigr){\boldsymbol \phi}_H^{\rm s} \\
& \quad - \int_{\Gamma^{\rm i}[\vec X_H^{n-1}]} \mathbb{T}^{\rm f}\bigl({\vec v}^{\mathrm f,n}_h, p^{\mathrm f,n}_h\bigr)\vec{n} \cdot
  \bigl({\boldsymbol \phi}_h^{\rm f}- {\boldsymbol \phi}_H^{\rm s} \circ (\vec X^{n-1}_H)^{-1}\bigr) \\
& \quad - \int_{\Gamma^{\rm i}[\vec X_H^{n-1}]} \mathbb{T}^{\rm f}\bigl({\boldsymbol \phi}_h^{\rm f}, -q_h\bigr)\vec{n} \cdot
  \bigl({\vec v}^{\mathrm f,n}_h - {\vec v}^{\mathrm s,n-\frac12}_H \circ (\vec X^{n-1}_H)^{-1}\bigr) \\
& \quad + \int_{\Gamma^{\rm i}[\vec X_H^{n-1}]} \frac{\gamma\mu^{\rm f}}{h}
  \bigl({\vec v}^{\mathrm f,n}_h - {\vec v}^{\mathrm s,n-\frac12}_H \circ (\vec X^{n-1}_H)^{-1}\bigr) \cdot
  \bigl({\boldsymbol \phi}_h^{\rm f}- {\boldsymbol \phi}_H^{\rm s} \circ (\vec X^{n-1}_H)^{-1}\bigr) \\
& \quad + s_{\mathrm v,h}^{n-1}\bigl({\vec v}_h^{\mathrm f,n-1}; {\vec v}_h^{\mathrm f,n}, {\vec v}_h\bigr)
  + s_{\mathrm p,h}^{n-1}\bigl({\vec v}_h^{\mathrm f,n-1}; p_h^{\mathrm f,n}, q_h\bigr)
  + g_h^{\rm u}\bigl({\vec v}^{\mathrm f,n}_h, {\vec v}_h\bigr)
  + g_h^{\rm p}\bigl(p_h^{\mathrm f,n}, q_h\bigr) \\
 &\quad  - \int_{\Gamma^{\rm i}_{0,H}} \lambda_H^n {\boldsymbol \phi}_H^{\rm s} \cdot \vec{n}_{\rm b} 
 +  \int_{\Gamma_{0,H}^{\rm i}} \left(\lambda_H^n  + \left| \frac{\gamma_{\rm C} E }{H} \big(\vec{u}^{\mathrm s,n}_H\cdot \bs n_{\rm b}  - {g}_\epsilon \big)- \lambda_H^n  \right|_+  \right) \xi_H   = 0
\end{aligned}
\end{equation}
for all $({\boldsymbol \phi}_h^{\rm f}, q_h, {\boldsymbol \phi}_H^{\rm s},\xi_H) \in {\vec V}_h^{\rm f} \times Q_h^{\rm f} \times {\vec V}_H^{\rm s}\times \Lambda_H$.
\item Update interface and fluid domain:
$
\vec X_H^n = {\vec I} + {\vec u}^{\mathrm s,n}_H$, 
  $\Gamma^{\rm i}[\vec X_H^n] \eqd \vec X_H^n\bigl(\Gamma_{0,H}^{\rm i}\bigr), $
$\Omega^{\rm s}[\vec X_H^n] = \vec X_H^n\bigl(\Omega^{\rm s}_{0,H}\bigr)$,  $\Omega^{\rm f}[\vec X_H^n] = \Omega_h \setminus \Omega^{\rm s}[\vec X_H^n].$
\end{enumerate}}
\caption{[CF1] Mixed Lagrangian-Eulerian CutFEM with Nitsche coupling and relaxed contact.}
\label{alg:fsi-unfitted-mf}
\end{algorithm}

Note that an augmented Lagrangian 
formulation is used for contact, where 
 $E$ denotes the Young modulus of the solid, $\gamma_{\rm C} = 10$ and the contact relaxation parameter is taken as 
 $\epsilon = 10^{-6} / 2^i  $ for each  level  of refinement $i=0,1,2$, respectively. The bi-linear forms $g_h^{\rm u}$ and $g_h^{\rm p}$ are numerical  ghost-penalty operators   (see, e.g., \cite{Bu10,BH14,OLSHANSKII2025113983}):
\begin{equation}\label{eqn:GPForms}
g_h^{\rm u}({\vec v}^{\rm f}_h, {\boldsymbol \phi}_h^{\rm f})
\eqd \sum_{e\in {\mathcal E}_h} h\int_e \jump{\partial_{\vec n_e}{\vec v}^{\rm f}_h} \cdot \jump{\partial_{\vec n_e}{\boldsymbol \phi}_h^{\rm f}},\qquad
g_h^{\rm p}(p_h^{\rm f}, q_h)
\eqd \sum_{e\in {\mathcal E}_h} h^3 \int_e \jump{\partial_{\vec n_e} p_h^{\rm f}}\,\jump{\partial_{\vec n_e} q_h},
\end{equation}
with $\partial_{\vec n_e}$ denoting the normal derivative on the edge $e$, in the
set ${\mathcal E}_h$ of all internal edges of ${\mathcal T}_h$. Note that, in the
spirit of \cite{OLSHANSKII2025113983,fernandez-li-olshanskii-26}, we use ghost-penalty terms acting on the
whole computational domain $\Omega$. 
 The purpose of the ghost-penalty 
terms \eqref{eqn:GPForms} is precisely to guarantee the well-posedness of the resulting fluid sub-problem. 
{In order to cope with the
lack of inf-sup stability of the discrete spaces
and with
the numerical instabilities arising for large local Reynolds numbers, the
continuous interior penalty (CIP) stabilization method is used (see, e.g., \cite{burman-fernandez-hansbo-06,burman-fernandez-07}). The
associated symmetric velocity and pressure stabilization operators are given by:
$$
\begin{aligned}
s_{\mathrm v,h}^{n}({\vec z}_h; {\vec v}_h^{\rm f}, {\boldsymbol \phi}_h^{\rm f})
&\eqd \gamma_{\rm v}\, h^2 \sum_{e \in \mathcal{E}_{h}^{i,n}}
  \xi\bigl({\rm Re}_e(\vec z_h)\bigr)\, \|\vec z_h \cdot {\vec n}_e\|_{L^\infty(e)}
  \int_e \jump{\partial_{\vec n_e}{\vec v}^{\rm f}_h} \cdot \jump{\partial_{\vec n_e} {\boldsymbol \phi}_h^{\rm f}}, \\
s_{\mathrm p,h}^{n}({\vec z}_h; p_h^{\rm f}, q_h)
&\eqd \gamma_{\rm p}\, h^2 \sum_{e \in \mathcal{E}_{h}^{i,n}}
  \frac{\xi\bigl({\rm Re}_e({\vec z}_h)\bigr)}{\|{\vec z}_h\|_{L^\infty(e)}}
  \int_e \jump{\partial_{\vec n_e} p_h^{\rm f}}\,\jump{\partial_{\vec n_e} q_h},
\end{aligned}
$$
where $\mathcal{E}_{h}^{i,n}$ denotes the set of interior edges of
$\mathcal{T}_{h}$ which either lie on the physical or interfacial regions,
${\rm Re}_e({\vec z}_h) \eqd \varrho^{\rm f}\, \|{\vec z}_h\|_{L^\infty(e)}\, h / \mu^{\rm f}$
denotes the local Reynolds number, $\xi(x) \eqd \min\{1, x\}$ is a cut-off
function, and $\gamma = 10^2$ and $\gamma_{\rm p}= \gamma_{\rm v} = 10^{-2}$ are user-defined dimensionless
stabilization parameters.  
In the numerical results reported in Section~\ref{results}, the 
 coupled problem \eqref{eq:discrete-semi-fsi} is solved using a Dirichlet-Neumann partitioned solution 
procedure, based on Newton-GMRES interfacial iterations (see, e.g., 
\cite{fernandez:hal-00705110}).  
All the computations have been performed with the  \texttt{FELiScE} finite element library\footnote{\url{https://gitlab.inria.fr/felisce/felisce}},  using \texttt{PETSc/MUMPS} as linear solver.
 }

In Models II and III a porous layer is added on $\Gamma_{\text{bo{t}}}$ and $\Gamma_t^{\rm i}$, respectively, as shown in \eqref{eq:darcy}-\eqref{eq:darcy-coupling} and \eqref{eq:darcy-moving}-\eqref{eq:por:int}.  
Specifically for [CF2] (based on Model II), the additional discrete 
trial and test functions   $P_{\mathrm l,h}^n,Q_{\mathrm l,h}:\Gamma_{\text{bo{t}}}\to\mathbb{R}$  and the following terms are added to Algorithm~\ref{alg:fsi-unfitted-mf} (see  \cite{BurmanFernandezFreiGerosa2022}):
\begin{multline}\label{varFor_porous}
    \int_{\Gamma_{\text{bo{t}}} } \left( P_{\mathrm l,h}^n+ \frac{\epsilon_{\rm p} K_{\boldsymbol n}^{-1} }{4} \vec{v}_h^{\mathrm f,n}\cdot\vec{n}\right)  \vec\phi_h^{\rm f}\cdot \vec{n} 
		+ \int_{\Gamma_{\text{bo{t}}}} \epsilon_{\rm p} K_{\boldsymbol \tau} {\boldsymbol \nabla}_{\boldsymbol \tau} P_{\mathrm l,h}^n  \cdot {\boldsymbol \nabla}_{\boldsymbol \tau} Q_{\mathrm l,h}\\
		- \int_{\Gamma_{\text{bo{t}}}} \vec{v}_h^{\mathrm f,n}\cdot\vec{n} Q_{\mathrm l,h}
		+\frac{\alpha}{\sqrt{K_{\boldsymbol \tau}\epsilon_{\rm p}}} 
        \int_{\Gamma_{\text{bo{t}}}} \vec{v}_{h}^{\mathrm f,n}\cdot {\boldsymbol \tau}  {\boldsymbol \phi}_{h}^{\rm f} \cdot {\boldsymbol \tau} .
\end{multline}
Here, $P_{\mathrm l,h}^n,Q_{\mathrm l,h}$ are taken in the piecewise affine trace space $\{q_h\vert_{\Gamma_{\text{bo{t}}}}\,:\, q_h \in Q_h^{\rm f} \}$. Finally, 
 for [CF3] (based on Model III), the  interface Nitsche coupling terms in Algorithm~\ref{alg:fsi-unfitted-mf}  are replaced by the following terms (see, e.g., \cite{champion:hal-05626925})
\begin{multline*}
\int_{\Gamma^{\rm i}[\vec{X}_H^{n-1}]} P_{\mathrm l,H}^n 
\big( {\boldsymbol \phi}_h^{\rm f}  - {\boldsymbol \phi}_H^{\rm s} )
\cdot \vec{n}  - \int_{\Gamma^{\rm i}
[\vec{X}_H^{n-1}]} 
\big(  \vec{v}_h^{\mathrm f,n}    
- \vec{v}^{\mathrm s,n-\frac12}_H\big) \cdot \vec{n}  Q_{\mathrm l,H} \\ 
+ \epsilon_{\rm p} K_{\boldsymbol \tau} 
\int_{\Gamma^{\rm i}[\vec{X}_H^{n-1}]} 
{\boldsymbol \nabla}_{\boldsymbol \tau} P_{\mathrm l,H}^n \cdot 
{\boldsymbol \nabla}_{\boldsymbol \tau} Q_{\mathrm l,H }
+ \frac{\alpha}{\sqrt{K_{\boldsymbol \tau}\epsilon_{\rm p}}}  \int_{\Gamma^{\rm i}[\vec{X}_H^{n-1}]}  \big(  \vec{v}_h^{\mathrm f,n}   
- \vec{v}_H^{\mathrm s,n-\frac12}\big) \cdot \boldsymbol  \tau  \big({\boldsymbol \phi}_h^{\rm f} - {\boldsymbol \phi}_H^{\rm s})\cdot {\boldsymbol  \tau}, 
\end{multline*}
where the  continuous piece-wise affine 
trial and test functions   $P_{\mathrm l,H}^n,Q_{\mathrm l,H}$ are now defined on the moving  interface $\Gamma^{\rm i}[\vec{X}_H^{n-1}]$.  
Note that this essentially amounts to enforcing the fluid–solid interface coupling in the normal direction using $P_{\mathrm l,H}^n$ as  (stabilized) Lagrange multiplier, and in the tangential direction through the penalty friction term. 

\subsubsection*{[Fr] Frei}
We use a Fully Eulerian approach with a fitted finite element discretization. In particular, we use a fixed coarse background mesh consisting of \textit{patches} which is identical for all time steps and apply a refinement in each time step that resolves the interface $\Gamma_t^{\rm i}$ linearly in each time step (\textit{locally modified finite elements}~\cite{FreiRichter2014}). This induces a splitting into (possibly anisotropic) solid and fluid cells on the finest mesh. For the coupling of fluid and solid equations, Nitsche's method is used, see below. The interface is captured using the Initial Point Set Method~\cite{Dunne2006}.

\paragraph{Spatial discretization.} We exploit symmetry of the configuration and simulate on half of the domain only. $Q_1$ finite elements are used on a mesh consisting of quadrilaterals for the four variables $(v_h^{\rm f}, v_h^{\rm s}, u_h, p_h) \in {\cal X}_h^n\eqd{\cal V}_f^{h,n} \times {\cal V}_s^{h,n} \times {\cal U}_h^n \times {\cal Q}_h^n$, where ${\cal V}_f^{h,n}$ and ${\cal Q}_h^n$ are defined on the set of fluid cells at time $t_n$ and ${\cal V}_s^{h,n}, {\cal U}_h^n$ in solid cells only. Due to the lack of inf-sup stability, we use an anisotropic variant of the interior penalty pressure stabilization (see~\cite{FreiPressure})
\begin{align*}
S^n(p_h, q_h)\eqd\alpha_{p,1} \sum_{e \in {\cal E}_h^{\text{int},n}} 
 \int_e \{h_n^3 {\boldsymbol \nabla}p_h\cdot {\boldsymbol \nabla}q_h\}_e \, do
 + \alpha_{p,2}\sum_{e \in {\cal E}_h^{\mathrm f,n}} \int_e \{h_n^3\}_e \jump{{\boldsymbol \nabla}p_h}_e\cdot \jump{{\boldsymbol \nabla}q_h}_e \, do
\end{align*}
where ${\cal E}_h^{\text{int},n}$ is the set of (possibly anisotropic) fluid edges that belong to a patch cut by the interface at time $t_n$, while  
${\cal E}_h^{\mathrm f,n}$ is the set of the remaining (interior) fluid edges. Moreover, $h_{n|K}\eqd\frac{|K|}{|e|}$ is the cell size of a cell $K$ in the direction normal 
to $e$ and $\{v_h\}_e$
is the mean value of the two cells $K_1,K_2$ sharing the edge $e$.

The mesh is highly refined around the central part of $\Gamma_{\text{bot}}$, in order to resolve the contact dynamics accurately. Coarse mesh cells are used close to the lateral boundaries and the top. The number of mesh cells is the same over all time steps and is given by $24\,929$, $49\,665$ and $99\,009$ on the three meshes used.

\begin{figure}
\begin{minipage}{0.51\textwidth}
\includegraphics[width=\textwidth]{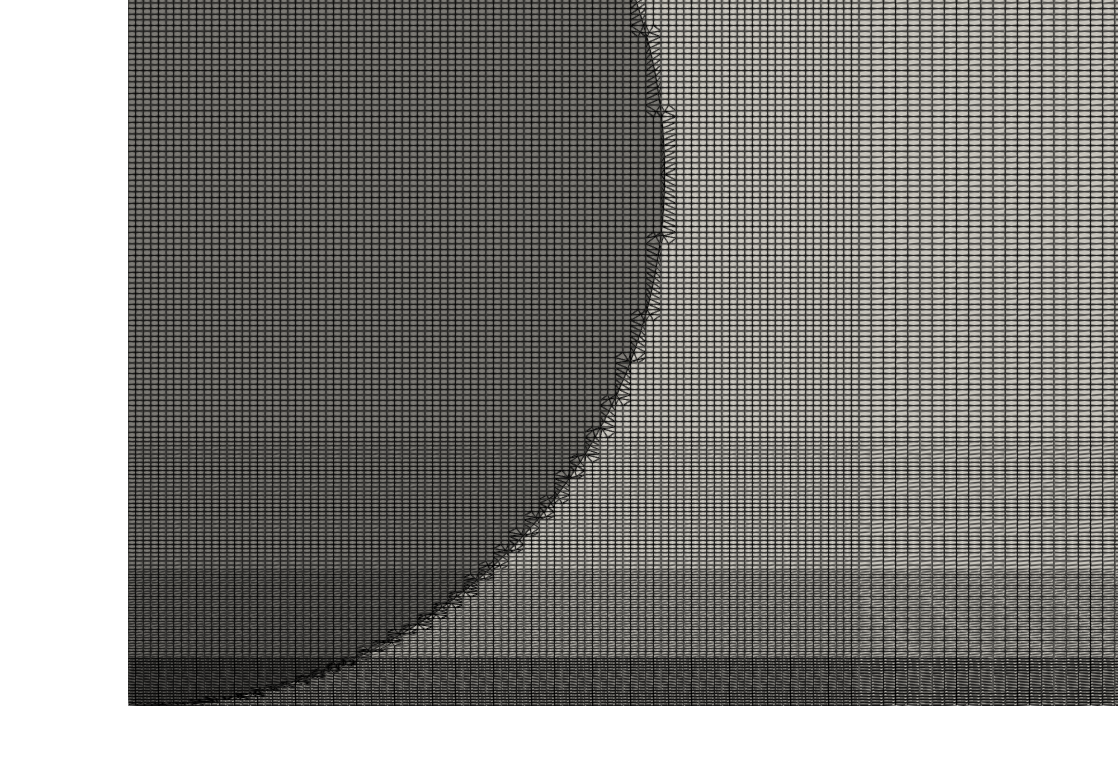}
\end{minipage}
\hfil
\begin{minipage}{0.49\textwidth}
\vspace*{-0.2cm}
\includegraphics[width=\textwidth]{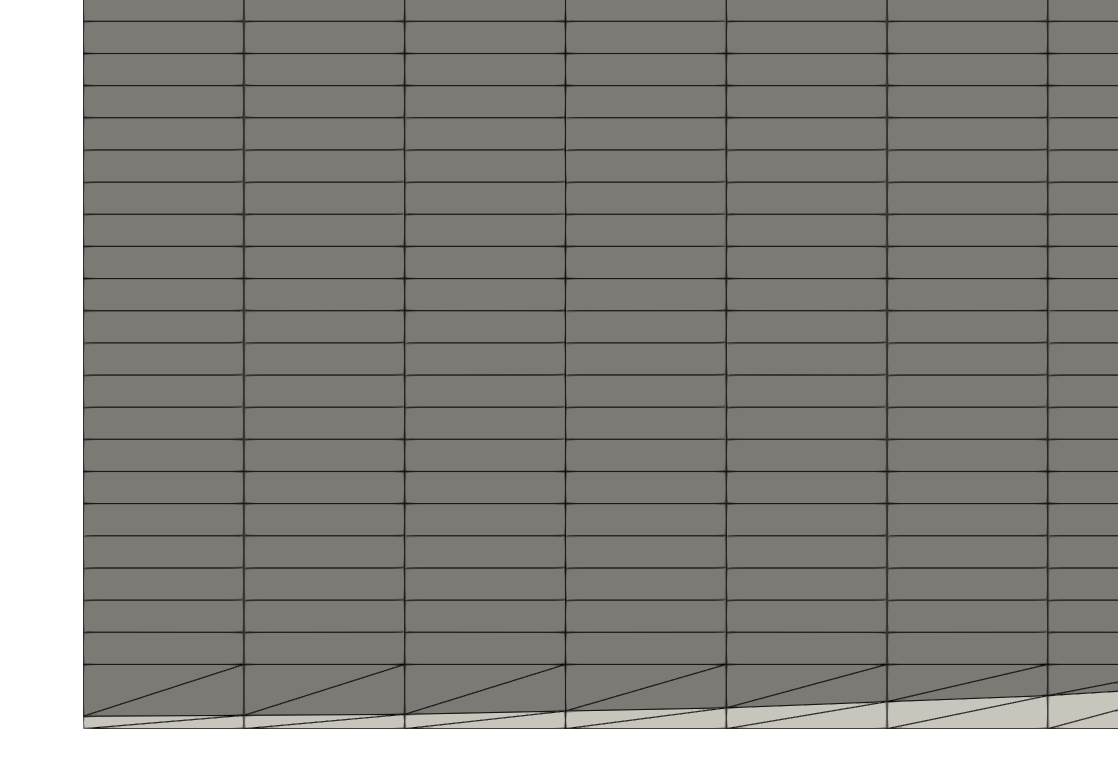}
\end{minipage}
    \caption{{[Fr] Illustration of the layered mesh at contact. Note that half the domain is used for simulation only. Left: Larger picture showing almost the entire half-disk: $2$cm $\times 1.5$cm. Right: zoom-in  around the contact region: $1$mm $\times 0.7$mm. Note that on the finest mesh level the cells are adjusted to resolve the interface. Only one layer of mesh cells is used in the contact region.}}
\end{figure}

\paragraph{Nitsche's method.}
The coupling in the Eulerian method is enforced weakly by using Nitsche's method, see e.g.~\cite{BurmanFernandezFrei2020, BurmanFernandezFreiGerosa2022, FreiKnokeSteinbachWenskeWick2026}. The semi-discrete variational formulation reads for $t\in [0,T]$: Find $(\vec v_h^{\rm f}, \vec v_h^{\rm s}, \vec u_h, p_h) \in {\cal X}_h$, such that
\begin{align}
  \sum_{i\in \{f,s\}} \left(\int_{\Omega^i_t} \tilde\varrho^i \pder{\vec{v}_h^i}{t}\cdot {\boldsymbol \phi}^i\right) +{\cal A}^t_{\text{Eu,N}}(\vec{v}_h^{\rm f}, \vec{v}_h^{\rm s}, \vec{u}_h, p; {\boldsymbol \phi}_h^{\rm f}, {\boldsymbol \phi}_h^{\rm s}, q, \vec{z}_h^{\rm s})  &=  \sum_{i\in \{f,s\}} \int_{\Omega^i_t} \tilde\varrho^i\vec{g}\cdot{\boldsymbol \phi}_h^i
  \end{align} 
  for all $({\boldsymbol \phi}_h^{\rm f}, {\boldsymbol \phi}_h^{\rm s}, q_h, \vec{z}_h^{\rm s}) \in  {\cal X}_h$, where
  \begin{align}\label{eq:Eu:Nitsche}
  \begin{split}
  {\cal A}^t_{\text{Eu,N}}(\vec{v}^{\rm f}, \vec{v}^{\rm s}, \vec{u}, p; {\boldsymbol \phi}^{\rm f}, {\boldsymbol \phi}^{\rm s}, q, \vec{z}^{\rm s}) \eqd  \sum_{i\in \{f,s\}} \left(
  \int_{\Omega^i_t} \tilde\varrho^i ( {\boldsymbol \nabla}\vec{v}^i) \vec{v}^i \cdot {\boldsymbol \phi}^i + \int_{\Omega^i_t} \mathbb{T}^i : {\boldsymbol \nabla}{\boldsymbol \phi}^i \right)  -\int_{\Omega^{\rm f}_t} q\,\div\,\vec{v}^{\rm f}  \\
 + \gamma_{\text{fsi}}\mu^{\rm f} h^{-1} \int_{\Gamma_t^{\rm i}} (\vec{v}^{\rm f}-\vec{v}^{\rm s})\cdot({\boldsymbol \phi}^{\rm f}-{\boldsymbol \phi}^{\rm s}) - \int_{\Gamma_t^{\rm i}}\mathbb{T}^{\rm f} \vec{n}\cdot ({\boldsymbol \phi}^{\rm f}-{\boldsymbol \phi}^{\rm s}) - \int_{\Gamma_t^{\rm i}}  (\vec{v}^{\rm f}-\vec{v}^{\rm s})\cdot\mathbb{T}^{\rm f}(\phi^{\rm f},-q) \vec{n} \\ 
 +  \int_{\Omega^{\rm s}_t} \left({ ({\boldsymbol \nabla}\vec{u}^{\rm s}) \vec{v}^{\rm s} - \vec{v}^{\rm s}} \right)\cdot\vec{z}^{\rm s}
  \end{split}
\end{align}
 for $\gamma_{\text{fsi}}>0$.
 
\paragraph{Time discretization.} For time discretization, we use a variant of the backward Euler method that
follows characteristics (see~\cite{FreiRichter2017}). To this purpose, we define the discrete time derivative as
\[
\widetilde{\partial}_{\delta t} \vec w_h^n
=
\frac{
\vec w_h^n -
\left(
\vec w_h^{n-1}\circ \Psi
\right)
}{\delta t}
-
\partial_t \Psi \cdot {\boldsymbol \nabla}\vec w_h^n ,
\]
where $\Psi$ is an (arbitrary) map defined in $\Omega$ that maps
$\Omega_f^n$ to $\Omega_f^{n-1}$ and $\Omega_s^n$ to
$\Omega_s^{n-1}$, respectively. In this way the variables $\vec w_h^{n-1}$ at the previous time step are only taken on the domain $\Omega_f^{n-1}$ resp. $\Omega_s^{n-1}$, where they are physically defined. In order to avoid unphysical pressure oscillations due to the locally modified meshes, we add the additional stabilization term $S_{\delta t}^n(p_h^n, p_h^{n-1}; q_h) \eqd \alpha_{\delta t} h(p_h^n-p_h^{n-1}, q_h)_{\Gamma_{\rm int}^n}$ on the interface ${\Gamma_{\rm int}^n}$ at time $t_n$.
The time step is chosen as
$\delta t= 2\cdot 10^{-3}\,\mathrm{s}$ in the beginning and decreased
adaptively down to $\delta t= 7.8125\cdot 10^{-6}\,\mathrm{s}$ depending on the
minimal distance between disk and $\Gamma_{{\rm bot}}$. As a consequence, a very small time step is used, shortly before, during and shortly after contact. 

\paragraph{[Fr1] Relaxed contact approach.} The first variant uses no-slip conditions on $\Gamma_{\text{bot}}$ in combination with a relaxed contact approach, see~\ref{subsec.contact}. This means that contact conditions get
active once the minimum distance between disk and $\Gamma_{\text{bot}}$
lies below $\varepsilon = c_0h$, where $c_0=\frac{1}{7}$. We use the Alart-Curnier trick to include the  corresponding contact conditions in combination with the FSI coupling conditions, see~\eqref{AlartCurnier}. The Lagrange multiplier $\lambda$ is eliminated and the condition is imposed  weakly by using Nitsche's method, where $\gamma_c=1$. For the details, we refer to~\cite{BurmanFernandezFrei2020}. The resulting algorithm is given in Algorithm~\ref{alg:Fr}. 
Note that the contact term is defined slightly different compared to ~\eqref{AlartCurnier}. This is because in Eulerian coordinates the no-penetration condition in \eqref{eq:contact_ineq2} reads 
\begin{align*}
(\vec{u}_h^n - \big(\vec{u}_h^{n-1}\circ \Psi\big))\cdot \vec{n} - \tilde{g}_\epsilon^n \leq 0
\end{align*} 
where $\tilde{g}_\epsilon^n$ denotes the current distance to the lower wall $\Gamma_{\rm bot}$ minus $\epsilon$. We refer to~\cite{BurmanFernandezFreiGerosa2022} for further details.

\begin{algorithm}[t]
	{For $n\geq 1$:
	\begin{enumerate}
		\item Update the domain affiliations  
		\begin{align*}
\Omega^{{\rm s},n} \eqd\left\{ \bx\in \Omega, \, \big| \, \bx-\vec u_h^{n-1} \in \Omega^{\rm s}\right\}, \quad \Gamma_{\rm int}^n \eqd\left\{ x\in \Omega, \, \big| \, \bx-\vec u_h^{n-1} \in \Gamma_{\rm int}^0\right\}, \quad \Omega^{\mathrm f,n} \eqd \Omega \setminus \left(\Omega^{{\rm s},n}\cup \Gamma_{\rm int}^n\right).  
\end{align*}
		\item Find $\vec y_h^n\eqd(\vec v_h^{\mathrm f,n}, \vec v_h^{\mathrm s,n}, \vec u_h^n, p_h^n) \in {\cal X}_h^n$ with
		\begin{multline}\label{varf:thick}
		\rho_{\rm f} \big( {\tilde{\partial}_{\delta t}} \vec v_h^{\mathrm f,n},  \vec \phi_h^{\rm f} \big)_{\Omega^{\mathrm f,n}} + \tilde{\rho}_{\rm s} \big( {\tilde{\partial}_{\delta t}} \vec v_h^{\mathrm s,n},  \vec \phi_h^{\rm s} \big)_{\Omega^{\mathrm s,n}} +  \big( {\tilde{\partial}_{\delta t}} \vec u_h^{n},  \vec z_h^{\rm s} \big)_{\Omega^{\mathrm s,n}}+ {\cal A}_{Eu,N}^{t_n}((\vec y_h^n; {\boldsymbol \phi}_h)  \\
		+\frac{\gamma_{{c}} \lambda^{\rm s}}{h} \big(\big[{P}_{\gamma_c}  \big]_+,\vec \phi_h^{\rm s}  \cdot \vec n \big)_{\Gamma_{\rm int}^n}
        +S^n(p_h^n, q_h) + S_{\delta t}^n(p_h^n, p_h^{n-1}; q_h)=  \sum_{i\in \{f,s\}} \int_{\Omega^i_t} \tilde\varrho^i\vec{g}\cdot{\boldsymbol \phi}_h^i	
		\end{multline}
		for all $\vec\Phi_h^n\eqd({\boldsymbol \phi}_h^{\rm f}, {\boldsymbol \phi}_h^{\rm s}, q_h, \vec{z}_h^{\rm s}) \in {\cal X}_h^n$ with the contact term ${P}_{\gamma_c}\eqd \big(\vec{u}_h^n - \big(\vec{u}_h^{n-1}\circ \Psi\big)\big)\cdot \vec{n} - \tilde{g}_\epsilon^n -\gamma_c^{-1}h  \vec n^\top (\mathbb{T}^{\mathrm f,n} -\mathbb{T}^{\mathrm s,n})\vec{n}$.
	\end{enumerate}}
    \caption{[Fr1]: Nitsche-based fully Eulerian finite element with relaxed contact condition\label{alg:Fr}}
    \end{algorithm}

    \paragraph{[Fr2] Porous layer contact model.}
In the second approach a porous layer is added on $\Gamma_{\text{bot}}$. We proceed exactly as [CF2] by adding the terms $a_p( \vec{v}_h^{\mathrm f,n},P_{\mathrm l,h}^n; {\boldsymbol \phi}_{h}^{\rm f}, q_{\mathrm l,h})$ defined in~\eqref{varFor_porous}. 
Both the FSI coupling and the contact conditions are incorporated in
a weak sense using Nitsche's method applied to~\eqref{AlartCurnier}. 


\paragraph{Numerical Solution.}
Both approaches use a semi-smooth Newton method to solve the non-linear system of equations, including the non-smooth term $([P_{\gamma_c}]_+, \phi_h^{\rm s}\cdot \vec n)_{\Gamma_{\rm int}^n}$, in combination with UMFPACK for the solution of the linearized equations.

\subsubsection*{[KWF] Knoke, Wick, Frei}

We use a fully Eulerian approach with the cut finite element method, where we allow the interface to cut freely through mesh cells resulting in overlapping computational subdomains $\Omega_{h,n}^{\rm f}$ and $\Omega_{h,n}^{\rm s}$. To extend coercivity of the bilinear form {onto} the interface region $\Omega_{h,n}^{\rm f} \cap \Omega_{h,n}^{\rm s}$ we add the ghost penalty terms 
{\begin{align} 
    \begin{split} \label{eq:def_gp}
    &G_h(\vec{v}_h^{\mathrm f,n}, p_h^n, \vec{v}_h^{s, n}, \vec{u}_h^{s, n};{\boldsymbol \phi}_h^{\rm f}, \vec{z}_h^{\rm s}, {\boldsymbol \phi}_h^{\rm s}, q_h) \eqd 
    2 \mu^{\rm f} \gamma_{\vec{v}^{\rm f}} \sum_{F \in \mathcal{F}_G^{\rm f}} \left(h (\llbracket \partial_\vec{n} \vec{v}_h^{\mathrm f,n} \rrbracket, \llbracket \partial_\vec{n} {\boldsymbol \phi}_h^{\rm f} \rrbracket)_F + \frac{h^3}{4} (\llbracket \partial_\vec{n}^2 \vec{v}_h^{\mathrm f,n} \rrbracket, \llbracket \partial_\vec{n}^2 {\boldsymbol \phi}_h^{\rm f} \rrbracket)_F\right) \\ 
    &+\gamma_{p} \sum_{F \in \mathcal{F}_G^{\rm f}} h^3 (\llbracket \partial_\vec{n} p_h^n \rrbracket, \llbracket \partial_\vec{n} q_h\rrbracket)_F, 
    + \tilde{\varrho^{\rm s}} \gamma_{\vec{v}^{\rm s}} \sum_{F \in \mathcal{F}_G^{\rm s}} h^3 (\llbracket \partial_\vec{n} \delta_n \vec{v}_h^{s}\rrbracket, \llbracket \partial_\vec{n} {\boldsymbol \phi}_h^{\rm s} \rrbracket)_F + 2 \mu^{\rm s} \gamma_{\vec{u}^{\rm s}} \sum_{F \in \mathcal{F}_G^{\rm s}} h (\llbracket \partial_\vec{n} \vec{u}_h^{s, n} \rrbracket, \llbracket \partial_\vec{n} {\boldsymbol \phi}_h^{\rm s} \rrbracket)_F,
    \end{split}
\end{align}}
{where} $\gamma_{\vec{v}^{\rm f}}=\gamma_p=\gamma_{\vec{u}^{\rm s}}=10^{-2}, \gamma_{\vec{v}^{\rm s}}=10^{-8}$, $\delta_n \vec{v}_h^{s} \eqd \delta t^{-1} (\vec{v}_h^{s, n} - \vec{v}_h^{s, n-1})$ and $\mathcal{F}_G^{\rm f}$ resp.\,$\mathcal{F}_G^{\rm s}$ are the sets of cell edges {of $\Omega_{h,n}^{\rm f}$ resp.\,$\Omega_{h,n}^{\rm s}$ that do not lie on the respective} boundary and where at least one of the adjacent cells is cut by the interface. The interface conditions are imposed using Nitsche’s method (see \eqref{eq:Eu:Nitsche}). A stability  and convergence analysis of this approach has been presented (for the simplified case of a fixed interface) in \cite{FreiKnokeSteinbachWenskeWick2026}.

\begin{figure}[h!]
    \centering
    \begin{minipage}[b]{.45\linewidth}
    \centering
    \includegraphics[scale=0.15]{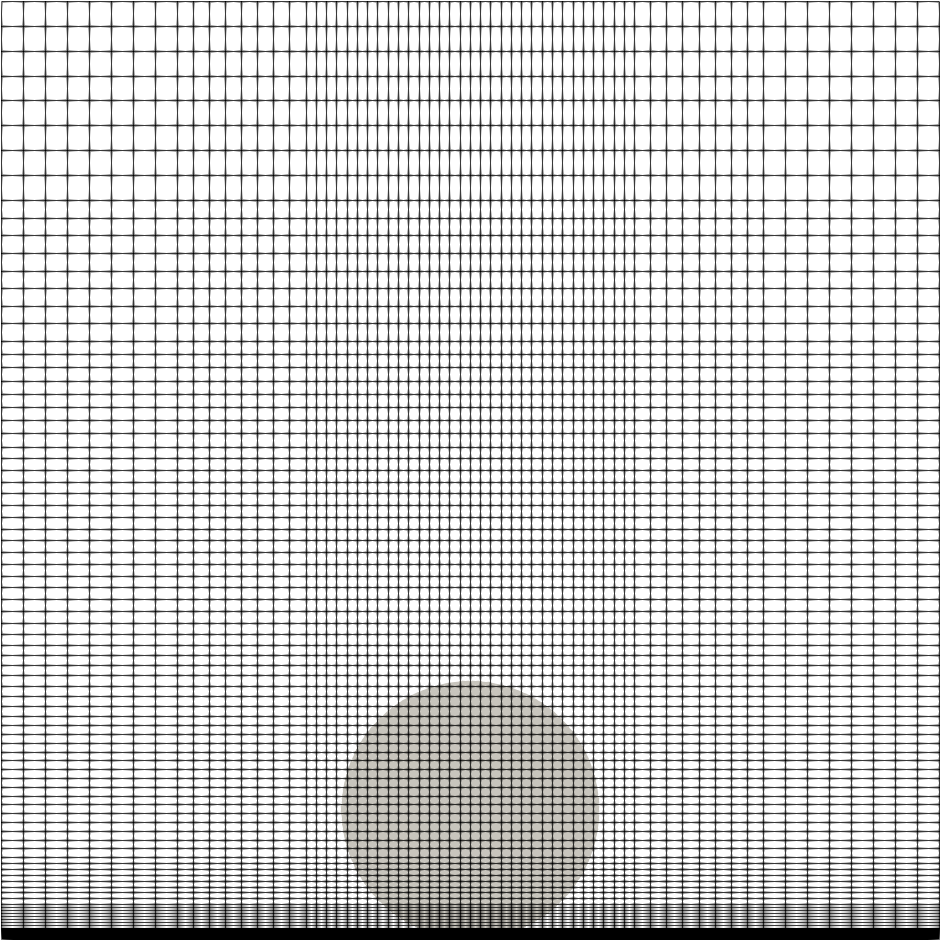}
    \end{minipage}
    \begin{minipage}[b]{.45\linewidth}
    \centering
    \includegraphics[scale=0.15]{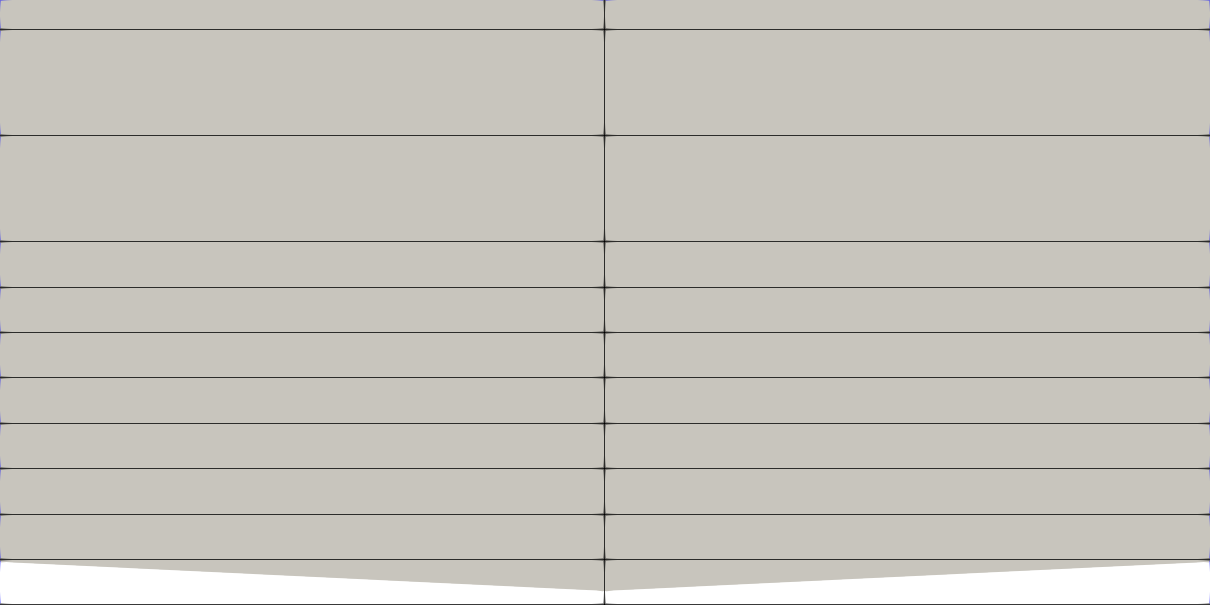}
    \end{minipage}
    \caption{Visualization of the global mesh (left) and detail view of a $1\,\mathrm{mm}\times 0.5\,\mathrm{mm}$ rectangle (right) of [KWF] at $t = 0.3580078\,\mathrm{s}$, when the disk is closest to the container boundary. The disk is shown in gray and the fluid region in white. 
    }
    \label{fig:KWF_mesh_ball_colored}
\end{figure}

\paragraph{Spatial discretization}
We use an anisotropic quadrilateral mesh where the cell size decreases vertically toward the bottom and increases horizontally toward the boundaries resulting in a maximum cell size of $1.875 \cdot 10^{-3} \times 2.125 \cdot 10^{-3}$ and a minimum cell size of $8.75 \cdot 10^{-4} \times 3.75 \cdot 10^{-5}$ (see Fig. \ref{fig:KWF_mesh_ball_colored}). For the discretization we use continuous {$Q_2$} elements for the
fluid velocity and continuous {$Q_1$} elements for the remaining solution components. Since the position of the interface and thus the computational subdomains depend on time the number of degrees of freedom varies in the range $[55\,144, 56\,071]$.

\paragraph{Temporal discretization.}
We use the A-stable first-order backward Euler method with initial time step sizes $\delta t\in \{10^{-3}, 5\cdot 10^{-4}, 2.5 \cdot 10^{-4}\}$ that are halved five times at the time thresholds $0.3, 0.32, 0.33, 0.34$ and $0.35$ {in order to resolve the contact dynamics accurately} and coarsened up again to the initial size when $t \ge 0.36$.

\paragraph{Implicit extensions.}
To deal with the moving interface we extend the computational subdomains by one cell layer and compute extensions of the solution by applying additional ghost penalty terms {on the newly included edges analogously to \eqref{eq:def_gp}. We denote the resulting stabilization by $G_h^\text{ext}$ with parameters $\gamma_p^\text{ext}=\gamma_{\vec{v}^{\rm s}}^\text{ext}=\gamma_{\vec{u}^{\rm s}}^\text{ext}=3 \cdot 10^{-3}, \gamma_{\vec{v}^{\rm f}}^\text{ext}=0.6$.} The position of the interface is represented by a level set function \[\Phi(x,t_n)\eqd\left\{\begin{array}{ll} \Phi_0(x-\vec{u}_{h,\text{ext}}^{\rm s}(x,t_n)), & x\in \Omega_{h,n}^{s,\text{ext}} \\ 1, & x\not\in \Omega_{h,n}^{s,\text{ext}} \end{array}\right.,\] where $\vec{u}_{h,\text{ext}}^{\rm s}$ denotes the extended displacement solution on the extended solid subdomain $\Omega_{h,n}^{s,\text{ext}}$ and $\Phi_0$ describes the initial configuration of the interface. For {more} details we refer to \cite{FrKnStWeWi25}.

\begin{algorithm}[h!]
    For $n \ge 1$:
    \begin{enumerate}
        \item Find $\vec{y}_h^n\eqd(\vec{v}_h^{\mathrm f,n}, p_h^n, \vec{v}_h^{s, n}, \vec{u}_h^{s, n}) \in \mathcal{X}_{h{,n}}$ such that
		{\begin{multline}\label{varf:thick2}
		\varrho^{\rm f} \big( {{\partial}_{\delta t}} \vec v_h^{\mathrm f,n},  \vec \phi_h^{\rm f} \big)_{\Omega^{\mathrm f,n}} + \tilde{\varrho}^{\rm s} \big( {{\partial}_{\delta t}} \vec v_h^{\mathrm s,n},  \vec \phi_h^{\rm s} \big)_{\Omega^{\mathrm s,n}} +  \big( {{\partial}_{\delta t}} \vec u_h^{n},  \vec z_h^{\rm s} \big)_{\Omega^{\mathrm s,n}}+ {\cal A}_\text{Eu,N}^{t_n}((\vec y_h^n; {\boldsymbol \phi}_h))  \\
		+G_h^\text{ext}(\vec{y}_h^n;{\Phi}_h^n)+G_h(\vec{y}_h^n;{\Phi}_h^n)+\frac{\gamma_{{c}} \lambda^{\rm s}}{h} \big(\big[{P}_{\gamma_c}  \big]_+,\vec \phi_h^{\rm s}  \cdot \vec n \big)_{\Gamma_n^I}
        =  \sum_{i\in \{f,s\}} \int_{\Omega^i_t} \tilde\varrho^i\vec{g}\cdot{\boldsymbol \phi}_h^i	
		\end{multline}}
		 for all $\vec\Phi_h^n\eqd({\boldsymbol \phi}_h^{\rm f}, {\boldsymbol \phi}_h^{\rm s}, q_h, \vec{z}_h^{\rm s}) \in {\cal X}_h^n$, 
        where $P_{\gamma_c}$ is defined as in Alg. \ref{alg:Fr} with $\gamma_c \lambda^{\rm s}=10^9$ and {$\partial_{\delta t}\vec{w} \eqd (\delta t)^{-1}(\vec{w}_h^n-\vec{w}_h^{n-1})$.}
        \item Update the interface and the subdomains:
        \begin{align*}
            \Gamma_n^I\eqd\{x \in \Omega \ | \ \Phi(x,t_n)=0\}, \ \Omega^{\mathrm s,n}\eqd\{x \in \Omega \ | \ \Phi(x,t_n)<0\}, \ \Omega^{\mathrm f,n}\eqd\{x \in \Omega \ | \ \Phi(x,t_n)>0\}
        \end{align*}
        and based on that the computational and extended subdomains.
    \end{enumerate}
\caption{[KWF]:  Fully Eulerian CutFEM with Nitsche coupling and relaxed contact  \label{alg:KWF}}
\end{algorithm}

\paragraph{No-penetration condition.}

In our implementation, we use a relaxed formulation of the contact conditions, see Section~\ref{subsec.contact}, but with a very small $\epsilon=10^{-8}$. Thus, in the simulations conducted for this article, the contact condition never got active and we obtained contactless rebounds.

\paragraph{Numerical solution.}
The nonlinear system is solved applying Newton's method with backtracking line search in each time step. For the linear subproblems we use a direct solver (MUMPS).

\section*{Data availability}
The time histories of the QoIs provided by all eight approaches are available as supplementary material on Mendeley Data~\cite{BenchmarkData}. 

\section*{Acknowledgements}
Stefan Frei and Thomas Wick gratefully acknowledge the financial support from the German Research Foundation (DFG) with the grant number 548064929. 
Jakub Fara, Sebastian Schwarzacher and Karel T\r{u}ma have been supported by the ERC-CZ grant LL2105 of the Ministry of Education, Youth, and Sport of the Czech Republic and by Charles University Research program No. UNCE/SCI/023.

\bibliographystyle{abbrv}
\bibliography{reference}

\end{document}